\documentclass[11pt]{article}

\usepackage[utf8]{inputenc}
\usepackage[T1]{fontenc}
\usepackage[english]{babel}

\usepackage{amsmath, amssymb, amsthm}
\usepackage{ytableau}        % Young tableaux
\usepackage{float}
\usepackage{caption}
\usepackage{hyperref}
\usepackage{cleveref}
\usepackage{mathdots}
\usepackage{lineno}

\usepackage{xcolor}
\usepackage{soul}            % \hl highlighting
\usepackage{tikz}
\usetikzlibrary{
  arrows.meta,
  arrows,
  calc,
  decorations.pathmorphing,
  decorations.pathreplacing,
  positioning,
  shapes
}
\usepackage{caption}
\usepackage{subcaption}

\usepackage{enumitem}
\usepackage{mathtools}

\usepackage[backend=biber, style=numeric,,maxnames=99]{biblatex}
\usepackage{csquotes}

\usepackage{todonotes}
\usepackage{comment}

\usepackage{geometry}
\theoremstyle{definition}
\newtheorem{theorem}{Theorem}[section]
\newtheorem*{theorem-non}{Theorem}
\newtheorem{lemma}[theorem]{Lemma}
\newtheorem{proposition}[theorem]{Proposition}

\newtheorem{definition}[theorem]{Definition}
\theoremstyle{remark}

\newtheorem{example}[theorem]{Example}
\newtheorem{remark}[theorem]{Remark}

\theoremstyle{definition}
\newtheorem*{theorem*}{Theorem}
\numberwithin{equation}{section}

\newcommand{\type}{\operatorname{type}}

\newcommand{\coloneq}{\mathrel{\mathop:}=}
\definecolor{turkoise}{RGB}{64,224,208} % light turquoise

\newcommand\qand{\quad\text{and}\quad}
\definecolor{darkblue}{rgb}{0.0,0,0.7}
\renewcommand\emph[1]{\textcolor{darkblue}{\it #1}}

\newcommand{\DiagCellg}{%
  \tikz[baseline, shift={(0,-0.5em)}]{%
    \begin{scope}[x=1.5em,y=1.5em]
      \fill[green] (0,0) rectangle (1,1);
      \draw[red, line width=1pt] (0,0) -- (1,1);
    \end{scope}
  }%
}
\newcommand{\DiagCellb}{%
  \tikz[baseline, shift={(0,-0.5em)}]{%
    \begin{scope}[x=1.5em,y=1.5em]
      \fill[blue] (0,0) rectangle (1,1);
      \draw[red, line width=1pt] (0,0) -- (1,1);
    \end{scope}
  }%
}

\title{$A_1^{(1)}$-Grounded partitions at levels $1$ and $2$,\\ Part II: affine crystal structure}
\author{Benedek Dombos\thanks{\href{mailto:benedek.dombos@unige.ch}{benedek.dombos@unige.ch}, University of Geneva} \and Jihyeug Jang\thanks{\href{mailto:jihyeugjang@gmail.com}{jihyeugjang@gmail.com}, Pusan National University}}
\date{November 12, 2026}

\makeatletter
\renewcommand{\@fnsymbol}[1]{%
  \ifcase#1
  \or *
  \or **%
  \or \dagger
  \or \ddagger
  \or \mathsection
  \or \mathparagraph
  \else
  \@ctrerr
  \fi}
\makeatother

\begin{document}

\maketitle

% Prevent clear page between includes
\let\oldclearpage\clearpage
\renewcommand{\clearpage}{}

\begin{abstract}
Grounded partitions are coloured partitions with difference conditions, introduced by Dousse and Konan, that are in bijection with the paths of the Kyoto path model of Kang, Kashiwara, Misra, Miwa, Nakashima and Nakayashiki. We describe the affine crystal structure on the grounded partitions arising from the level \(1\) and level \(2\) perfect crystals of type \(A_1^{(1)}\). At level \(1\), grounded partitions are naturally identified with \(2\)-restricted partitions, and the affine crystal structure is given by Kleshchev's modular branching rule for symmetric groups in characteristic \(2\). At level \(2\), we give explicit bracketing rules which make grounded partitions into affine crystal graphs of highest weights \(2\Lambda_0\), \(\Lambda_0+\Lambda_1\), and \(2\Lambda_1\). We prove this by constructing explicit graph isomorphisms with the asymptotic model of Ariki and Mathas, whose vertices are Kleshchev multipartitions arising from Ariki--Koike algebras. Finally, we erase the \(0\)-arrows, which leaves the crystal graphs of the restrictions to the Levi subalgebra, and decompose them into finite strings. This decomposition gives new \(q\)-series identities at levels \(1\) and \(2\).
\end{abstract}

% Intro
\section{Introduction}

In the first part of this work~\cite{DJ25}, we studied grounded partitions of type \(A_1^{(1)}\) at levels \(1\) and \(2\) from a combinatorial point of view. We now describe the affine crystal structure on these grounded partitions by explicit bracketing rules. The main construction is at level \(2\): we identify the resulting graphs with the corresponding affine crystals by constructing explicit isomorphisms with the \textit{asymptotic model}. We first recall the algebraic motivation and the definition of grounded partitions.

A \emph{partition} of a nonnegative integer \( n \) is a non-increasing sequence \(\lambda = (\lambda_1, \lambda_2, \dots, \lambda_k)\) of positive integers whose sum is \( n \). The integers \( \lambda_i \) are its \emph{parts}, and their sum is its \emph{size}, denoted by $|\lambda|$. The empty sequence is the unique partition of size $0$, called the \emph{empty partition}.

For example, there are five partitions of size 4: \((4), (3,1), (2,2), (2,1,1) \), and \( (1,1,1,1)\).
We identify partitions with their \emph{Young diagrams}:
\begin{align*}
\overset{\textcolor{purple}{(4)}}{\scalebox{0.7}{\begin{ytableau}
          \null & \null & \null & \null
        \end{ytableau}}}
\quad
\overset{\textcolor{purple}{(3,1)}}{\scalebox{0.7}{\begin{ytableau}
          \null & \null & \null \\
          \null
        \end{ytableau}}}
\quad
\overset{\textcolor{purple}{(2,2)}}{\scalebox{0.7}{\begin{ytableau}
          \null & \null \\
          \null & \null
        \end{ytableau}}}
\quad
\overset{\textcolor{purple}{(2,1,1)}}{\scalebox{0.7}{\begin{ytableau}
          \null & \null \\
          \null \\
          \null
        \end{ytableau}}}
\quad
\overset{\textcolor{purple}{(1,1,1,1)}}{\scalebox{0.7}{\begin{ytableau}
          \null \\
          \null \\
          \null \\
          \null
        \end{ytableau}}}
\end{align*}

\begin{figure}[H]
\centering
\scalebox{0.6}{%
\begin{tikzpicture}[>=stealth, thick, x=2.0cm,y=1.6cm]
  \tikzset{mynode/.style={inner sep=0pt, font=\scriptsize}}
  \tikzset{arrow style/.style={->, shorten <=2pt, shorten >=2pt, line width=0.8pt}}

  % level 0
  \node[mynode] (e) at (0,0) {\(\emptyset\)};

  % level 1
  \node[mynode] (1) at (0,-1) {\(\scalebox{0.8}{\begin{ytableau}\null\end{ytableau}}\)};
  \draw[arrow style] (e) -- (1);

  % level 2
  \node[mynode] (2)  at (-0.8,-2) {\(\scalebox{0.8}{\begin{ytableau}\null&\null\end{ytableau}}\)};
  \node[mynode] (11) at ( 0.8,-2) {\(\scalebox{0.8}{\begin{ytableau}\null\\\null\end{ytableau}}\)};
  \draw[arrow style] (1) -- (2);
  \draw[arrow style] (1) -- (11);

  % level 3
  \node[mynode] (3)   at (-1.6,-3) {\(\scalebox{0.8}{\begin{ytableau}\null&\null&\null\end{ytableau}}\)};
  \node[mynode] (21)  at ( 0.0,-3) {\(\scalebox{0.8}{\begin{ytableau}\null&\null\\\null\end{ytableau}}\)};
  \node[mynode] (111) at ( 1.6,-3) {\(\scalebox{0.8}{\begin{ytableau}\null\\\null\\\null\end{ytableau}}\)};
  \draw[arrow style] (2)  -- (3);
  \draw[arrow style] (2)  -- (21);
  \draw[arrow style] (11) -- (21);
  \draw[arrow style] (11) -- (111);

  % level 4
  \node[mynode] (4)    at (-3.2,-4) {\(\scalebox{0.8}{\begin{ytableau}\null&\null&\null&\null\end{ytableau}}\)};
  \node[mynode] (31)   at (-1.6,-4) {\(\scalebox{0.8}{\begin{ytableau}\null&\null&\null\\\null\end{ytableau}}\)};
  \node[mynode] (22)   at ( 0.0,-4) {\(\scalebox{0.8}{\begin{ytableau}\null&\null\\\null&\null\end{ytableau}}\)};
  \node[mynode] (211)  at ( 1.6,-4) {\(\scalebox{0.8}{\begin{ytableau}\null&\null\\\null\\\null\end{ytableau}}\)};
  \node[mynode] (1111) at ( 3.2,-4) {\(\scalebox{0.8}{\begin{ytableau}\null\\\null\\\null\\\null\end{ytableau}}\)};

  \draw[arrow style] (3)   -- (4);
  \draw[arrow style] (3)   -- (31);
  \draw[arrow style] (21)  -- (31);
  \draw[arrow style] (21)  -- (22);
  \draw[arrow style] (21)  -- (211);
  \draw[arrow style] (111) -- (211);
  \draw[arrow style] (111) -- (1111);

  % level 5
  \node[mynode] (5)     at (-4.8,-5.5) {\(\scalebox{0.8}{\begin{ytableau}\null&\null&\null&\null&\null\end{ytableau}}\)};
  \node[mynode] (41)    at (-3.2,-5.5) {\(\scalebox{0.8}{\begin{ytableau}\null&\null&\null&\null\\\null\end{ytableau}}\)};
  \node[mynode] (32)    at (-1.6,-5.5) {\(\scalebox{0.8}{\begin{ytableau}\null&\null&\null\\\null&\null\end{ytableau}}\)};
  \node[mynode] (311)   at ( 0.0,-5.5) {\(\scalebox{0.8}{\begin{ytableau}\null&\null&\null\\\null\\\null\end{ytableau}}\)};
  \node[mynode] (221)   at ( 1.6,-5.5) {\(\scalebox{0.8}{\begin{ytableau}\null&\null\\\null&\null\\\null\end{ytableau}}\)};
  \node[mynode] (2111)  at ( 3.2,-5.5) {\(\scalebox{0.8}{\begin{ytableau}\null&\null\\\null\\\null\\\null\end{ytableau}}\)};
  \node[mynode] (11111) at ( 4.8,-5.5) {\(\scalebox{0.8}{\begin{ytableau}\null\\\null\\\null\\\null\\\null\end{ytableau}}\)};

  \draw[arrow style] (4)    -- (5);
  \draw[arrow style] (4)    -- (41);

  \draw[arrow style] (31)   -- (41);
  \draw[arrow style] (31)   -- (32);
  \draw[arrow style] (31)   -- (311);

  \draw[arrow style] (22)   -- (32);
  \draw[arrow style] (22)   -- (221);

  \draw[arrow style] (211)  -- (311);
  \draw[arrow style] (211)  -- (221);
  \draw[arrow style] (211)  -- (2111);

  \draw[arrow style] (1111) -- (2111);
  \draw[arrow style] (1111) -- (11111);

\end{tikzpicture}%
}
\caption{The branching graph of symmetric groups over a field of characteristic $0$}
\label{fig:young-lattice}
\end{figure}
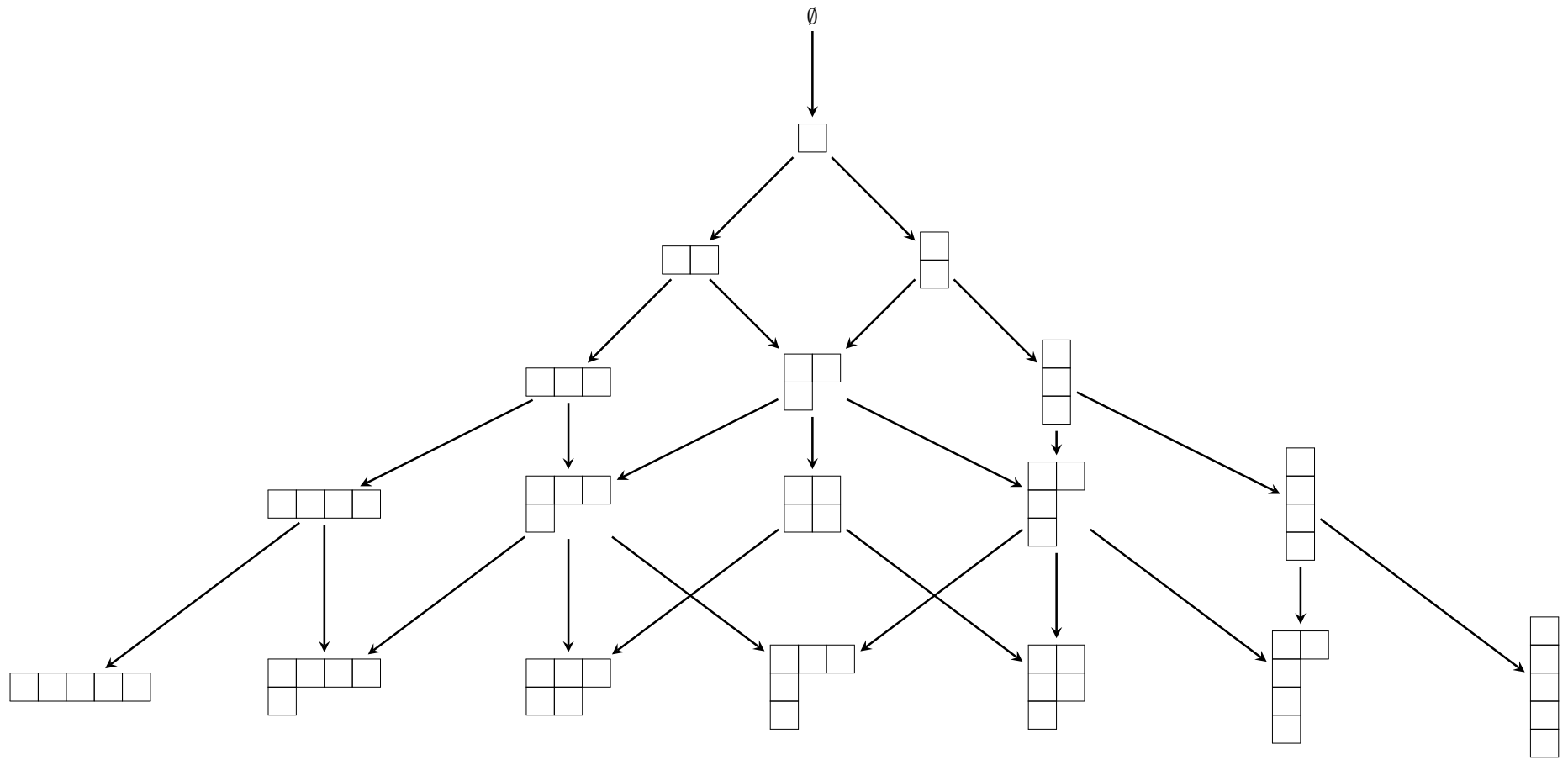

Let \(F\) be a field of characteristic \(0\), and let \(S_n\) denote the symmetric group on \(n\) letters. The irreducible \(F S_n\)-modules are indexed by partitions of \(n\). For a partition $\lambda$ of $n$, write \(S^\lambda\) for the corresponding irreducible \(F S_n\)-module. The \emph{branching graph} of the symmetric groups over \(F\) is the directed graph whose vertex set is the set of partitions, with an arrow from \(\mu\) to \(\lambda\) whenever \(S^\mu\) is an irreducible constituent of the restriction \(\operatorname{Res}^{S_n}_{S_{n-1}}S^\lambda\). Young's branching rule~\cite[Theorem~VII, p.~269]{You28} states that
\begin{align*}
\operatorname{Res}^{S_n}_{S_{n-1}} S^\lambda
\cong
\bigoplus_{\mu} S^\mu,
\end{align*}
where the sum is over all partitions $\mu$ of $n-1$ obtained from \(\lambda\) by removing one cell from the diagram. The branching graph of the symmetric groups over $F$ is shown in Figure~\ref{fig:young-lattice}.

For an integer \(e\geq2\), a partition \(\lambda=(\lambda_1,\ldots,\lambda_k)\) is called \emph{\(e\)-restricted} if \(\lambda_i-\lambda_{i+1}<e\) for \(1\leq i\leq k\), where \(\lambda_{k+1}\coloneq0\). For example, the \(2\)-restricted partitions of \(4\) are \((2,1,1)\) and \((1,1,1,1)\).

Let $p$ be a prime and let \(F\) be a field of characteristic \(p\). The irreducible \(F S_n\)-modules are labelled by \(p\)-restricted partitions of \(n\); write \(D^\lambda\) for the module labelled by \(\lambda\). Modules need not be direct sums of irreducibles in this setting, so the branching graph is defined using the \emph{socle}, the sum of all irreducible submodules. Its vertices are the \(p\)-restricted partitions, with an arrow \(\mu\to\lambda\) when the socle of \(\operatorname{Res}^{S_n}_{S_{n-1}}D^\lambda\) contains a submodule isomorphic to \(D^\mu\). Kleshchev~\cite[Theorem~0.5, p.~166]{Kl95} gave a combinatorial description of this graph, stated there for $p$-regular partitions, the conjugates of $p$-restricted partitions. We recall it for \(p=2\) in Definitions~\ref{def:partition_corners} and~\ref{def:partition_model_1}.

\begin{remark}
The notation \( \mathcal{B}_{1,b}^{\mathcal{P}} \) that we use for this graph in Definition~\ref{def:partition_model_1} will be justified below: it is chosen to match the notation of grounded partitions, which are the focus of this paper.
\end{remark}

\begin{definition}\label{def:partition_corners}
Given a partition \(\lambda\), we consider the blue--green
checkerboard pattern on \(\lambda\), with a green cell in the top-left
corner. A \emph{removable blue corner} is a blue cell inside
\( \lambda \) such that removing it still gives a partition. An
\emph{addable blue corner} is a blue cell outside \( \lambda \) such
that adding it gives a partition. \emph{Removable green corners} and
\emph{addable green corners} are defined in the same way.
\end{definition}

\begin{definition}\label{def:partition_model_1}
  The directed graph \( \mathcal{B}_{1,b}^{\mathcal{P}} \), whose
  vertices are partitions and whose edges are coloured blue and green,
  is defined recursively starting from the empty partition. For each
  vertex \( \lambda \), the outgoing edges are determined by the
  following bracketing rules:
  \begin{enumerate}
  \item \textbf{Blue arrow.} Assign a left bracket ``$($'' to each
    addable blue corner of \( \lambda \), and a right bracket ``$)$''
    to each removable blue corner. Read brackets from bottom to top.
    If there exists an unmatched addable blue corner, add a cell at the corner corresponding to the leftmost unmatched left bracket. A
    blue arrow is drawn from \( \lambda \) to the resulting partition.
  \item \textbf{Green arrow.} Defined in the same way, with addable and
    removable green corners in place of blue ones.
  \end{enumerate}
\end{definition}

\begin{example}\label{ex:addable_corner_bracketing}
We show a few examples of the bracketing rules from Definition~\ref{def:partition_model_1}. The first and third partitions have no unmatched addable blue corners, whereas the second has two:
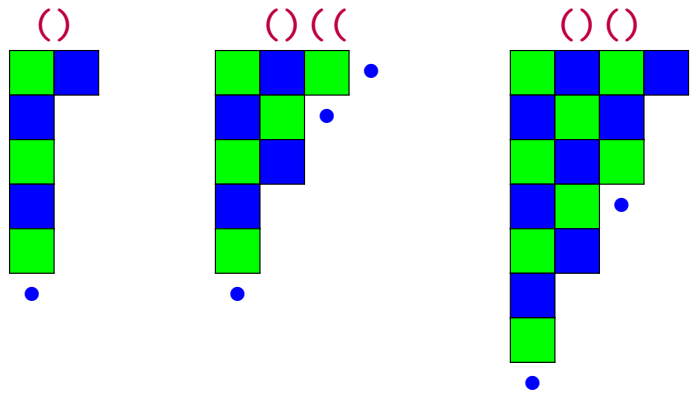
\begin{align*}
\overset{\textcolor{purple}{\texttt{()}}}{\scalebox{0.5}{\begin{ytableau}
  *(green)&*(blue)\\
  *(blue)\\
  *(green)\\
  *(blue)\\
  *(green)\\
  \none[{\color{blue}\scalebox{1.25}{$\bullet$}}]
\end{ytableau}}}
\qquad\overset{\textcolor{purple}{\texttt{()((}}}{\scalebox{0.5}{\begin{ytableau}
  *(green)&*(blue)&*(green)&\none[{\color{blue}\scalebox{1.25}{$\bullet$}}]\\
  *(blue)&*(green)&\none[{\color{blue}\scalebox{1.25}{$\bullet$}}]\\
  *(green)&*(blue)\\
  *(blue)\\
  *(green)\\
  \none[{\color{blue}\scalebox{1.25}{$\bullet$}}]
\end{ytableau}}}
\qquad
\overset{\textcolor{purple}{\texttt{()()}}}{\scalebox{0.5}{\begin{ytableau}
  *(green)&*(blue)&*(green)&*(blue)\\
  *(blue)&*(green)&*(blue)\\
  *(green)&*(blue)&*(green)\\
  *(blue)&*(green)&\none[{\color{blue}\scalebox{1.25}{$\bullet$}}]\\
  *(green)&*(blue)\\
  *(blue)\\
  *(green)\\
  \none[{\color{blue}\scalebox{1.25}{$\bullet$}}]
\end{ytableau}}}
\end{align*}
\end{example}

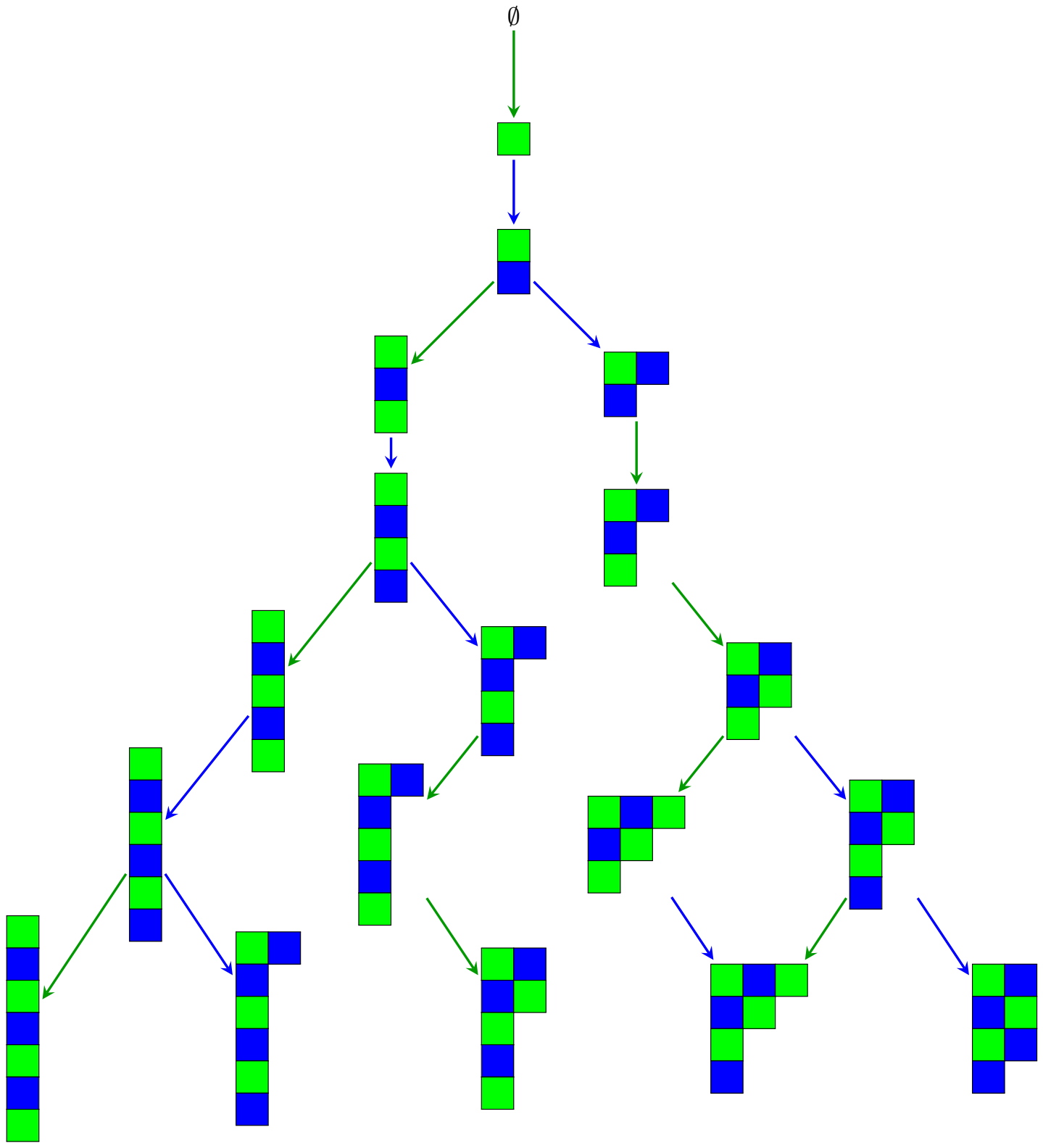
\begin{figure}[h]
\centering
\scalebox{0.8}{\begin{tikzpicture}[>=stealth, thick, scale=0.7]
  %%%%%%%%%%%%%%
  % Coordinates (unchanged)
  %%%%%%%%%%%%%%
  % Row 1: 1 vertex
  \coordinate (R1C1) at (0,0);
  
  % Row 2: 1 vertex
  \coordinate (R2C1) at (0,-2);
  
  % Row 3: 1 vertex
  \coordinate (R3C1) at (0,-4);
  
  % Row 4: 2 vertices
  \coordinate (R4C1) at (-2,-6);
  \coordinate (R4C2) at (2,-6);
  
  % Row 5: 2 vertices
  \coordinate (R5C1) at (-2,-8.5);
  \coordinate (R5C2) at (2,-8.5);
  
  % Row 6: 3 vertices
  \coordinate (R6C1) at (-4,-11);
  \coordinate (R6C2) at (0,-11);
  \coordinate (R6C3) at (4,-11);
  
  % Row 7: 4 vertices
  \coordinate (R7C1) at (-6,-13.5);
  \coordinate (R7C2) at (-2,-13.5);
  \coordinate (R7C3) at (2,-13.5);
  \coordinate (R7C4) at (6,-13.5);
    
  % Row 8: 5 vertices
  \coordinate (R8C1) at (-8,-16.5);
  \coordinate (R8C2) at (-4,-16.5);
  \coordinate (R8C3) at (0,-16.5);
  \coordinate (R8C4) at (4,-16.5);
  \coordinate (R8C5) at (8,-16.5);
  
  %%%%%%%%%%%%%%
  % Nodes (transposed shapes)
  %%%%%%%%%%%%%%
  \tikzset{mynode/.style={inner sep=0pt, font=\scriptsize}}
  
  \node[mynode] (N1C1) at (R1C1) {\(\emptyset\)};
  
  \node[mynode] (N2C1) at (R2C1) {\(\scalebox{0.8}{\begin{ytableau}
*(green)
\end{ytableau}}\)};
  
  \node[mynode] (N3C1) at (R3C1) {\(\scalebox{0.8}{\begin{ytableau}
*(green)\\
*(blue)
\end{ytableau}}\)};
  
  \node[mynode] (N4C1) at (R4C1) {\(\scalebox{0.8}{\begin{ytableau}
*(green)\\
*(blue)\\
*(green)
\end{ytableau}}\)};
  \node[mynode] (N4C2) at (R4C2) {\(\scalebox{0.8}{\begin{ytableau}
*(green) & *(blue)\\
*(blue)
\end{ytableau}}\)};
  
  \node[mynode] (N5C1) at (R5C1) {\(\scalebox{0.8}{\begin{ytableau}
*(green)\\
*(blue)\\
*(green)\\
*(blue)
\end{ytableau}}\)};
  \node[mynode] (N5C2) at (R5C2) {\(\scalebox{0.8}{\begin{ytableau}
*(green) & *(blue)\\
*(blue)\\
*(green)
\end{ytableau}}\)};
  
  \node[mynode] (N6C1) at (R6C1) {\(\scalebox{0.8}{\begin{ytableau}
*(green)\\
*(blue)\\
*(green)\\
*(blue)\\
*(green)
\end{ytableau}}\)};
  \node[mynode] (N6C2) at (R6C2) {\(\scalebox{0.8}{\begin{ytableau}
*(green) & *(blue)\\
*(blue)\\
*(green)\\
*(blue)
\end{ytableau}}\)};
  \node[mynode] (N6C3) at (R6C3) {\(\scalebox{0.8}{\begin{ytableau}
*(green) & *(blue)\\
*(blue) & *(green)\\
*(green)
\end{ytableau}}\)};
  
  \node[mynode] (N7C1) at (R7C1) {\(\scalebox{0.8}{\begin{ytableau}
*(green)\\
*(blue)\\
*(green)\\
*(blue)\\
*(green)\\
*(blue)
\end{ytableau}}\)};
  \node[mynode] (N7C2) at (R7C2) {\(\scalebox{0.8}{\begin{ytableau}
*(green) & *(blue)\\
*(blue)\\
*(green)\\
*(blue)\\
*(green)
\end{ytableau}}\)};
  \node[mynode] (N7C3) at (R7C4) {\(\scalebox{0.8}{\begin{ytableau}
*(green) & *(blue)\\
*(blue) & *(green)\\
*(green)\\
*(blue)
\end{ytableau}}\)};
  \node[mynode] (N7C4) at (R7C3) {\(\scalebox{0.8}{\begin{ytableau}
*(green) & *(blue) & *(green)\\
*(blue) & *(green)\\
*(green)
\end{ytableau}}\)};

  \node[mynode] (N8C1) at (R8C1) {\(\scalebox{0.8}{\begin{ytableau}
*(green)\\
*(blue)\\
*(green)\\
*(blue)\\
*(green)\\
*(blue)\\
*(green)
\end{ytableau}}\)};
  \node[mynode] (N8C2) at (R8C2) {\(\scalebox{0.8}{\begin{ytableau}
*(green) & *(blue)\\
*(blue)\\
*(green)\\
*(blue)\\
*(green)\\
*(blue)
\end{ytableau}}\)};
  \node[mynode] (N8C3) at (R8C3) {\(\scalebox{0.8}{\begin{ytableau}
*(green) & *(blue)\\
*(blue) & *(green)\\
*(green)\\
*(blue)\\
*(green)
\end{ytableau}}\)};
  \node[mynode] (N8C4) at (R8C5) {\(\scalebox{0.8}{\begin{ytableau}
*(green) & *(blue)\\
*(blue) & *(green)\\
*(green) & *(blue)\\
*(blue)
\end{ytableau}}\)};
  \node[mynode] (N8C5) at (R8C4) {\(\scalebox{0.8}{\begin{ytableau}
*(green) & *(blue) & *(green)\\
*(blue) & *(green)\\
*(green)\\
*(blue)
\end{ytableau}}\)};
  
  %%%%%%%%%%%%%%
  % Arrows (unchanged)
  %%%%%%%%%%%%%%
  \tikzset{arrow style/.style={->, shorten <=1pt, shorten >=1pt, line width=0.8pt}}
  \draw[arrow style, green!60!black] (N1C1) -- (N2C1);
  \draw[arrow style, blue]                  (N2C1) -- (N3C1);
  \draw[arrow style, green!60!black] (N3C1) -- (N4C1);
  \draw[arrow style, blue]                  (N3C1) -- (N4C2);
  \draw[arrow style, blue]                  (N4C1) -- (N5C1);
  \draw[arrow style, green!60!black] (N4C2) -- (N5C2);
  \draw[arrow style, green!60!black] (N5C1) -- (N6C1);
  \draw[arrow style, blue]                  (N5C1) -- (N6C2);
  \draw[arrow style, green!60!black] (N5C2) -- (N6C3);
  \draw[arrow style, blue]                  (N6C1) -- (N7C1);
  \draw[arrow style, green!60!black] (N6C2) -- (N7C2);
  \draw[arrow style, blue]                  (N6C3) -- (N7C3);
  \draw[arrow style, green!60!black] (N6C3) -- (N7C4);
  \draw[arrow style, green!60!black] (N7C1) -- (N8C1);
  \draw[arrow style, blue]                  (N7C1) -- (N8C2);
  \draw[arrow style, green!60!black] (N7C2) -- (N8C3);
  \draw[arrow style, blue]                  (N7C3) -- (N8C4);
  \draw[arrow style, green!60!black] (N7C3) -- (N8C5);
  \draw[arrow style, blue]                  (N7C4) -- (N8C5);
\end{tikzpicture}}
\caption{The directed graph \( \mathcal{B}_{1,b}^{\mathcal{P}} \)}
\label{fig:checkerboard}
\end{figure}

The graph \( \mathcal{B}_{1,b}^{\mathcal{P}} \) is shown in Figure~\ref{fig:checkerboard}. One can verify that the vertex set of the graph
\( \mathcal{B}_{1,b}^{\mathcal{P}} \) is given by the $2$-restricted partitions. 

The same combinatorial rule makes sense with \(e\) colours for any integer \(e\geq2\): colour the cell in row \(i\) and column \(j\) by \(j-i\) modulo \(e\). Brundan~\cite[Theorem~2.6, p.~554]{Br98} identified the resulting graph with the branching graph for irreducible modules of the Hecke algebra of type \(A\) with quantum characteristic \(e\).

Ariki--Koike algebras, introduced by Ariki and Koike~\cite[Definition~3.1, p.~217]{AK94}, are cyclotomic Hecke algebras deforming the group algebra of the complex reflection group $G(r,1,n)=(\mathbb{Z}/r\mathbb{Z})\wr S_n$. Ariki~\cite[Theorem~4.4, p.~806]{Ari96} identified the Grothendieck groups for Ariki--Koike algebras with integrable highest-weight modules of affine type $A_{e-1}^{(1)}$ at level $r$, where $e$ is the order of the Hecke parameter $q$, with induction and restriction corresponding to the Chevalley operators. Following Jimbo, Misra, Miwa and Okado~\cite[Theorems~3.8 and~3.13, pp.~554 and~556]{JMMO91}, Foda, Leclerc, Okado, Thibon and Welsh~\cite[Theorem~2.10 and Proposition~2.11]{FLOTW99} gave a combinatorial model for affine crystal graphs of type $A_{e-1}^{(1)}$ at arbitrary level, known as the \textit{FLOTW model}. Ariki and Mathas~\cite[Theorem~2.8 and Corollary~2.10, pp.~620--621]{AM00} described another combinatorial model for affine crystal graphs of type $A_{e-1}^{(1)}$ at arbitrary level, which we call the \emph{asymptotic model}. Its vertices are certain multipartitions, called \emph{Kleshchev multipartitions}~\cite[p.~605]{AM00}.

We use the level-$2$ asymptotic model of type $A_1^{(1)}$, defined in Section~\ref{section:partition_model}. Its vertices are certain pairs of $2$-restricted partitions, each with a blue--green checkerboard colouring. Ariki, Kreiman and Tsuchioka~\cite[Theorem~9.5]{AKT08} characterised Kleshchev bipartitions in crystal-theoretic terms. Chern, Li, Stanton, Xue and Yee~\cite[Theorem~1.2 and Corollary~1.3]{CLSXY24} studied Kleshchev multipartitions from a combinatorial point of view.

Crystal graphs, introduced by Kashiwara~\cite[\S4, p.~252, and \S6, p.~255]{Kas90} (see also~\cite{Kas91}), are edge-coloured directed graphs associated with representations of quantum groups, with colours indexed by the set \(I\) of simple roots. A representation \(V\) need not have a basis on which every Chevalley generator \(e_i\) and \(f_i\), \(i\in I\), maps basis elements to basis elements or to \(0\). This is the role of crystal bases: one chooses a lattice \(L\subset V\) over the ring of rational functions regular at \(q=0\), together with a basis \(B\) of \(L/qL\), such that the Kashiwara operators \(\widetilde e_i\) and \(\widetilde f_i\) map \(B\) into \(B\cup\{0\}\). The associated crystal graph has vertex set \(B\) and an \(i\)-arrow \(b\to b'\) whenever \(\widetilde f_i b=b'\). A `combinatorial model' is an edge-coloured directed graph on combinatorial objects, with arrows given by explicit rules, that is isomorphic to this crystal graph.

In type $A$, there is a combinatorial model for the crystals of highest weight irreducible representations whose vertices are semistandard Young tableaux and whose arrows are given by a simple bracketing rule. The bracketing rule is due to Greene and Kleitman~\cite[Theorems~1 and~7, pp.~82 and~86]{GK76}, in a different context. This model goes back to Lascoux--Schützenberger’s plactic theory~\cite[2.9 and 4.3, pp.~133--134 and~144]{LS81}; its crystal interpretation was described by Kashiwara and Nakashima~\cite[Theorem~3.4.2, p.~309]{KN94}. In affine type $A_{e-1}^{(1)}$, Misra and Miwa gave a partition model for the crystal graph of the basic representation~\cite[Propositions~3.3 and~4.3 and Theorem~4.7, pp.~85, 86 and~88]{MM90}. Lascoux, Leclerc and Thibon~\cite[Proposition~5.2, p.~220]{LLT96} observed that, for \(e=p\), the Misra--Miwa partition model coincides with Kleshchev's modular branching graph for symmetric groups in characteristic \(p\), defined in Definitions~\ref{def:partition_corners} and~\ref{def:partition_model_1} for $p=2$.

Perfect crystals are finite crystals from which affine highest weight crystals can be constructed recursively. Kang, Kashiwara, Misra, Miwa, Nakashima and Nakayashiki introduced them in~\cite{KMN2a} (see also~\cite[Definition~1.1.1]{KMN92}). Given a perfect crystal, they constructed the associated crystal of semi-infinite paths and proved that it is isomorphic, as an edge-coloured directed graph, to the corresponding affine highest weight crystal~\cite[Theorem~4.5.2]{KMN2a} (see also~\cite[Proposition~1.1.2]{KMN92} and Hong and Kang~\cite[Theorem~10.6.4, p.~254]{HK02}). This is the \emph{Kyoto path model}. Grounded partitions, introduced by Dousse and Konan~\cite[Definition~3.2 and Proposition~3.4, p.~15]{DK19b}, are coloured partitions with difference conditions, in bijection with the paths of this model.

Dousse, Hardiman and Konan~\cite{DHK25} defined grounded partitions of type \( A_1^{(1)} \) at all levels. We consider only grounded partitions of this type and henceforth refer to them simply as grounded partitions.

\begin{definition}\label{def:grounded_partitions}
  Let \( \ell \) be a nonnegative integer, and let \(C_\ell=\{c_0,c_1,\ldots,c_\ell\}\) be a set of \emph{labels}.
  Define \(M_\ell\) to be the square matrix of size \(\ell+1\), with
  rows and columns indexed by \( 0,1, \dots, \ell \), whose
  \((i,j)\)-entry is \(|\ell-i-j|\).
  Let \( c \) be a label in \( C_\ell \). A \emph{grounded partition}
  at level \(\ell\) with \emph{ground} \(0_{c}\) is a sequence
  \[
    \pi = (\pi_0, \pi_1, \dots, \pi_s),
  \]
  of integers, each carrying a label in \( C_\ell \), with the initial
  part \( \pi_0 = 0_c \) and \( \pi_1, \dots, \pi_s \) positive, satisfying the exact difference condition:
  for \( 0 \leq j < s \), if two consecutive parts \(\pi_{j} \) and
  \( \pi_{j+1}\) have sizes \( |\pi_j| \) and \( |\pi_{j+1}| \) with
  labels \(c_{i_j}\) and \(c_{i_{j+1}}\), respectively, then
  \[
    |\pi_{j+1}|-|\pi_{j}| = \text{the \( (i_{j+1},i_j) \)-entry of } M_\ell.
  \]

  We denote by \(\mathcal{P}_{\ell,i}\) the set of all grounded
  partitions at level \( \ell \) with ground \(0_{c_i}\). For a
  grounded partition \( \pi \), the number \( s \) of parts other than the ground is called the \emph{length} of \( \pi \), denoted by \( \ell(\pi) \). The \emph{size} of \(\pi\) is defined as the size of the underlying
  partition obtained by forgetting the labels, denoted \(|\pi|\). For
  convenience, we omit the ground \(\pi_0\), and write
  \( \pi = (\pi_0, \pi_1, \pi_2, \dots, \pi_s) \) simply as
  \(\pi_1 \pi_2 \cdots\pi_s\).
\end{definition}

The generating function for \( \mathcal{P}_{\ell,i} \) is the
principally specialised character of the irreducible highest-weight
module \( L(\lambda) \) of type \( A_1^{(1)} \) with highest weight
\(\lambda=i\Lambda_0+(\ell-i)\Lambda_1 \), and therefore has a nice
infinite product expression given by Lepowsky’s formula~\cite{Lep78}; see Lepowsky~\cite[(24c), p.~144]{Lepowskyproduct}.

\begin{remark}\label{rem:2}
  Partitions are typically written in weakly decreasing order. For
  grounded partitions, however, the sizes of the parts are weakly
  increasing, since the label of the initial \( 0 \)-part is
  significant. Equivalently, one could define grounded partitions in
  weakly decreasing order by reversing the difference condition and
  setting \( \pi_s = 0_c \) in Definition~\ref{def:grounded_partitions}.
\end{remark}

\begin{remark}\label{rem:3}
  For type \( A_1^{(1)} \), the level of a highest weight
  \( a \Lambda_0 + b \Lambda_1 \) is defined to be \( a+b \). Since
  the highest weight \(i\Lambda_0+(\ell-i)\Lambda_1 \) has level
  \( \ell \), we refer to \( \mathcal{P}_{\ell,i} \) as the set of
  grounded partitions at level \( \ell \).
\end{remark}

It is convenient to abbreviate some finite and infinite products by \emph{$q$-Pochhammer symbols}. For formal variables $a$ and $q$, define
\begin{align*}
    (a;q)_k &\coloneq (1 - a)(1 - aq)(1 - aq^2) \cdots (1 - aq^{k-1}), \quad\text{ for \( k \geq 1 \)},\\
     (a;q)_0 &\coloneq 1, \quad (a;q)_\infty \coloneq \prod_{k=0}^{\infty} (1 - aq^k).
\end{align*}

\begin{theorem}[{\cite[Theorem~1.6]{DHK25}}]\label{thm:gf_grounded}
  Let \( \ell \) and \( i \) be integers with \( 0 \leq i \leq \ell \).
  Then, the generating function for grounded partitions in
  \( \mathcal{P}_{\ell,i} \) is given by
\[
  \sum_{\pi \in \mathcal{P}_{\ell, i}} q^{|\pi|} = \frac{(q^{i+1};q^{\ell+2})_\infty(q^{\ell-i+1};q^{\ell+2})_\infty(q^{\ell+2};q^{\ell+2})_\infty}{(q;q^2)_\infty (q;q)_\infty}.
\]
\end{theorem}

For clarity of notation, we denote the labels by
\( (c_0, c_1, c_2) = (a,b,c) \) throughout this paper, and accordingly
write
\( \mathcal{P}_{\ell,a}, \mathcal{P}_{\ell,b}, \mathcal{P}_{\ell,c} \)
for
\( \mathcal{P}_{\ell,0}, \mathcal{P}_{\ell,1}, \mathcal{P}_{\ell,2}
\). We consider only the following matrices for $M_\ell$:
\begin{align}
\label{matrix:M}
M_1=\begin{array}{c|cc}
    & \textcolor{black}{a} & \textcolor{black}{b} \\ \hline
\textcolor{black}{a} & 1 & 0 \\
\textcolor{black}{b} & 0 & 1 \\
\end{array}
&&
M_2=\begin{array}{c|ccc}
    & \textcolor{black}{a} & \textcolor{black}{b} & \textcolor{black}{c} \\ \hline
\textcolor{black}{a} & 2 & 1 & 0 \\
\textcolor{black}{b} & 1 & 0 & 1 \\
\textcolor{black}{c} & 0 & 1 & 2 \\
\end{array}
\end{align}
The ground is always clear from context, as we specify the set of grounded partitions in advance. We write \( 0_b \in \mathcal{P}_{2,b} \) for the empty partition, using the subscript \( b \) to indicate the ground. This is the unique grounded partition of size zero.

\vspace{1em}
\begin{example}\label{ex:grounded_partitions}
\begin{enumerate}[
  label=(\roman*),
  leftmargin=0pt,
  labelwidth=1.4em,   % just wide enough for (i)
  labelsep=0.3em,     % half the usual space after (i)
  itemindent=!,
  align=left,
  itemsep=0pt
]
\item We list the elements of $\mathcal{P}_{1,b}$ of size $7$. Here the ground is $0_b$.
  \begin{align*}
    \textcolor{black}{1_b} \, \textcolor{black}{1_a} \, \textcolor{black}{1_b} \, \textcolor{black}{1_a} \, \textcolor{black}{1_b} \, \textcolor{black}{1_a} \, \textcolor{black}{1_b},\quad
    \textcolor{black}{1_b} \, \textcolor{black}{1_a} \, \textcolor{black}{1_b} \, \textcolor{black}{1_a} \, \textcolor{black}{1_b} \, \textcolor{black}{2_b},\quad
    \textcolor{black}{1_b} \, \textcolor{black}{1_a} \, \textcolor{black}{1_b} \, \textcolor{black}{2_b} \, \textcolor{black}{2_a},\quad \textcolor{black}{1_b} \, \textcolor{black}{2_b} \, \textcolor{black}{2_a} \, \textcolor{black}{2_b},\quad
    \textcolor{black}{1_b} \, \textcolor{black}{1_a} \, \textcolor{black}{2_a} \, \textcolor{black}{3_a}
  \end{align*}
  In this case, the labels are completely determined by the sizes of
  the parts, and by forgetting the labels, we recover the $2$-restricted partitions.
\item There is a simple size-preserving bijection between $\mathcal{P}_{1,b}$ and $\mathcal{P}_{1,a}$, by swapping the labels $a$ and $b$.
\item We list the elements of $\mathcal{P}_{2,b}$ of size \( 5 \):
\begin{align*}
\textcolor{black}{1_a} \textcolor{black}{1_c} \textcolor{black}{1_a} \textcolor{black}{1_c} \textcolor{black}{1_a},\ \textcolor{black}{1_c} \textcolor{black}{1_a} \textcolor{black}{1_c} \textcolor{black}{1_a} \textcolor{black}{1_c},\ \textcolor{black}{1_a} \textcolor{black}{1_c} \textcolor{black}{1_a} \textcolor{black}{2_b},\ \textcolor{black}{1_c} \textcolor{black}{1_a} \textcolor{black}{1_c} \textcolor{black}{2_b},\ \textcolor{black}{1_a} \textcolor{black}{2_b} \textcolor{black}{2_b},\ \textcolor{black}{1_c} \textcolor{black}{2_b} \textcolor{black}{2_b},\ \textcolor{black}{1_a} \textcolor{black}{1_c} \textcolor{black}{3_c},\ \textcolor{black}{1_c} \textcolor{black}{1_a} \textcolor{black}{3_a}.
\end{align*}
\item We list the elements of $\mathcal{P}_{2,a}$ of size \( 5 \):
 \begin{align*}
\textcolor{black}{1_b} \textcolor{black}{2_a} \textcolor{black}{2_c},\ \textcolor{black}{2_a} \textcolor{black}{3_b},\ \textcolor{black}{1_b} \textcolor{black}{2_c} \textcolor{black}{2_a},\ \textcolor{black}{1_b} \textcolor{black}{1_b} \textcolor{black}{1_b} \textcolor{black}{2_a},\ \textcolor{black}{1_b} \textcolor{black}{1_b} \textcolor{black}{1_b} \textcolor{black}{2_c},\ \textcolor{black}{1_b} \textcolor{black}{1_b} \textcolor{black}{1_b} \textcolor{black}{1_b} \textcolor{black}{1_b}.
\end{align*}
\end{enumerate}
\end{example}

Forgetting the labels identifies \( \mathcal{P}_{1,b} \) with the
$2$-restricted partitions. Thus the graph
$\mathcal{B}_{1,b}^{\mathcal{P}}$ in
Definition~\ref{def:partition_model_1} gives the affine crystal
structure on \( \mathcal{P}_{1,b} \); see
Figure~\ref{fig:checkerboard}. However, at level $2$, an explicit
affine crystal structure on \( \mathcal{P}_{2,b} \) was previously
unknown.

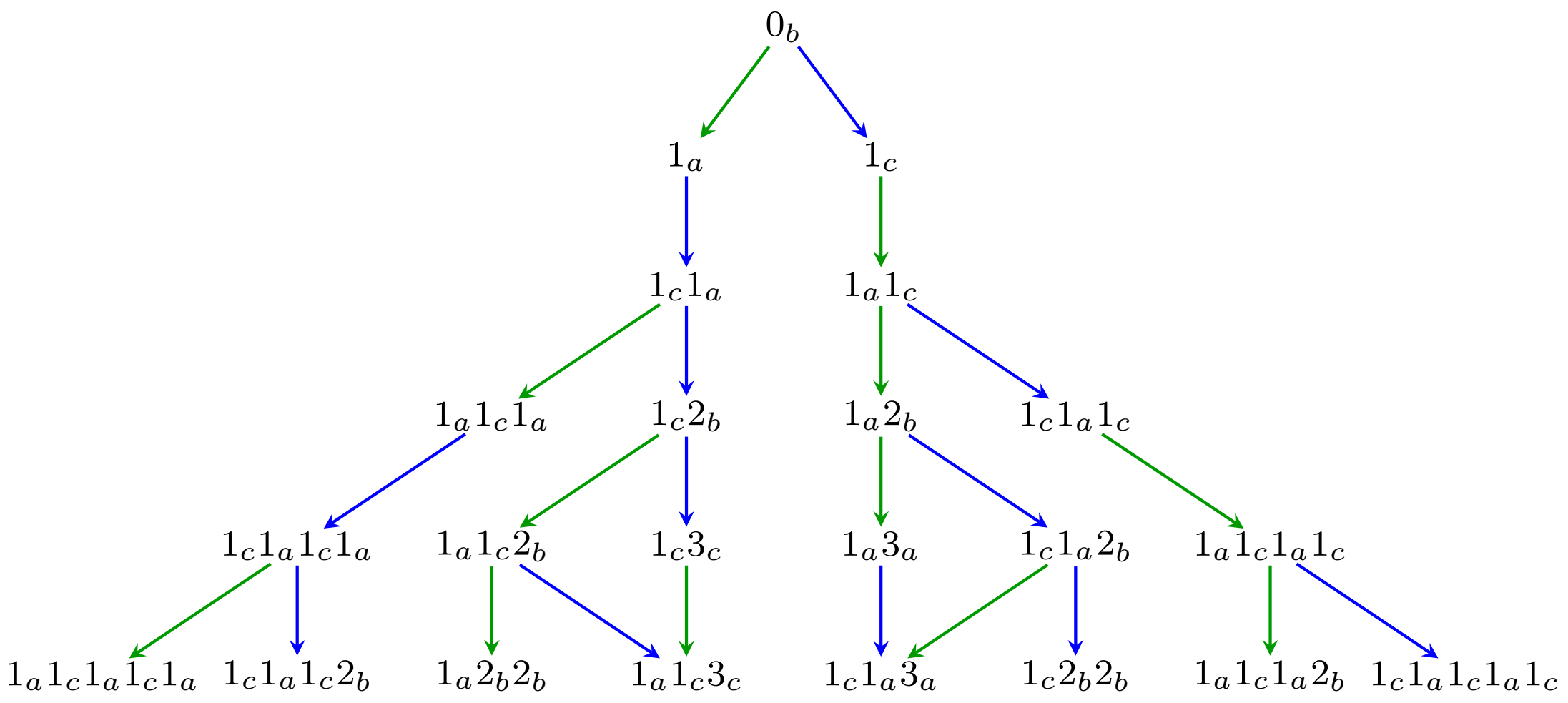
\begin{figure}
\centering
{\color{black}\begin{tikzpicture}[>=stealth, thick, scale=0.6]
  %%%%%%%%%%%%%%
  % Coordinates (unchanged)
  %%%%%%%%%%%%%%
  % Row 1: 1 vertex
  \coordinate (R1C1) at (0,0);
  
  % Row 2: 1 vertex
  \coordinate (R2C1) at (-1.5,-2);
  \coordinate (R2C2) at (1.5,-2);
  
  % Row 3: 1 vertex
  \coordinate (R3C1) at (-1.5,-4);
  \coordinate (R3C2) at (1.5,-4);
  
  % Row 4: 2 vertices
  \coordinate (R4C1) at (-4.5,-6);
  \coordinate (R4C2) at (-1.5,-6);
  \coordinate (R4C3) at (1.5,-6);
  \coordinate (R4C4) at (4.5,-6);
  
  % Row 5: 2 vertices
  \coordinate (R5C1) at (-7.5,-8);
  \coordinate (R5C2) at (-4.5,-8);
  \coordinate (R5C3) at (-1.5,-8);
  \coordinate (R5C4) at (1.5,-8);
  \coordinate (R5C5) at (4.5,-8);
  \coordinate (R5C6) at (7.5,-8);

    % Row 6: 2 vertices
  \coordinate (R6C1) at (-10.5,-10);
  \coordinate (R6C2) at (-7.5,-10);
  \coordinate (R6C3) at (-4.5,-10);
  \coordinate (R6C4) at (-1.5,-10);
  \coordinate (R6C5) at (1.5,-10);
  \coordinate (R6C6) at (4.5,-10);
  \coordinate (R6C7) at (7.5,-10);
  \coordinate (R6C8) at (10.5,-10);
  
  %%%%%%%%%%%%%%
  % Nodes (transposed shapes)
  %%%%%%%%%%%%%%
  \tikzset{mynode/.style={inner sep=0pt, font=\scriptsize}}
  
  \node[mynode] (N1C1) at (R1C1) {\(\scalebox{1.2}{$0_b$}\)};
  
  \node[mynode] (N2C1) at (R2C1) {\(\scalebox{1.2}{$1_a$}\)};

  \node[mynode] (N2C2) at (R2C2) {\(\scalebox{1.2}{$1_c$}\)};
  
  \node[mynode] (N3C1) at (R3C1) {\(\scalebox{1.2}{$1_c 1_a$}\)};

  \node[mynode] (N3C2) at (R3C2) {\(\scalebox{1.2}{$1_a 1_c$}\)};
  
  \node[mynode] (N4C1) at (R4C1) {\(\scalebox{1.2}{$1_a 1_c 1_a$}\)};

  \node[mynode] (N4C2) at (R4C2) {\(\scalebox{1.2}{$1_c 2_b$}\)};

  \node[mynode] (N4C3) at (R4C3) {\(\scalebox{1.2}{$1_a 2_b$}\)};

  \node[mynode] (N4C4) at (R4C4) {\(\scalebox{1.2}{$1_c 1_a 1_c$}\)};
  
  \node[mynode] (N5C1) at (R5C1) {\(\scalebox{1.2}{$1_c 1_a 1_c 1_a$}\)};

  \node[mynode] (N5C3) at (R5C3) {\(\scalebox{1.2}{$1_c 3_c$}\)};

  \node[mynode] (N5C2) at (R5C2) {\(\scalebox{1.2}{$1_a 1_c 2_b$}\)};

  \node[mynode] (N5C5) at (R5C5) {\(\scalebox{1.2}{$1_c 1_a 2_b$}\)};

  \node[mynode] (N5C4) at (R5C4) {\(\scalebox{1.2}{$1_a 3_a$}\)};

  \node[mynode] (N5C6) at (R5C6) {\(\scalebox{1.2}{$1_a 1_c 1_a 1_c$}\)};

  \node[mynode] (N6C1) at (R6C1) {\(\scalebox{1.2}{$1_a 1_c 1_a 1_c 1_a$}\)};

  \node[mynode] (N6C2) at (R6C2) {\(\scalebox{1.2}{$1_c 1_a 1_c 2_b$}\)};

  \node[mynode] (N6C4) at (R6C4) {\(\scalebox{1.2}{$1_a 1_c 3_c$}\)};

  \node[mynode] (N6C3) at (R6C3) {\(\scalebox{1.2}{$1_a 2_b 2_b$}\)};

  \node[mynode] (N6C6) at (R6C6) {\(\scalebox{1.2}{$1_c 2_b 2_b$}\)};

  \node[mynode] (N6C5) at (R6C5) {\(\scalebox{1.2}{$1_c 1_a 3_a$}\)};

  \node[mynode] (N6C7) at (R6C7) {\(\scalebox{1.2}{$1_a 1_c 1_a 2_b$}\)};

  \node[mynode] (N6C8) at (R6C8) {\(\scalebox{1.2}{$1_c 1_a 1_c 1_a 1_c$}\)};
  
  %%%%%%%%%%%%%%
  % Arrows (unchanged)
  %%%%%%%%%%%%%%
  \tikzset{arrow style/.style={->, shorten <=1pt, shorten >=1pt, line width=0.8pt}}
  \draw[arrow style, blue] (N1C1) -- (N2C2);
  \draw[arrow style, green!60!black] (N1C1) -- (N2C1);
  
  \draw[arrow style, blue] (N2C1) -- (N3C1);
\draw[arrow style, green!60!black] (N2C2) -- (N3C2);
  
  \draw[arrow style, blue] (N3C1) -- (N4C2);
  \draw[arrow style, green!60!black] (N3C1) -- (N4C1);
  \draw[arrow style, blue] (N3C2) -- (N4C4);
  \draw[arrow style, green!60!black] (N3C2) -- (N4C3);
  
  \draw[arrow style, blue] (N4C1) -- (N5C1);
  \draw[arrow style, blue] (N4C2) -- (N5C3);
  \draw[arrow style, green!60!black] (N4C2) -- (N5C2);
  \draw[arrow style, blue] (N4C3) -- (N5C5);
\draw[arrow style, green!60!black] (N4C3) -- (N5C4);
\draw[arrow style, green!60!black] (N4C4) -- (N5C6);

  \draw[arrow style, green!60!black] (N5C1) -- (N6C1);
  \draw[arrow style, blue] (N5C1) -- (N6C2);
  \draw[arrow style, green!60!black] (N5C2) -- (N6C3);
  \draw[arrow style, blue] (N5C2) -- (N6C4);
  \draw[arrow style, green!60!black] (N5C3) -- (N6C4);
  \draw[arrow style, blue] (N5C5) -- (N6C6);
  \draw[arrow style, green!60!black] (N5C5) -- (N6C5);
  \draw[arrow style, blue] (N5C4) -- (N6C5);
  \draw[arrow style, green!60!black] (N5C6) -- (N6C7);
  \draw[arrow style, blue] (N5C6) -- (N6C8);
\end{tikzpicture}}
\caption{Affine crystal graph of highest weight \( \Lambda_0 + \Lambda_1 \)}
\label{fig:level2_grounded_crystal}
\end{figure}

In Section~\ref{section:partition_model}, we construct an affine crystal structure of highest weight \(\Lambda_0+\Lambda_1\) on \(\mathcal{P}_{2,b}\), shown in Figure~\ref{fig:level2_grounded_crystal}. Slight modifications of \(\mathcal{P}_{2,c}\) and \(\mathcal{P}_{2,a}\) give models for highest weights \(2\Lambda_0\) and \(2\Lambda_1\), respectively. The usual bracketing rule on addable and removable corners does not apply directly to these grounded partitions. We instead associate a binary word, called the \(\type\), to each grounded partition and use this word to determine the arrows. We prove that these explicit rules define the desired affine crystals (Theorems~\ref{thm:grounded_crystal_b} and~\ref{thm:grounded_crystal_ac}) by constructing explicit graph isomorphisms with the asymptotic model.

To obtain $q$-series identities, we erase the green arrows, which are the $0$-arrows, from the affine crystal graphs on $\mathcal{P}_{1,b}$ and $\mathcal{P}_{2,b}$. For the crystal $\mathcal{B}(\Lambda)$ of $L(\Lambda)$, this gives the crystal of the restriction of $L(\Lambda)$ to the Levi subalgebra~\cite[Proposition~4.3]{Kas95}. Its connected components are finite blue strings. We call a vertex with no incoming blue arrow the \emph{starting point} of its string.

To describe the generating functions that arise from this restriction,
we need some notation. For an integer \( n \geq 0 \), the
\emph{$q$-integer} $[n]_q$  and the \emph{$q$-factorial} $[n]_q!$ are defined by
\[
  [n]_q\coloneq 1+q+\cdots+q^{n-1}, \qand [n]_q! \coloneq
  [n]_q\,[n-1]_q\cdots [1]_q.
\]
For two integers $n\geq k\geq0$, the \emph{$q$-binomial coefficient} is 
\[\begin{bmatrix} n \\ k \end{bmatrix}_q \coloneq \frac{[n]_q!}{[k]_q![n-k]_q!},\]
and, for \( n \geq k \geq 0 \), the \emph{$q$-Yamanouchi number} $Y_q(n,k)$ is defined by
\[
Y_q(n,k) \coloneq \frac{[n-k+1]_q}{[n+1]_q}\begin{bmatrix} n+k \\ n \end{bmatrix}_q.
\]

In Section~\ref{section:decomp}, the decomposition into strings gives the following infinite-sum expressions for the generating functions of \( \mathcal{P}_{1,b} \) and \( \mathcal{P}_{2,b} \). The $q$-Yamanouchi numbers are, up to a power of $q$, the generating functions of the modified major index on Yamanouchi words, defined in Subsection~\ref{subsection:major index}.

\begin{theorem}\label{thm:decomp1} We have
\begin{align*}
  (-q;q)_\infty
  &= \sum_{\substack{n\geq k\geq 0}} q^{n(2n-1) + 3k}\, Y_{q^2}(2n-k,k)\,(q^2;q^2)_{2n}^{-1}\,[2(n-k)+1]_q \\
  &= \sum_{\substack{n\geq k\geq 0}} q^{(2n+1)n + 3k}\, Y_{q^2}(2n+1-k,k)\,(q^2;q^2)_{2n+1}^{-1}\,[2(n-k)+2]_q.
\end{align*}
\end{theorem}

\begin{theorem}\label{thm:decomp2} We have
\begin{align*}
\frac{(-q;q^2)_\infty}{(q;q^2)_\infty} = \sum_{\substack{n\geq k\geq 0}}q^{2k+(n-k)^2}\,Y_{q^4}(n,k)\,(q^2;q^2)_{n+k}^{-1}\,[2(n-k)+2]_q.
\end{align*}
\end{theorem}

\begin{remark}
For $M\geq0$, the \emph{string function} $G^{\Lambda}_M(q)$ is the generating function, by size, of the starting points of the blue strings with $M$ arrows. The sizes along such a string are $|\pi|,|\pi|+1,\dots,|\pi|+M$, so
\begin{align}\label{eq:string_decomp}
\sum_{\pi\in\mathcal{B}(\Lambda)} q^{|\pi|} = \sum_{M\geq0} G^{\Lambda}_M(q)\,(1+q+\dots+q^{M}).
\end{align}
The \emph{charge} of Lascoux and Schützenberger~\cite{LS78} is a statistic on semistandard Young tableaux; on standard Young tableaux its generating function is, up to a power of $q$, that of the major index~\cite[Chapter~III, Section~6, Example~2]{Mac95}, \cite[Corollary~7.21.5]{Sta1}. The \emph{Kostka--Foulkes polynomial} $K_{\lambda\mu}(q)$ is the generating function of charge on the semistandard Young tableaux of shape $\lambda$ and content $\mu$~\cite{LS78}. Within the Demazure subcrystal $\mathcal{B}_{(s_1s_0)^N}(\Lambda_0)$, the generating function of the starting points by the number of green arrows is, up to a power of $q$, a Kostka--Foulkes polynomial~\cite[Theorem~5.1]{KMOTU00}, and its limit as $N\to\infty$ is the string function~\cite[Theorem~1.1]{NY97}.

Theorem~\ref{thm:decomp1} makes the string functions of~\cite[Theorem~1.1, p.~550, and Corollary~4.8, p.~572]{NY97} explicit for $\Lambda_0$ and $\Lambda_1$. Theorem~\ref{thm:decomp2} gives those of $\Lambda_0+\Lambda_1$, a highest weight not of the form $\ell\Lambda_i$.
\end{remark}

The remainder of the paper is organised as follows. In Section~\ref{section:partition_model}, we construct affine crystal graphs of type~$A_1^{(1)}$ and highest weights $2\Lambda_0$, $\Lambda_0+\Lambda_1$, and $2\Lambda_1$, whose vertices are grounded partitions at level \( 2 \). In Section~\ref{section:decomp}, we prove the $q$-series identities in Theorems~\ref{thm:decomp1} and~\ref{thm:decomp2} by decomposing the grounded partition model at levels~$1$ and~$2$ into crystal graphs of finite type \( A_1 \).

\section*{Acknowledgements}

The authors warmly thank Jehanne Dousse for her guidance throughout as well as her valuable comments on the earlier versions of this paper. We are also grateful to Thomas Gerber, Frédéric Jouhet, Christian Krattenthaler, Philippe Nadeau and Ali Uncu for valuable discussions and insightful comments that contributed to the development of this work. The authors are funded by the SNSF Eccellenza grant of Jehanne Dousse, PCEFP2 202784.

\section{Affine crystal graphs on grounded partitions at level \( 2 \)}\label{section:partition_model}

\begin{figure}
    \centering
{\color{black}\begin{tikzpicture}[>=stealth, thick, scale=0.55]
  %%%%%%%%%%%%%%
  % Coordinates (unchanged)
  %%%%%%%%%%%%%%
  % Row 1: 1 vertex
  \coordinate (R1C1) at (0,0);
  
  % Row 2: 1 vertex
  \coordinate (R2C1) at (-2,-2.5);
  \coordinate (R2C2) at (2,-2.5);
  
  % Row 3: 1 vertex
  \coordinate (R3C1) at (-2,-5);
  \coordinate (R3C2) at (2,-5);
  
  % Row 4: 2 vertices
  \coordinate (R4C1) at (-6,-7.5);
  \coordinate (R4C2) at (-2,-7.5);
  \coordinate (R4C3) at (2,-7.5);
  \coordinate (R4C4) at (6,-7.5);
  
  % Row 5: 2 vertices
  \coordinate (R5C1) at (-10,-10);
  \coordinate (R5C2) at (-6,-10);
  \coordinate (R5C3) at (-2,-10);
  \coordinate (R5C4) at (2,-10);
  \coordinate (R5C5) at (6,-10);
  \coordinate (R5C6) at (10,-10);

    % Row 6: 2 vertices
  \coordinate (R6C1) at (-14,-13);
  \coordinate (R6C2) at (-10,-13);
  \coordinate (R6C3) at (-6,-13);
  \coordinate (R6C4) at (-2,-13);
  \coordinate (R6C5) at (2,-13);
  \coordinate (R6C6) at (6,-13);
  \coordinate (R6C7) at (10,-13);
  \coordinate (R6C8) at (14,-13);
  
  %%%%%%%%%%%%%%
  % Nodes (transposed shapes)
  %%%%%%%%%%%%%%
  \tikzset{mynode/.style={inner sep=0pt, font=\scriptsize}}
  
  \node[mynode] (N1C1) at (R1C1) {\(\emptyset\hspace{0.5em}\emptyset\)};
  
  \node[mynode] (N2C1) at (R2C1) {\(\scalebox{0.5}{\begin{ytableau}
*(green)
\end{ytableau}} \hspace{0.5em}\emptyset\)};

  \node[mynode] (N2C2) at (R2C2) {\(\emptyset\hspace{0.5em}\scalebox{0.5}{\begin{ytableau}
*(blue)
\end{ytableau}}\)};
  
  \node[mynode] (N3C1) at (R3C1) {\(\scalebox{0.5}{\begin{ytableau}
*(green)
\end{ytableau}}\hspace{0.5em}\scalebox{0.5}{\begin{ytableau}
*(blue)
\end{ytableau}}\)};

  \node[mynode] (N3C2) at (R3C2) {\(\emptyset\hspace{0.5em}\scalebox{0.5}{\begin{ytableau}
*(blue) \\ *(green)
\end{ytableau}}\)};
  
  \node[mynode] (N4C1) at (R4C1) {\(\scalebox{0.5}{\begin{ytableau}
*(green)
\end{ytableau}}\hspace{0.5em}\scalebox{0.5}{\begin{ytableau}
*(blue) \\ *(green)
\end{ytableau}}\)};

  \node[mynode] (N4C2) at (R4C2) {\(\scalebox{0.5}{\begin{ytableau}
*(green) \\ *(blue)
\end{ytableau}}\hspace{0.5em}\scalebox{0.5}{\begin{ytableau}
*(blue)
\end{ytableau}}\)};

  \node[mynode] (N4C3) at (R4C3) {\(\emptyset\hspace{0.5em}\scalebox{0.5}{\begin{ytableau}
*(blue) & *(green) \\ *(green)
\end{ytableau}}\)};

  \node[mynode] (N4C4) at (R4C4) {\(\emptyset\hspace{0.5em}\scalebox{0.5}{\begin{ytableau}
*(blue) \\ *(green) \\ *(blue)
\end{ytableau}}\)};
  
  \node[mynode] (N5C1) at (R5C1) {\(\scalebox{0.5}{\begin{ytableau}
*(green)
\end{ytableau}}\hspace{0.5em}\scalebox{0.5}{\begin{ytableau}
*(blue) \\ *(green) \\ *(blue)
\end{ytableau}}\)};

  \node[mynode] (N5C3) at (R5C3) {\(\scalebox{0.5}{\begin{ytableau}
*(green) & *(blue) \\ *(blue)
\end{ytableau}}\hspace{0.5em}\scalebox{0.5}{\begin{ytableau}
*(blue)
\end{ytableau}}\)};

  \node[mynode] (N5C2) at (R5C2) {\(\scalebox{0.5}{\begin{ytableau}
*(green) \\ *(blue)
\end{ytableau}}\hspace{0.5em}\scalebox{0.5}{\begin{ytableau}
*(blue) \\ *(green)
\end{ytableau}}\)};

  \node[mynode] (N5C4) at (R5C4) {\(\scalebox{0.5}{\begin{ytableau}
*(green)
\end{ytableau}}\hspace{0.5em}\scalebox{0.5}{\begin{ytableau}
*(blue)  & *(green) \\ *(green)
\end{ytableau}}\)};

  \node[mynode] (N5C5) at (R5C5) {\(\emptyset\hspace{0.5em}\scalebox{0.5}{\begin{ytableau}
*(blue) & *(green) \\ *(green) \\ *(blue)
\end{ytableau}}\)};

  \node[mynode] (N5C6) at (R5C6) {\(\emptyset\hspace{0.5em}\scalebox{0.5}{\begin{ytableau}
*(blue) \\ *(green) \\ *(blue) \\ *(green)
\end{ytableau}}\)};

  \node[mynode] (N6C1) at (R6C1) {\(\scalebox{0.5}{\begin{ytableau}
*(green)
\end{ytableau}}\hspace{0.5em}\scalebox{0.5}{\begin{ytableau}
*(blue) \\ *(green) \\ *(blue) \\ *(green)
\end{ytableau}}\)};

  \node[mynode] (N6C2) at (R6C2) {\(\scalebox{0.5}{\begin{ytableau}
*(green) \\ *(blue)
\end{ytableau}}\hspace{0.5em}\scalebox{0.5}{\begin{ytableau}
*(blue) \\ *(green) \\ *(blue)
\end{ytableau}}\)};

  \node[mynode] (N6C4) at (R6C4) {\(\scalebox{0.5}{\begin{ytableau}
*(green) & *(blue) \\ *(blue)
\end{ytableau}}\hspace{0.5em}\scalebox{0.5}{\begin{ytableau}
*(blue) \\ *(green)
\end{ytableau}}\)};

  \node[mynode] (N6C3) at (R6C3) {\(\scalebox{0.5}{\begin{ytableau}
*(green) \\ *(blue)
\end{ytableau}}\hspace{0.5em}\scalebox{0.5}{\begin{ytableau}
*(blue) & *(green) \\ *(green)
\end{ytableau}}\)};

  \node[mynode] (N6C6) at (R6C6) {\(\emptyset\hspace{0.5em}\scalebox{0.5}{\begin{ytableau}
*(blue) & *(green) \\ *(green) & *(blue) \\ *(blue)
\end{ytableau}}\)};

  \node[mynode] (N6C5) at (R6C5) {\(\scalebox{0.5}{\begin{ytableau}
*(green)
\end{ytableau}}\hspace{0.5em}\scalebox{0.5}{\begin{ytableau}
*(blue) & *(green) \\ *(green) \\ *(blue)
\end{ytableau}}\)};

  \node[mynode] (N6C7) at (R6C7) {\(\emptyset\hspace{0.5em}\scalebox{0.5}{\begin{ytableau}
*(blue) & *(green) \\ *(green) \\ *(blue) \\ *(green)
\end{ytableau}}\)};

  \node[mynode] (N6C8) at (R6C8) {\(\emptyset\hspace{0.5em}\scalebox{0.5}{\begin{ytableau}
*(blue) \\ *(green) \\ *(blue) \\ *(green) \\ *(blue)
\end{ytableau}}\)};
  
  %%%%%%%%%%%%%%
  % Arrows (unchanged)
  %%%%%%%%%%%%%%
  \tikzset{arrow style/.style={->, shorten <=1pt, shorten >=1pt, line width=0.8pt}}
  \draw[arrow style, blue] (N1C1) -- (N2C2);
  \draw[arrow style, green!60!black] (N1C1) -- (N2C1);
  
  \draw[arrow style, blue] (N2C1) -- (N3C1);
\draw[arrow style, green!60!black] (N2C2) -- (N3C2);
  
  \draw[arrow style, blue] (N3C1) -- (N4C2);
  \draw[arrow style, green!60!black] (N3C1) -- (N4C1);
  \draw[arrow style, blue] (N3C2) -- (N4C4);
  \draw[arrow style, green!60!black] (N3C2) -- (N4C3);
  
  \draw[arrow style, blue] (N4C1) -- (N5C1);
  \draw[arrow style, blue] (N4C2) -- (N5C3);
  \draw[arrow style, green!60!black] (N4C2) -- (N5C2);
  \draw[arrow style, blue] (N4C3) -- (N5C5);
\draw[arrow style, green!60!black] (N4C3) -- (N5C4);
\draw[arrow style, green!60!black] (N4C4) -- (N5C6);

  \draw[arrow style, green!60!black] (N5C1) -- (N6C1);
  \draw[arrow style, blue] (N5C1) -- (N6C2);
  \draw[arrow style, green!60!black] (N5C2) -- (N6C3);
  \draw[arrow style, blue] (N5C2) -- (N6C4);
  \draw[arrow style, green!60!black] (N5C3) -- (N6C4);
  \draw[arrow style, blue] (N5C5) -- (N6C6);
  \draw[arrow style, green!60!black] (N5C5) -- (N6C5);
  \draw[arrow style, blue] (N5C4) -- (N6C5);
  \draw[arrow style, green!60!black] (N5C6) -- (N6C7);
  \draw[arrow style, blue] (N5C6) -- (N6C8);
\end{tikzpicture}}
    \caption{The graph $\mathcal{B}_{2,b}^{\mathcal{A}}$}
    \label{fig:pairs_partition_2}
\end{figure}
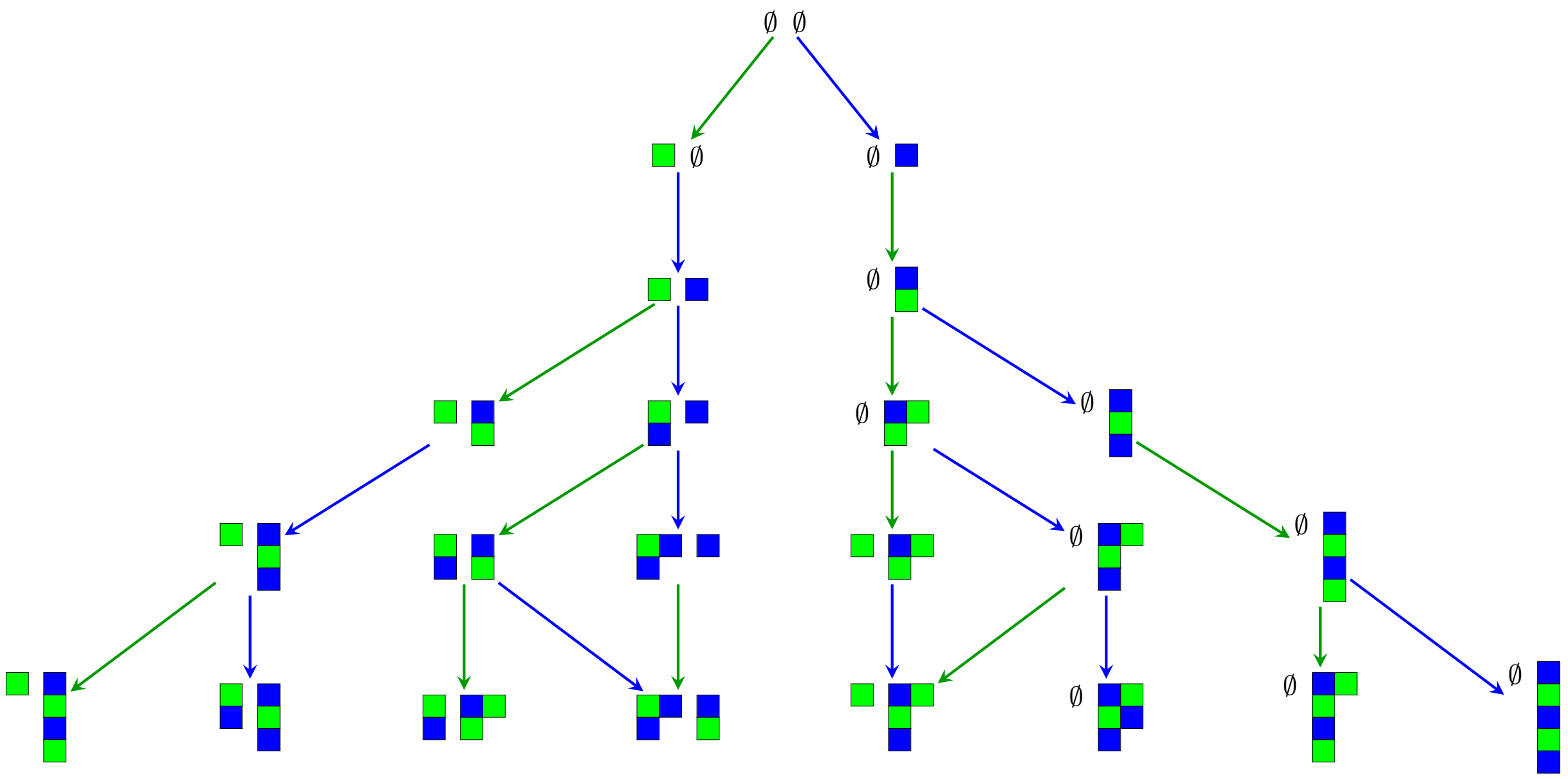

In this section, we construct affine crystal graphs of type \( A_1^{(1)} \) and highest weights \( 2\Lambda_0 \), \( \Lambda_0 + \Lambda_1 \), and \( 2\Lambda_1 \) on grounded partitions at level \( 2 \). We first treat \( \Lambda_0 + \Lambda_1 \), whose grounded partition model and asymptotic model appear in Figures~\ref{fig:level2_grounded_crystal} and~\ref{fig:pairs_partition_2}. We define explicit rules for the arrows and construct an isomorphism between these models, then adapt the construction to the other two highest weights.

First, we reformulate the description of the affine crystal graph \(\mathcal{B}_{1,b}^{\mathcal{P}}\) in Definition~\ref{def:partition_model_1}. 

\begin{definition}\label{def:partition_type}
  Let \( \lambda \) be a Young diagram with a blue--green checkerboard colouring, whose top-left cell may be green or blue. Its \emph{type}, denoted \( \type(\lambda) \), is the following binary word. Read the columns from left to right, writing \( 0 \) for a green cell at the bottom and \( 1 \) for a blue cell at the bottom. If the first row ends in green, append one additional \( 0 \).

  For the empty partition \( \emptyset \), we set
  \( \type(\emptyset) \coloneq \emptyset \) if the top-left cell is green, and
  \( \type(\emptyset) \coloneq 0 \) if the top-left cell is blue.
\end{definition}

\begin{example}
  A few examples of $2$-restricted partitions and their types, with the top-left cell green,
  are shown below:
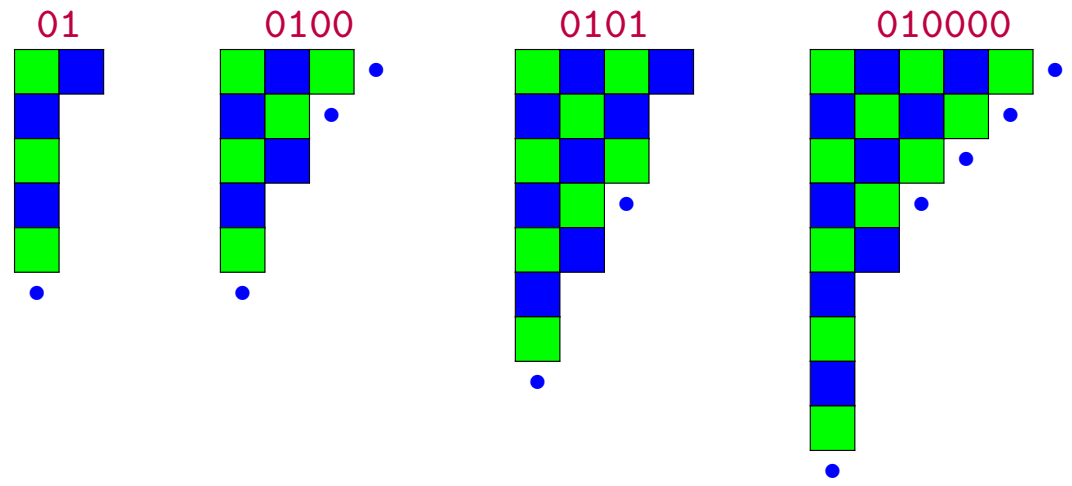
\begin{align*}
\overset{\textcolor{purple}{\texttt{01}}}{\scalebox{0.5}{\begin{ytableau}
  *(green)&*(blue)\\
  *(blue)\\
  *(green)\\
  *(blue)\\
  *(green)\\
  \none[{\color{blue}\scalebox{1.25}{$\bullet$}}]
\end{ytableau}}}
\qquad\overset{\textcolor{purple}{\texttt{0100}}}{\scalebox{0.5}{\begin{ytableau}
  *(green)&*(blue)&*(green)&\none[{\color{blue}\scalebox{1.25}{$\bullet$}}]\\
  *(blue)&*(green)&\none[{\color{blue}\scalebox{1.25}{$\bullet$}}]\\
  *(green)&*(blue)\\
  *(blue)\\
  *(green)\\
  \none[{\color{blue}\scalebox{1.25}{$\bullet$}}]
\end{ytableau}}}
\qquad
\overset{\textcolor{purple}{\texttt{0101}}}{\scalebox{0.5}{\begin{ytableau}
  *(green)&*(blue)&*(green)&*(blue)\\
  *(blue)&*(green)&*(blue)\\
  *(green)&*(blue)&*(green)\\
  *(blue)&*(green)&\none[{\color{blue}\scalebox{1.25}{$\bullet$}}]\\
  *(green)&*(blue)\\
  *(blue)\\
  *(green)\\
  \none[{\color{blue}\scalebox{1.25}{$\bullet$}}]
\end{ytableau}}}
\qquad
\overset{\textcolor{purple}{\texttt{010000}}}{\scalebox{0.5}{\begin{ytableau}
  *(green)&*(blue)&*(green)&*(blue)&*(green)&\none[{\color{blue}\scalebox{1.25}{$\bullet$}}]\\
  *(blue)&*(green)&*(blue)&*(green)&\none[{\color{blue}\scalebox{1.25}{$\bullet$}}]\\
  *(green)&*(blue)&*(green)&\none[{\color{blue}\scalebox{1.25}{$\bullet$}}]\\
  *(blue)&*(green)&\none[{\color{blue}\scalebox{1.25}{$\bullet$}}]\\
  *(green)&*(blue)\\
  *(blue)\\
  *(green)\\
  *(blue)\\
  *(green)\\
  \none[{\color{blue}\scalebox{1.25}{$\bullet$}}]
\end{ytableau}}}
\end{align*}
  For a $2$-restricted partition \(\lambda\), each column has a removable or addable blue corner. Consequently, the sequence of brackets
  associated to \(\type(\lambda)\) is the same as the sequence of
  brackets obtained from addable and removable blue corners of
  \(\lambda\) according to Definition~\ref{def:partition_model_1}.
\end{example}

\begin{definition}\label{def:level1_type_crystal}
  The directed graph \( \mathcal{B}_{1,b}^{\mathcal{P}} \), whose vertices are partitions and whose edges are coloured blue and green,
  is defined recursively starting from the empty partition
  \( \emptyset \). For each vertex \( \lambda \), consider the blue--green checkerboard where the top-left cell is green. The outgoing arrows
  are determined by the following bracketing rules:
\begin{enumerate}
\item \textbf{Blue arrow.} Assign a left bracket ``('' to each \( 0 \)
  in \(\mathrm{type}(\lambda)\), and a right bracket ``)'' to each
  \( 1 \). Take the leftmost unmatched \( 0 \), and add a blue cell at
  the bottom of the corresponding column of \( \lambda \). A blue arrow is drawn from \( \lambda \) to the resulting partition.
\item \textbf{Green arrow.} Defined in the same way, using the binary
  word obtained by interchanging the roles of blue and green cells in the
  definition of \(\mathrm{type}(\lambda)\).
\end{enumerate}
\end{definition}

As mentioned below Definition~\ref{def:partition_model_1}, the vertex
set of the graph \( \mathcal{B}_{1,b}^{\mathcal{P}} \) is the set of
$2$-restricted partitions.

\begin{definition}\label{def:pairs_crystal}
  Consider a pair of partitions \( (\lambda^{(1)}, \lambda^{(2)}) \),
  each with a blue--green checkerboard colouring. The top-left cell is
  green in \( \lambda^{(1)} \) and blue in \( \lambda^{(2)} \). Define
  the \emph{type} of the pair by reading the type of
  \( \lambda^{(2)} \) before that of \( \lambda^{(1)} \):
  \[
    \type(\lambda^{(1)}, \lambda^{(2)}) \coloneq \type(\lambda^{(2)})\type(\lambda^{(1)}).
  \]
  
  The directed graph \( \mathcal{B}_{2,b}^{\mathcal{A}} \), whose
  vertices are pairs of partitions and whose edges are coloured blue
  and green, is defined recursively starting from the empty pair
  \( (\lambda^{(1)},\lambda^{(2)}) = (\emptyset,\emptyset) \). For each vertex
  \( (\lambda^{(1)},\lambda^{(2)}) \), the outgoing arrows are determined by
  the following bracketing rules:
\begin{enumerate}
\item \textbf{Blue arrow.} Assign a left bracket ``('' to each letter
  \( 0 \) in \( \type(\lambda^{(1)}, \lambda^{(2)}) \), and a right bracket
  ``)'' to each \( 1 \). The \( j \)th letter of \( \type(\lambda^{(2)}) \)
  corresponds to the \( j \)th column of \( \lambda^{(2)} \), and the
  \( j \)th letter of \( \type(\lambda^{(1)}) \) to the \( j \)th column
  of \( \lambda^{(1)} \). Take the leftmost unmatched \( 0 \), and add
  a blue cell at the bottom of the corresponding column. A blue arrow is drawn from \( (\lambda^{(1)}, \lambda^{(2)}) \) to the resulting pair of partitions.
\item \textbf{Green arrow.} Defined in the same way, using the binary
  word obtained by interchanging the roles of blue and green cells in
  the definition of \( \type(\lambda^{(1)}, \lambda^{(2)}) \).
\end{enumerate}
\end{definition}

\begin{remark}
  \begin{enumerate}
  \item It follows from the signature rule~\cite[Section~2.4]{BS17} for the tensor product of crystals that the graph
    \(\mathcal{B}_{2,b}^{\mathcal{A}}\) is isomorphic to a connected
    component of the tensor product
    \(\mathcal{B}(\Lambda_0)\otimes\mathcal{B}(\Lambda_1)\). By
    Wakimoto~\cite[Corollary~7.4.6]{Wak01}, every connected component
    of \(\mathcal{B}(\Lambda_0)\otimes\mathcal{B}(\Lambda_1)\) is
    isomorphic to \(\mathcal{B}(\Lambda_0+\Lambda_1)\).
  \item If a pair \((\lambda^{(1)},\lambda^{(2)})\) appears as a vertex of
    \(\mathcal{B}_{2,b}^{\mathcal{A}}\), then both \(\lambda^{(1)}\) and
    \(\lambda^{(2)}\) are $2$-restricted partitions. However, not
    every pair of $2$-restricted partitions appears as a vertex of
    \(\mathcal{B}_{2,b}^{\mathcal{A}}\).
  \end{enumerate}
\end{remark}

We represent a grounded partition \( \pi \in \mathcal{P}_{2,b} \) by
the Young diagram of its underlying partition, filled with entries
$a,b,c$. We fill the rightmost cell in each row with the label ($a$,
$b$ or $c$) of the corresponding part. Note that all odd columns have
an entry $a$ or $c$ at the bottom: the bottom cell of the
\( (2i+1) \)th column is assigned the label of the first occurrence of
the part \( (2i+1) \) in \( \pi \). In each column with an odd index,
we fill the remaining entries such that the entries along the column
alternate between \( a \) and \( c \). In each column with an even
index, we put an entry \( b \) in each cell. This completely
determines the entries in the diagram. Given such a diagram, we can
recover the grounded partition $\pi$ by reading the rows from bottom
to top to obtain the parts and by reading the last entry in each row
to obtain the label of each part. From now on, we identify
\( \mathcal{P}_{2,b} \) with the set of such diagrams.

We extend this description to Young diagrams of arbitrary shape. The
entries determine a blue--green colouring, which in turn allows us to
define the \emph{type} of a diagram filled with entries in
\(\{a,b,c\}\).

\begin{definition}\label{def:P2_type}
  Let \( \pi \) be a partition whose Young diagram is filled with
  entries in \(\{a,b,c\}\). The following rules specify the allowed
  entries, the blue--green colouring of the cells, and define the type
  of the resulting diagram.
  \begin{enumerate}
  \item (Entries) All cells in columns with odd indices have entry \( a \)
    or \( c \), while all cells in columns with even indices have entry
    \( b \). As an initial condition, all cells in the \( 0 \)th column
    have entry \( b \).
  \item (Colours) For convenience, we call a cell with entry
    \( x \in \{a,b,c\} \) an \( x \)-cell. Each cell is coloured either
    green or blue according to its entry, as follows:
    \begin{itemize}
    \item An \( a \)-cell is coloured green, and a \( c \)-cell is coloured blue.
    \item A \( b \)-cell is coloured blue if the cell to its left is an \( a \)-cell, and green if the cell to its left is a \( c \)-cell.
    \end{itemize}
  \item (Type) The \emph{type} of \(\pi\), denoted \( \type(\pi) \), is
    the binary word defined as follows:
    \begin{itemize}
    \item Read the columns of \( \pi \) from left to right. Assign \( 0 \)
      to a column if its bottom cell is green, and \( 1 \) if its bottom
      cell is blue.
    \item If the rightmost cell in the first row is an \( a \)-cell,
      append \( 00 \) at the end of the resulting word. If it is a \( b \)-cell,
      append \( 0 \).
    \end{itemize}
  \end{enumerate}
\end{definition}

The convention that the \( 0 \)th column consists of \( b \)-cells
gives \( \type(\emptyset) = 0 \). Cells in columns \( 2i+1 \) and
\( 2i+2 \) in the same row have opposite colours. The entries
therefore determine the colouring; conversely, the colouring and
column indices uniquely determine the entries \( a,b,c \).

We now define a directed graph whose vertices are the diagrams with
these entries and colouring, and whose structure is determined by
their types.

\begin{definition}\label{def:P2_crystal}
  The directed graph \(\mathcal{B}_{2,b}^{\mathcal{P}}\), with blue
  and green arrows, is defined recursively starting from the empty
  partition \( \emptyset \). The vertices are diagrams with entries
  \( a,b,c \) and colouring as in Definition~\ref{def:P2_type},
  obtained inductively during the construction. For each vertex
  \(\pi\), the outgoing arrows are determined by the following
  bracketing rules:
\begin{enumerate}
\item \textbf{Blue arrow.} Given \( \type(\pi) \), assign a left
  bracket ``('' to each letter \( 0 \), and a right bracket ``)'' to
  each \( 1 \). The \( j \)th letter of \( \type(\pi) \) corresponds
  to the \( j \)th column. Take the leftmost unmatched \( 0 \), and
  add a blue cell at the bottom of the corresponding column. Since the
  appended letters are unmatched, the leftmost unmatched \( 0 \) is
  never the second of two appended letters. A blue arrow is drawn from
  $\pi$ to the resulting diagram.
\item \textbf{Green arrow.} Defined in the same way, using the binary
  word obtained by interchanging the roles of blue and green cells in
  the definition of \(\mathrm{type}(\pi)\); and interchanging the
  roles of \(a\)-cells and \(c\)-cells at the end of the first row.
\end{enumerate}
Here, the entry of the newly added cell is determined by
Definition~\ref{def:P2_type}: if the column has odd index, the entry
is either \(a\) or \(c\), depending on the colour of the cell; if the
column has even index, the entry is \( b \). The graph
\( \mathcal{B}_{2,b}^{\mathcal{P}} \) is shown in
Figure~\ref{fig:level2_grounded_crystal_diags}.
\end{definition}
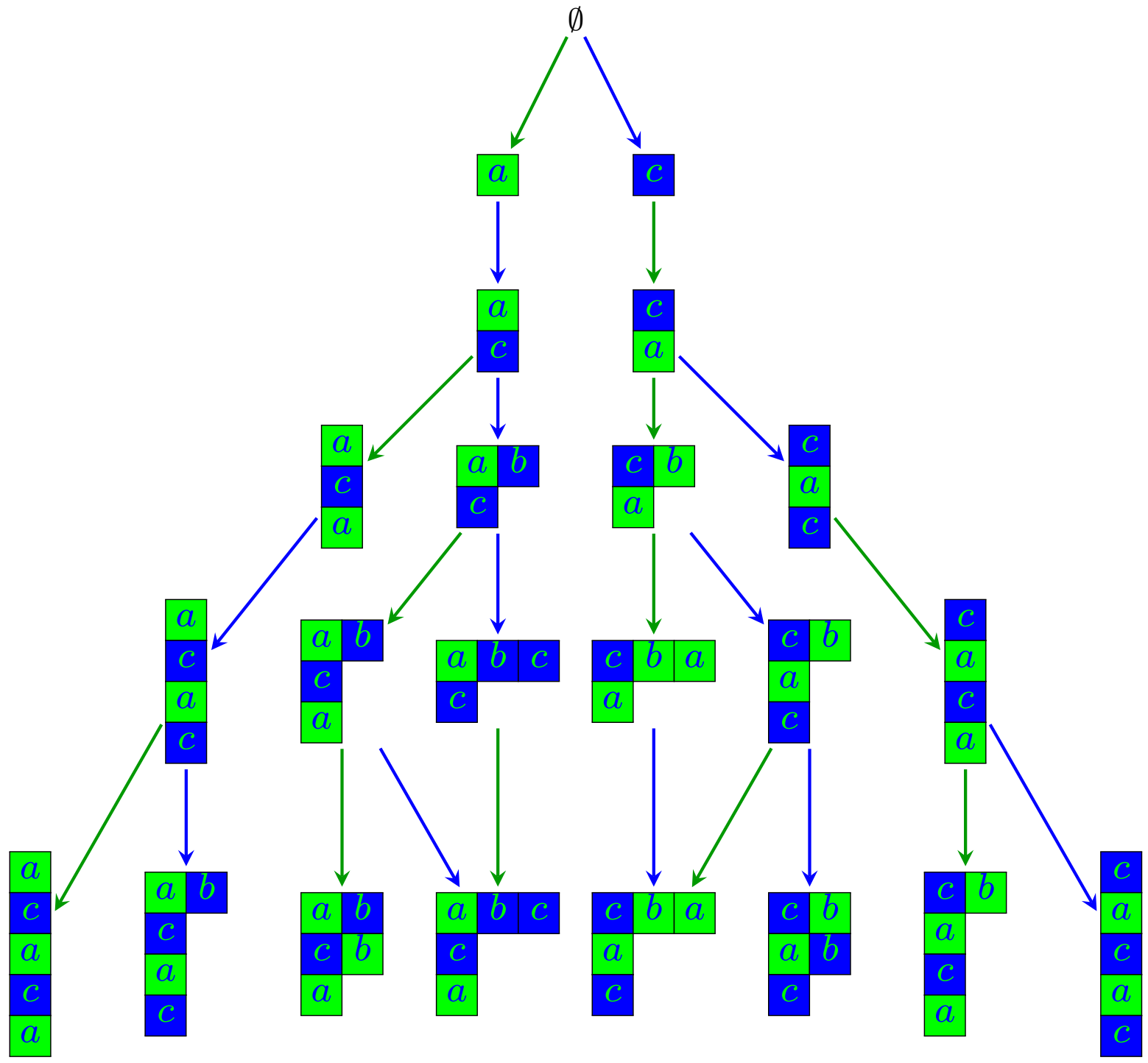
\begin{figure}
\centering
{\color{black}\begin{tikzpicture}[>=stealth, thick, scale=0.7]
  %%%%%%%%%%%%%%
  % Coordinates (unchanged)
  %%%%%%%%%%%%%%
  % Row 1: 1 vertex
  \coordinate (R1C1) at (0,0);
  
  % Row 2: 1 vertex
  \coordinate (R2C1) at (-1,-2);
  \coordinate (R2C2) at (1,-2);
  
  % Row 3: 1 vertex
  \coordinate (R3C1) at (-1,-4);
  \coordinate (R3C2) at (1,-4);
  
  % Row 4: 2 vertices
  \coordinate (R4C1) at (-3,-6);
  \coordinate (R4C2) at (-1,-6);
  \coordinate (R4C3) at (1,-6);
  \coordinate (R4C4) at (3,-6);
  
  % Row 5: 2 vertices
  \coordinate (R5C1) at (-5,-8.5);
  \coordinate (R5C2) at (-3,-8.5);
  \coordinate (R5C3) at (-1,-8.5);
  \coordinate (R5C4) at (1,-8.5);
  \coordinate (R5C5) at (3,-8.5);
  \coordinate (R5C6) at (5,-8.5);

    % Row 6: 2 vertices
  \coordinate (R6C1) at (-7,-12);
  \coordinate (R6C2) at (-5,-12);
  \coordinate (R6C3) at (-3,-12);
  \coordinate (R6C4) at (-1,-12);
  \coordinate (R6C5) at (1,-12);
  \coordinate (R6C6) at (3,-12);
  \coordinate (R6C7) at (5,-12);
  \coordinate (R6C8) at (7,-12);
  
  %%%%%%%%%%%%%%
  % Nodes (transposed shapes)
  %%%%%%%%%%%%%%
  \tikzset{mynode/.style={inner sep=0pt, font=\scriptsize}}
  
  \node[mynode] (N1C1) at (R1C1) {\(\emptyset\)};
  
  \node[mynode] (N2C1) at (R2C1) {\(\scalebox{0.8}{\begin{ytableau}
*(green) \scalebox{1.5}{\color{blue}$a$}
\end{ytableau}}\)};

  \node[mynode] (N2C2) at (R2C2) {\(\scalebox{0.8}{\begin{ytableau}
*(blue) \scalebox{1.5}{\color{green}$c$}
\end{ytableau}}\)};
  
  \node[mynode] (N3C1) at (R3C1) {\(\scalebox{0.8}{\begin{ytableau}
*(green) \scalebox{1.5}{\color{blue}$a$}\\
*(blue) \scalebox{1.5}{\color{green}$c$}
\end{ytableau}}\)};

  \node[mynode] (N3C2) at (R3C2) {\(\scalebox{0.8}{\begin{ytableau}
*(blue) \scalebox{1.5}{\color{green}$c$}\\
*(green) \scalebox{1.5}{\color{blue}$a$}
\end{ytableau}}\)};
  
  \node[mynode] (N4C1) at (R4C1) {\(\scalebox{0.8}{\begin{ytableau}
*(green) \scalebox{1.5}{\color{blue}$a$}\\
*(blue) \scalebox{1.5}{\color{green}$c$}\\
*(green) \scalebox{1.5}{\color{blue}$a$}
\end{ytableau}}\)};

  \node[mynode] (N4C2) at (R4C2) {\(\scalebox{0.8}{\begin{ytableau}
*(green) \scalebox{1.5}{\color{blue}$a$} & *(blue) \scalebox{1.5}{\color{green}$b$}\\
*(blue) \scalebox{1.5}{\color{green}$c$}
\end{ytableau}}\)};

  \node[mynode] (N4C3) at (R4C3) {\(\scalebox{0.8}{\begin{ytableau}
*(blue) \scalebox{1.5}{\color{green}$c$} & *(green) \scalebox{1.5}{\color{blue}$b$}\\
*(green) \scalebox{1.5}{\color{blue}$a$}
\end{ytableau}}\)};

  \node[mynode] (N4C4) at (R4C4) {\(\scalebox{0.8}{\begin{ytableau}
*(blue) \scalebox{1.5}{\color{green}$c$}\\
*(green) \scalebox{1.5}{\color{blue}$a$}\\
*(blue) \scalebox{1.5}{\color{green}$c$}
\end{ytableau}}\)};
  
  \node[mynode] (N5C1) at (R5C1) {\(\scalebox{0.8}{\begin{ytableau}
*(green) \scalebox{1.5}{\color{blue}$a$}\\
*(blue) \scalebox{1.5}{\color{green}$c$}\\
*(green) \scalebox{1.5}{\color{blue}$a$}\\
*(blue) \scalebox{1.5}{\color{green}$c$}
\end{ytableau}}\)};

  \node[mynode] (N5C3) at (R5C3) {\(\scalebox{0.8}{\begin{ytableau}
*(green) \scalebox{1.5}{\color{blue}$a$} & *(blue) \scalebox{1.5}{\color{green}$b$} & *(blue) \scalebox{1.5}{\color{green}$c$}\\
*(blue) \scalebox{1.5}{\color{green}$c$}
\end{ytableau}}\)};

  \node[mynode] (N5C2) at (R5C2) {\(\scalebox{0.8}{\begin{ytableau}
*(green) \scalebox{1.5}{\color{blue}$a$} & *(blue) \scalebox{1.5}{\color{green}$b$} \\
*(blue) \scalebox{1.5}{\color{green}$c$}\\
*(green) \scalebox{1.5}{\color{blue}$a$}
\end{ytableau}}\)};

  \node[mynode] (N5C5) at (R5C5) {\(\scalebox{0.8}{\begin{ytableau}
*(blue) \scalebox{1.5}{\color{green}$c$} & *(green) \scalebox{1.5}{\color{blue}$b$} \\
*(green) \scalebox{1.5}{\color{blue}$a$}\\
*(blue) \scalebox{1.5}{\color{green}$c$}
\end{ytableau}}\)};

  \node[mynode] (N5C4) at (R5C4) {\(\scalebox{0.8}{\begin{ytableau}
*(blue) \scalebox{1.5}{\color{green}$c$} & *(green) \scalebox{1.5}{\color{blue}$b$} & *(green) \scalebox{1.5}{\color{blue}$a$}\\
*(green) \scalebox{1.5}{\color{blue}$a$}
\end{ytableau}}\)};

  \node[mynode] (N5C6) at (R5C6) {\(\scalebox{0.8}{\begin{ytableau}
*(blue) \scalebox{1.5}{\color{green}$c$} \\
*(green) \scalebox{1.5}{\color{blue}$a$} \\
*(blue) \scalebox{1.5}{\color{green}$c$} \\
*(green) \scalebox{1.5}{\color{blue}$a$}
\end{ytableau}}\)};

  \node[mynode] (N6C1) at (R6C1) {\(\scalebox{0.8}{\begin{ytableau}
*(green) \scalebox{1.5}{\color{blue}$a$} \\
*(blue) \scalebox{1.5}{\color{green}$c$} \\
*(green) \scalebox{1.5}{\color{blue}$a$} \\
*(blue) \scalebox{1.5}{\color{green}$c$} \\
*(green) \scalebox{1.5}{\color{blue}$a$}
\end{ytableau}}\)};

  \node[mynode] (N6C2) at (R6C2) {\(\scalebox{0.8}{\begin{ytableau}
*(green) \scalebox{1.5}{\color{blue}$a$} & *(blue) \scalebox{1.5}{\color{green}$b$}\\
*(blue) \scalebox{1.5}{\color{green}$c$} \\
*(green) \scalebox{1.5}{\color{blue}$a$} \\
*(blue) \scalebox{1.5}{\color{green}$c$}
\end{ytableau}}\)};

  \node[mynode] (N6C4) at (R6C4) {\(\scalebox{0.8}{\begin{ytableau}
*(green) \scalebox{1.5}{\color{blue}$a$} & *(blue) \scalebox{1.5}{\color{green}$b$} & *(blue) \scalebox{1.5}{\color{green}$c$} \\
*(blue) \scalebox{1.5}{\color{green}$c$} \\
*(green) \scalebox{1.5}{\color{blue}$a$}
\end{ytableau}}\)};

  \node[mynode] (N6C3) at (R6C3) {\(\scalebox{0.8}{\begin{ytableau}
*(green) \scalebox{1.5}{\color{blue}$a$} & *(blue) \scalebox{1.5}{\color{green}$b$} \\
*(blue) \scalebox{1.5}{\color{green}$c$} & *(green) \scalebox{1.5}{\color{blue}$b$}\\
*(green) \scalebox{1.5}{\color{blue}$a$}
\end{ytableau}}\)};

  \node[mynode] (N6C6) at (R6C6) {\(\scalebox{0.8}{\begin{ytableau}
*(blue) \scalebox{1.5}{\color{green}$c$} & *(green) \scalebox{1.5}{\color{blue}$b$} \\
*(green) \scalebox{1.5}{\color{blue}$a$} & *(blue) \scalebox{1.5}{\color{green}$b$}\\
*(blue) \scalebox{1.5}{\color{green}$c$}
\end{ytableau}}\)};

  \node[mynode] (N6C5) at (R6C5) {\(\scalebox{0.8}{\begin{ytableau}
*(blue) \scalebox{1.5}{\color{green}$c$} & *(green) \scalebox{1.5}{\color{blue}$b$} & *(green) \scalebox{1.5}{\color{blue}$a$}\\
*(green) \scalebox{1.5}{\color{blue}$a$} \\
*(blue) \scalebox{1.5}{\color{green}$c$}
\end{ytableau}}\)};

  \node[mynode] (N6C7) at (R6C7) {\(\scalebox{0.8}{\begin{ytableau}
*(blue) \scalebox{1.5}{\color{green}$c$} & *(green) \scalebox{1.5}{\color{blue}$b$}\\
*(green) \scalebox{1.5}{\color{blue}$a$} \\
*(blue) \scalebox{1.5}{\color{green}$c$} \\
*(green) \scalebox{1.5}{\color{blue}$a$}
\end{ytableau}}\)};

  \node[mynode] (N6C8) at (R6C8) {\(\scalebox{0.8}{\begin{ytableau}
*(blue) \scalebox{1.5}{\color{green}$c$} \\
*(green) \scalebox{1.5}{\color{blue}$a$} \\
*(blue) \scalebox{1.5}{\color{green}$c$} \\
*(green) \scalebox{1.5}{\color{blue}$a$} \\
*(blue) \scalebox{1.5}{\color{green}$c$}
\end{ytableau}}\)};
  
  %%%%%%%%%%%%%%
  % Arrows (unchanged)
  %%%%%%%%%%%%%%
  \tikzset{arrow style/.style={->, shorten <=1pt, shorten >=1pt, line width=0.8pt}}
  \draw[arrow style, blue] (N1C1) -- (N2C2);
  \draw[arrow style, green!60!black] (N1C1) -- (N2C1);
  
  \draw[arrow style, blue] (N2C1) -- (N3C1);
\draw[arrow style, green!60!black] (N2C2) -- (N3C2);
  
  \draw[arrow style, blue] (N3C1) -- (N4C2);
  \draw[arrow style, green!60!black] (N3C1) -- (N4C1);
  \draw[arrow style, blue] (N3C2) -- (N4C4);
  \draw[arrow style, green!60!black] (N3C2) -- (N4C3);
  
  \draw[arrow style, blue] (N4C1) -- (N5C1);
  \draw[arrow style, blue] (N4C2) -- (N5C3);
  \draw[arrow style, green!60!black] (N4C2) -- (N5C2);
  \draw[arrow style, blue] (N4C3) -- (N5C5);
\draw[arrow style, green!60!black] (N4C3) -- (N5C4);
\draw[arrow style, green!60!black] (N4C4) -- (N5C6);

  \draw[arrow style, green!60!black] (N5C1) -- (N6C1);
  \draw[arrow style, blue] (N5C1) -- (N6C2);
  \draw[arrow style, green!60!black] (N5C2) -- (N6C3);
  \draw[arrow style, blue] (N5C2) -- (N6C4);
  \draw[arrow style, green!60!black] (N5C3) -- (N6C4);
  \draw[arrow style, blue] (N5C5) -- (N6C6);
  \draw[arrow style, green!60!black] (N5C5) -- (N6C5);
  \draw[arrow style, blue] (N5C4) -- (N6C5);
  \draw[arrow style, green!60!black] (N5C6) -- (N6C7);
  \draw[arrow style, blue] (N5C6) -- (N6C8);
\end{tikzpicture}}
\caption{The graph $\mathcal{B}_{2,b}^{\mathcal{P}}$}
\label{fig:level2_grounded_crystal_diags}
\end{figure}

\begin{remark}
  The motivation for the definition of the graph
  \(\mathcal{B}_{2,b}^{\mathcal{P}}\) comes from the perfect crystal
  graph of type \(A_1^{(1)}\) at level \(2\), illustrated as
  
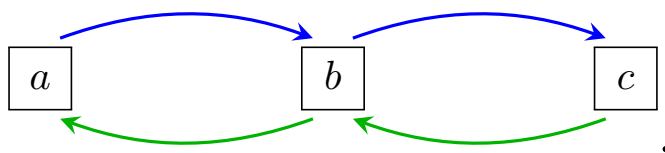
\[
    \begin{tikzpicture}[
      >=stealth,
      node distance=1.5cm,
      tableau/.style={inner sep=0pt, anchor=base}
      ]

      \node[tableau] (T1) {\scalebox{0.7}{\begin{ytableau} a \end{ytableau}}};
      \node[tableau, right=of T1] (T2) {\scalebox{0.7}{\begin{ytableau} b \end{ytableau}}};
      \node[tableau, right=of T2] (T3) {\scalebox{0.7}{\begin{ytableau} c \end{ytableau}}};

      \draw[->, blue, line width=0.6pt, shorten <=4pt, shorten >=4pt]   
      (T1.north) to[out=20,in=160] (T2.north);
      \draw[->, blue, line width=0.6pt, shorten <=4pt, shorten >=4pt]   
      (T2.north) to[out=20,in=160] (T3.north);
      \draw[->, green!70!black, line width=0.6pt, shorten <=4pt, shorten >=4pt] 
      (T3.south) to[out=-160,in=-20] (T2.south);
      \draw[->, green!70!black, line width=0.6pt, shorten <=4pt, shorten >=4pt] 
      (T2.south) to[out=-160,in=-20] (T1.south);

    \end{tikzpicture}.
  \]
  For the definition of perfect crystals, see~\cite[Chapter~10]{HK02}.
  The cell colours in Definitions~\ref{def:P2_type}
  and~\ref{def:P2_crystal} record the incoming arrows of the perfect
  crystal: $a$ has only an incoming green arrow, and $c$ only an
  incoming blue arrow, whereas the arrows $c\to b$ and $a\to b$ are
  green and blue, respectively. The appended $0$s record the outgoing
  blue arrows: two consecutive blue arrows start at $a$, and one
  starts at $b$. Thus the type ends with two additional $0$s when the
  first row ends in an $a$-cell, and with one when it ends in a
  $b$-cell.
\end{remark}

The rules for the arrows give a few useful properties of the vertices
of \( \mathcal{B}_{2,b}^{\mathcal{P}} \). A blue cell is always added
below a green cell, and a green cell below a blue cell, so the colours
alternate along each column. In particular, the entries \( a,c \)
alternate in odd-indexed columns. In an even-indexed column, each cell
has the opposite colour to the cell immediately to its left.
Consequently, the colour of any one cell in columns \( 2i+1 \) and
\( 2i+2 \) determines all the colours in those two columns.

Reading the rows of a vertex from bottom to top, and labelling each
part by the rightmost entry in its row, gives a partition with labels
\( a,b,c \). The next proposition identifies the resulting set with
\( \mathcal{P}_{2,b} \).

\begin{proposition}\label{prop:vertices_P2_crystal}
The set of vertices of the graph $\mathcal{B}_{2,b}^{\mathcal{P}}$ is $\mathcal{P}_{2,b}$.
\end{proposition}
\begin{proof}
  We first prove that every vertex of
  \(\mathcal{B}_{2,b}^{\mathcal{P}}\) is an element of
  \(\mathcal{P}_{2,b}\). By the assumption that all entries in the
  \( 0 \)th column are \( b \), we have that the empty
  partition \( \emptyset \) corresponds to
  \( 0_b \in \mathcal{P}_{2,b} \). Now let \(\pi\) be a nonempty vertex
  of \(\mathcal{B}_{2,b}^{\mathcal{P}}\). We need to show that the labelled partition obtained from \( \pi \) satisfies the exact difference condition in Definition~\ref{def:grounded_partitions}, for the matrix \( M_2 \) in \eqref{matrix:M}.

  We claim that if two consecutive columns of \(\pi\) have the same
  size, then their bottom cells are coloured either both blue or both
  green. Suppose, to the contrary, that the first column ends with a
  green cell and the second with a blue cell. Just before adding the
  last blue cell, both columns end with green cells. In the
  corresponding type, if the letter for the second column is an
  unmatched \( 0 \), then the letter for the first column is also an
  unmatched \( 0 \). Thus the second column cannot be the leftmost
  column with an unmatched \( 0 \), and hence, no cell can be added to
  the second column unless a cell is added to the first column. The
  same argument applies when the first column ends with a blue cell
  and the second with a green cell. This proves the claim.

  By the colouring rule for \( b \)-cells in Definition~\ref{def:P2_type}, there
  cannot be three consecutive cells of the same colour (blue or green)
  in any row of \(\pi\). The bottom cells of three consecutive columns of the same size lie in one row and, by the claim, have the same colour. Therefore, no three consecutive columns of $\pi$ have the same size. In other words, the difference between two consecutive parts of the labelled partition assigned to \(\pi\), including the ground \(0_b\), is at most \(2\).

\begin{itemize}
\item Since entries $a$ and $c$ alternate in columns with odd indices, the difference between two consecutive $a$-parts or two consecutive $c$-parts cannot be $0$. 
\item If the part $(2k)_b$ appears (possibly as $0_b$) and the first part of size $2k+1$ is $(2k+1)_a$, then we cannot add the blue cell $(2k+2)_b$ at the end of $(2k+1)_a$. Indeed, if the letter $0$ of this position were unmatched, then the letter $0$ of the preceding position, at the end of the top row $(2k)_b$, would also be unmatched, and it is further to the left. Similarly, if the part $(2k)_b$ appears (possibly as $0_b$) and the first part of size $2k+1$ is $(2k+1)_c$, then we cannot add the green cell $(2k+2)_b$ at the end of $(2k+1)_c$. Thus, the difference between two consecutive $b$-parts cannot be $2$.
\item If the part of size $2k-1$ appears, followed by a part of size $2k+1$, then the $(2k)$th and $(2k+1)$th columns must have the same colour at the bottom, since they have the same size. Since colours alternate along columns, the bottom cell in the $(2k+1)$th columns must have the same colour as the topmost $(2k-1)$-part. Thus, the difference between a consecutive $a$-part and $c$-part cannot be $2$.
\end{itemize}
  
Hence, in the labelled partition assigned to \( \pi \),
\begin{itemize}
\item the part following \( (2k+1)_a \) is \( (2k+1)_c \), \( (2k+2)_b \), or \( (2k+3)_a \),
\item the part following \( (2k)_b \) is \( (2k)_b \), \( (2k+1)_a \), or \( (2k+1)_c \),
\item the part following \( (2k+1)_c \) is \( (2k+1)_a \), \( (2k+2)_b \), or \( (2k+3)_c \).
\end{itemize}
These conditions coincide with the exact difference condition for
the matrix \( M_2 \). Hence \( \pi \) is an element of
\( \mathcal{P}_{2,b} \).

It remains to show that every element \(\lambda \in \mathcal{P}_{2,b}\) appears as a vertex of \(\mathcal{B}_{2,b}^{\mathcal{P}}\). The proof of the first inclusion shows that arrows map $\mathcal{P}_{2,b}$ into itself.

By Definition~\ref{def:P2_crystal}, a vertex \(\pi\) of \(\mathcal{B}_{2,b}^{\mathcal{P}}\) is the target of a green (respectively, blue) arrow if and only if there is an unmatched right bracket, where we associate a left bracket ``('' to each column of \(\pi\) that has a blue (respectively, green) cell at the bottom and a right bracket ``)'' to each column of \(\pi\) that has a green (respectively, blue) cell at the bottom. The source of the arrow is obtained by removing the bottom cell of the column of the rightmost unmatched right bracket. We extend this definition of inverse arrows to the set \(\mathcal{P}_{2,b}\). Assume \(\pi\in\mathcal{P}_{2,b}\) is the target of a blue arrow with source \(\mu\). We need to show that \(\mu\) also lies in \(\mathcal{P}_{2,b}\). The argument for green arrows is analogous.
\begin{itemize}
\item First, we show that the diagram of \(\mu\) has partition shape, i.e. the columns are weakly decreasing in size. We have seen that there are at most two consecutive columns in \(\pi\) of the same size. If the \(i\)th and \((i+1)\)th columns in \(\pi\) have the same size, then the bottom cells are either both blue or both green. Assume they are both blue. Then, if the letter \(1\) in \(\type(\pi)\) corresponding to the \(i\)th column is unmatched, so is the letter \(1\) corresponding to the \((i+1)\)th column, which is further to the right. In other words, the blue cell at the bottom of the \(i\)th column of \(\pi\) is still present in \(\mu\).
\item Next, we show that \(\mu\) cannot have two consecutive \(b\)-parts with difference \(2\). Assume the \(i\)th column of \(\pi\) has a \(c\)-cell at the bottom and the size of the \((i+1)\)th column is one less. Then the \(b\)-cell at the bottom of the \((i+1)\)th column is also blue, since colours along columns alternate. By the same argument as above, the cell at the bottom of the \(i\)th column of \(\pi\) is still present in \(\mu\). This is illustrated in the first diagram of Figure~\ref{fig:closed_for_inverse_arrow}.
\item Assume that both the \(i\)th and the \((i+2)\)th columns of \(\pi\) have a \(c\)-cell at the bottom; the \((i+1)\)th and the \((i+2)\)th columns have the same size; and the size of the \(i\)th column is one greater. Then the \(b\)-cell at the bottom of the \((i+1)\)th column is blue. By the same argument as above, the cell at the bottom of the \(i\)th column of \(\pi\) is still present in \(\mu\). This is illustrated in the second diagram of Figure~\ref{fig:closed_for_inverse_arrow}.
\begin{figure}[h]
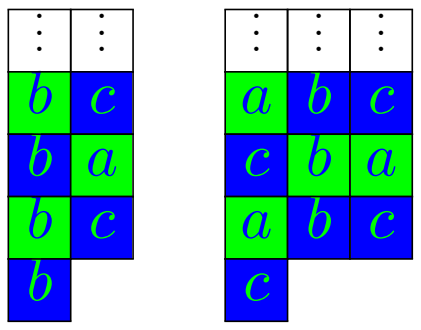

    \centering
    \scalebox{0.7}{%
    \begin{ytableau}
        \vdots & \vdots
        \\
        *(green) \scalebox{1.5}{\color{blue}$b$} & *(blue) \scalebox{1.5}{\color{green}$c$}\\
        *(blue) \scalebox{1.5}{\color{green}$b$} & *(green) \scalebox{1.5}{\color{blue}$a$}\\
        *(green) \scalebox{1.5}{\color{blue}$b$} & *(blue) \scalebox{1.5}{\color{green}$c$}\\
        *(blue) \scalebox{1.5}{\color{green}$b$}
    \end{ytableau}
    }
    \quad
        \scalebox{0.7}{%
    \begin{ytableau}
        \vdots & \vdots & \vdots
        \\
        *(green) \scalebox{1.5}{\color{blue}$a$} & *(blue) \scalebox{1.5}{\color{green}$b$} & *(blue) \scalebox{1.5}{\color{green}$c$}\\
        *(blue) \scalebox{1.5}{\color{green}$c$} & *(green) \scalebox{1.5}{\color{blue}$b$} & *(green) \scalebox{1.5}{\color{blue}$a$}\\
        *(green) \scalebox{1.5}{\color{blue}$a$} & *(blue) \scalebox{1.5}{\color{green}$b$} & *(blue) \scalebox{1.5}{\color{green}$c$}\\
        *(blue) \scalebox{1.5}{\color{green}$c$}
    \end{ytableau}
    }
    \caption{$\mathcal{P}_{2,b}$ is closed for inverse arrows}
    \label{fig:closed_for_inverse_arrow}
\end{figure}
\end{itemize}

By definition, the unique size-\(0\) element \(0_b\) of \(\mathcal{P}_{2,b}\) lies in \(\mathcal{B}_{2,b}^{\mathcal{P}}\). We proceed by induction on the size of \(\pi\). Given \(\pi \in \mathcal{P}_{2,b}\) of size greater than \(0\), we need to show that there is either a green or a blue arrow whose target is \(\pi\).

Assume \(\pi\) is not the target of a green arrow. Then, by definition of the inverse green arrow and because \(\pi\) is nonempty, the bottom cell in the first column is blue. Thus, the first bracket in the sequence is a right bracket, which remains unmatched. This means, by definition of the inverse blue arrow, that \(\pi\) is the target of a blue arrow.

Hence, every element \(\pi \in \mathcal{P}_{2,b}\) appears as a vertex of \(\mathcal{B}_{2,b}^{\mathcal{P}}\).
\end{proof}

\begin{figure}[h]
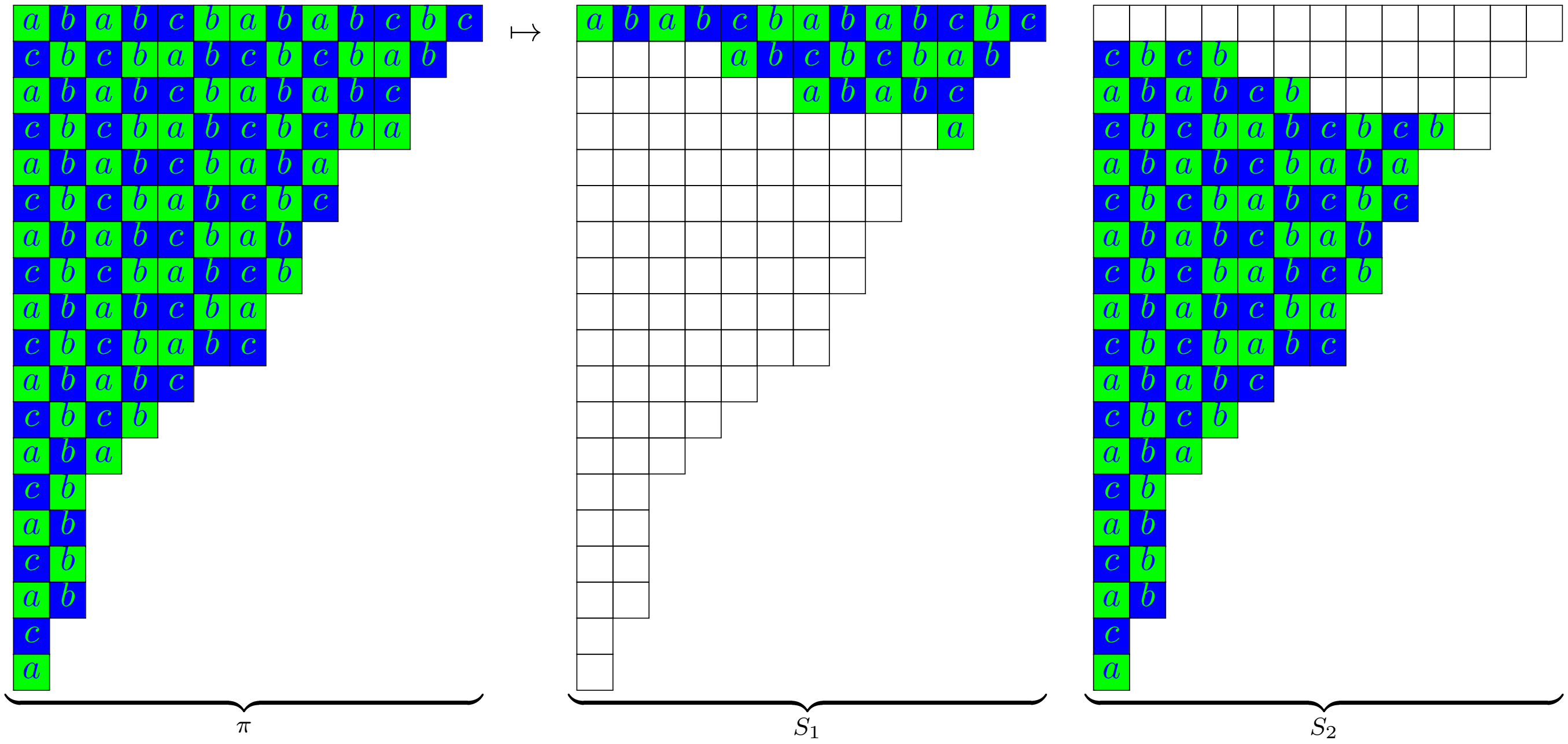

\centering
\[
\underbrace{\scalebox{0.7}{
\begin{ytableau}
*(green) \scalebox{1.5}{\color{blue}$a$}
& *(blue) \scalebox{1.5}{\color{green}$b$}
& *(green) \scalebox{1.5}{\color{blue}$a$}
& *(blue) \scalebox{1.5}{\color{green}$b$}
& *(blue) \scalebox{1.5}{\color{green}$c$}
& *(green) \scalebox{1.5}{\color{blue}$b$}
& *(green) \scalebox{1.5}{\color{blue}$a$}
& *(blue) \scalebox{1.5}{\color{green}$b$}
& *(green) \scalebox{1.5}{\color{blue}$a$}
& *(blue) \scalebox{1.5}{\color{green}$b$}
& *(blue) \scalebox{1.5}{\color{green}$c$}
& *(green) \scalebox{1.5}{\color{blue}$b$}
& *(blue) \scalebox{1.5}{\color{green}$c$} \\
*(blue) \scalebox{1.5}{\color{green}$c$}
& *(green) \scalebox{1.5}{\color{blue}$b$}
& *(blue) \scalebox{1.5}{\color{green}$c$}
& *(green) \scalebox{1.5}{\color{blue}$b$}
& *(green) \scalebox{1.5}{\color{blue}$a$}
& *(blue) \scalebox{1.5}{\color{green}$b$}
& *(blue) \scalebox{1.5}{\color{green}$c$}
& *(green) \scalebox{1.5}{\color{blue}$b$}
& *(blue) \scalebox{1.5}{\color{green}$c$}
& *(green) \scalebox{1.5}{\color{blue}$b$}
& *(green) \scalebox{1.5}{\color{blue}$a$}
& *(blue) \scalebox{1.5}{\color{green}$b$} \\
*(green) \scalebox{1.5}{\color{blue}$a$}
& *(blue) \scalebox{1.5}{\color{green}$b$}
& *(green) \scalebox{1.5}{\color{blue}$a$}
& *(blue) \scalebox{1.5}{\color{green}$b$}
& *(blue) \scalebox{1.5}{\color{green}$c$}
& *(green) \scalebox{1.5}{\color{blue}$b$}
& *(green) \scalebox{1.5}{\color{blue}$a$}
& *(blue) \scalebox{1.5}{\color{green}$b$}
& *(green) \scalebox{1.5}{\color{blue}$a$}
& *(blue) \scalebox{1.5}{\color{green}$b$}
& *(blue) \scalebox{1.5}{\color{green}$c$} \\
*(blue) \scalebox{1.5}{\color{green}$c$}
& *(green) \scalebox{1.5}{\color{blue}$b$}
& *(blue) \scalebox{1.5}{\color{green}$c$}
& *(green) \scalebox{1.5}{\color{blue}$b$}
& *(green) \scalebox{1.5}{\color{blue}$a$}
& *(blue) \scalebox{1.5}{\color{green}$b$}
& *(blue) \scalebox{1.5}{\color{green}$c$}
& *(green) \scalebox{1.5}{\color{blue}$b$}
& *(blue) \scalebox{1.5}{\color{green}$c$}
& *(green) \scalebox{1.5}{\color{blue}$b$}
& *(green) \scalebox{1.5}{\color{blue}$a$} \\
*(green) \scalebox{1.5}{\color{blue}$a$}
& *(blue) \scalebox{1.5}{\color{green}$b$}
& *(green) \scalebox{1.5}{\color{blue}$a$}
& *(blue) \scalebox{1.5}{\color{green}$b$}
& *(blue) \scalebox{1.5}{\color{green}$c$}
& *(green) \scalebox{1.5}{\color{blue}$b$}
& *(green) \scalebox{1.5}{\color{blue}$a$}
& *(blue) \scalebox{1.5}{\color{green}$b$}
& *(green) \scalebox{1.5}{\color{blue}$a$} \\
*(blue) \scalebox{1.5}{\color{green}$c$}
& *(green) \scalebox{1.5}{\color{blue}$b$}
& *(blue) \scalebox{1.5}{\color{green}$c$}
& *(green) \scalebox{1.5}{\color{blue}$b$}
& *(green) \scalebox{1.5}{\color{blue}$a$}
& *(blue) \scalebox{1.5}{\color{green}$b$}
& *(blue) \scalebox{1.5}{\color{green}$c$}
& *(green) \scalebox{1.5}{\color{blue}$b$}
& *(blue) \scalebox{1.5}{\color{green}$c$} \\
*(green) \scalebox{1.5}{\color{blue}$a$}
& *(blue) \scalebox{1.5}{\color{green}$b$}
& *(green) \scalebox{1.5}{\color{blue}$a$}
& *(blue) \scalebox{1.5}{\color{green}$b$}
& *(blue) \scalebox{1.5}{\color{green}$c$}
& *(green) \scalebox{1.5}{\color{blue}$b$}
& *(green) \scalebox{1.5}{\color{blue}$a$}
& *(blue) \scalebox{1.5}{\color{green}$b$} \\
*(blue) \scalebox{1.5}{\color{green}$c$}
& *(green) \scalebox{1.5}{\color{blue}$b$}
& *(blue) \scalebox{1.5}{\color{green}$c$}
& *(green) \scalebox{1.5}{\color{blue}$b$}
& *(green) \scalebox{1.5}{\color{blue}$a$}
& *(blue) \scalebox{1.5}{\color{green}$b$}
& *(blue) \scalebox{1.5}{\color{green}$c$}
& *(green) \scalebox{1.5}{\color{blue}$b$} \\
*(green) \scalebox{1.5}{\color{blue}$a$}
& *(blue) \scalebox{1.5}{\color{green}$b$}
& *(green) \scalebox{1.5}{\color{blue}$a$}
& *(blue) \scalebox{1.5}{\color{green}$b$}
& *(blue) \scalebox{1.5}{\color{green}$c$}
& *(green) \scalebox{1.5}{\color{blue}$b$}
& *(green) \scalebox{1.5}{\color{blue}$a$} \\
*(blue) \scalebox{1.5}{\color{green}$c$}
& *(green) \scalebox{1.5}{\color{blue}$b$}
& *(blue) \scalebox{1.5}{\color{green}$c$}
& *(green) \scalebox{1.5}{\color{blue}$b$}
& *(green) \scalebox{1.5}{\color{blue}$a$}
& *(blue) \scalebox{1.5}{\color{green}$b$}
& *(blue) \scalebox{1.5}{\color{green}$c$} \\
*(green) \scalebox{1.5}{\color{blue}$a$}
& *(blue) \scalebox{1.5}{\color{green}$b$}
& *(green) \scalebox{1.5}{\color{blue}$a$}
& *(blue) \scalebox{1.5}{\color{green}$b$}
& *(blue) \scalebox{1.5}{\color{green}$c$}\\
*(blue) \scalebox{1.5}{\color{green}$c$}
& *(green) \scalebox{1.5}{\color{blue}$b$}
& *(blue) \scalebox{1.5}{\color{green}$c$}
& *(green) \scalebox{1.5}{\color{blue}$b$}\\
*(green) \scalebox{1.5}{\color{blue}$a$}
& *(blue) \scalebox{1.5}{\color{green}$b$}
& *(green) \scalebox{1.5}{\color{blue}$a$} \\
*(blue) \scalebox{1.5}{\color{green}$c$}
& *(green) \scalebox{1.5}{\color{blue}$b$} \\
*(green) \scalebox{1.5}{\color{blue}$a$}
& *(blue) \scalebox{1.5}{\color{green}$b$} \\
*(blue) \scalebox{1.5}{\color{green}$c$}
& *(green) \scalebox{1.5}{\color{blue}$b$} \\
*(green) \scalebox{1.5}{\color{blue}$a$}
& *(blue) \scalebox{1.5}{\color{green}$b$} \\
*(blue) \scalebox{1.5}{\color{green}$c$} \\
*(green) \scalebox{1.5}{\color{blue}$a$}
\end{ytableau}}}_{\pi} \hspace{0.5em}\mapsto \hspace{0.5em}
\underbrace{\scalebox{0.7}{
\begin{ytableau}
*(green) \scalebox{1.5}{\color{blue}$a$}
& *(blue) \scalebox{1.5}{\color{green}$b$}
& *(green) \scalebox{1.5}{\color{blue}$a$}
& *(blue) \scalebox{1.5}{\color{green}$b$}
& *(blue) \scalebox{1.5}{\color{green}$c$}
& *(green) \scalebox{1.5}{\color{blue}$b$}
& *(green) \scalebox{1.5}{\color{blue}$a$}
& *(blue) \scalebox{1.5}{\color{green}$b$}
& *(green) \scalebox{1.5}{\color{blue}$a$}
& *(blue) \scalebox{1.5}{\color{green}$b$}
& *(blue) \scalebox{1.5}{\color{green}$c$}
& *(green) \scalebox{1.5}{\color{blue}$b$}
& *(blue) \scalebox{1.5}{\color{green}$c$} \\
\null
& \null
& \null
& \null
& *(green) \scalebox{1.5}{\color{blue}$a$}
& *(blue) \scalebox{1.5}{\color{green}$b$}
& *(blue) \scalebox{1.5}{\color{green}$c$}
& *(green) \scalebox{1.5}{\color{blue}$b$}
& *(blue) \scalebox{1.5}{\color{green}$c$}
& *(green) \scalebox{1.5}{\color{blue}$b$}
& *(green) \scalebox{1.5}{\color{blue}$a$}
& *(blue) \scalebox{1.5}{\color{green}$b$} \\
\null
& \null
& \null
& \null
& \null
& \null
& *(green) \scalebox{1.5}{\color{blue}$a$}
& *(blue) \scalebox{1.5}{\color{green}$b$}
& *(green) \scalebox{1.5}{\color{blue}$a$}
& *(blue) \scalebox{1.5}{\color{green}$b$}
& *(blue) \scalebox{1.5}{\color{green}$c$} \\
\null
& \null
& \null
& \null
& \null
& \null
& \null
& \null
& \null
& \null
& *(green) \scalebox{1.5}{\color{blue}$a$} \\
\null
& \null
& \null
& \null
& \null
& \null
& \null
& \null
& \null \\
\null
& \null
& \null
& \null
& \null
& \null
& \null
& \null
& \null \\
\null
& \null
& \null
& \null
& \null
& \null
& \null
& \null \\
\null
& \null
& \null
& \null
& \null
& \null
& \null
& \null \\
\null
& \null
& \null
& \null
& \null
& \null
& \null \\
\null
& \null
& \null
& \null
& \null
& \null
& \null \\
\null
& \null
& \null
& \null
& \null \\
\null
& \null
& \null
& \null \\
\null
& \null
& \null \\
\null
& \null \\
\null
& \null \\
\null
& \null \\
\null
& \null \\
\null \\
\null
\end{ytableau}}}_{S_1}\quad
\underbrace{\scalebox{0.7}{
\begin{ytableau}
\null
& \null
& \null
& \null
& \null
& \null
& \null
& \null
& \null
& \null
& \null
& \null
& \null \\
*(blue) \scalebox{1.5}{\color{green}$c$}
& *(green) \scalebox{1.5}{\color{blue}$b$}
& *(blue) \scalebox{1.5}{\color{green}$c$}
& *(green) \scalebox{1.5}{\color{blue}$b$}
& \null
& \null
& \null
& \null
& \null
& \null
& \null
& \null \\
*(green) \scalebox{1.5}{\color{blue}$a$}
& *(blue) \scalebox{1.5}{\color{green}$b$}
& *(green) \scalebox{1.5}{\color{blue}$a$}
& *(blue) \scalebox{1.5}{\color{green}$b$}
& *(blue) \scalebox{1.5}{\color{green}$c$}
& *(green) \scalebox{1.5}{\color{blue}$b$}
& \null
& \null
& \null
& \null
& \null \\
*(blue) \scalebox{1.5}{\color{green}$c$}
& *(green) \scalebox{1.5}{\color{blue}$b$}
& *(blue) \scalebox{1.5}{\color{green}$c$}
& *(green) \scalebox{1.5}{\color{blue}$b$}
& *(green) \scalebox{1.5}{\color{blue}$a$}
& *(blue) \scalebox{1.5}{\color{green}$b$}
& *(blue) \scalebox{1.5}{\color{green}$c$}
& *(green) \scalebox{1.5}{\color{blue}$b$}
& *(blue) \scalebox{1.5}{\color{green}$c$}
& *(green) \scalebox{1.5}{\color{blue}$b$}
& \null \\
*(green) \scalebox{1.5}{\color{blue}$a$}
& *(blue) \scalebox{1.5}{\color{green}$b$}
& *(green) \scalebox{1.5}{\color{blue}$a$}
& *(blue) \scalebox{1.5}{\color{green}$b$}
& *(blue) \scalebox{1.5}{\color{green}$c$}
& *(green) \scalebox{1.5}{\color{blue}$b$}
& *(green) \scalebox{1.5}{\color{blue}$a$}
& *(blue) \scalebox{1.5}{\color{green}$b$}
& *(green) \scalebox{1.5}{\color{blue}$a$} \\
*(blue) \scalebox{1.5}{\color{green}$c$}
& *(green) \scalebox{1.5}{\color{blue}$b$}
& *(blue) \scalebox{1.5}{\color{green}$c$}
& *(green) \scalebox{1.5}{\color{blue}$b$}
& *(green) \scalebox{1.5}{\color{blue}$a$}
& *(blue) \scalebox{1.5}{\color{green}$b$}
& *(blue) \scalebox{1.5}{\color{green}$c$}
& *(green) \scalebox{1.5}{\color{blue}$b$}
& *(blue) \scalebox{1.5}{\color{green}$c$} \\
*(green) \scalebox{1.5}{\color{blue}$a$}
& *(blue) \scalebox{1.5}{\color{green}$b$}
& *(green) \scalebox{1.5}{\color{blue}$a$}
& *(blue) \scalebox{1.5}{\color{green}$b$}
& *(blue) \scalebox{1.5}{\color{green}$c$}
& *(green) \scalebox{1.5}{\color{blue}$b$}
& *(green) \scalebox{1.5}{\color{blue}$a$}
& *(blue) \scalebox{1.5}{\color{green}$b$} \\
*(blue) \scalebox{1.5}{\color{green}$c$}
& *(green) \scalebox{1.5}{\color{blue}$b$}
& *(blue) \scalebox{1.5}{\color{green}$c$}
& *(green) \scalebox{1.5}{\color{blue}$b$}
& *(green) \scalebox{1.5}{\color{blue}$a$}
& *(blue) \scalebox{1.5}{\color{green}$b$}
& *(blue) \scalebox{1.5}{\color{green}$c$}
& *(green) \scalebox{1.5}{\color{blue}$b$} \\
*(green) \scalebox{1.5}{\color{blue}$a$}
& *(blue) \scalebox{1.5}{\color{green}$b$}
& *(green) \scalebox{1.5}{\color{blue}$a$}
& *(blue) \scalebox{1.5}{\color{green}$b$}
& *(blue) \scalebox{1.5}{\color{green}$c$}
& *(green) \scalebox{1.5}{\color{blue}$b$}
& *(green) \scalebox{1.5}{\color{blue}$a$} \\
*(blue) \scalebox{1.5}{\color{green}$c$}
& *(green) \scalebox{1.5}{\color{blue}$b$}
& *(blue) \scalebox{1.5}{\color{green}$c$}
& *(green) \scalebox{1.5}{\color{blue}$b$}
& *(green) \scalebox{1.5}{\color{blue}$a$}
& *(blue) \scalebox{1.5}{\color{green}$b$}
& *(blue) \scalebox{1.5}{\color{green}$c$} \\
*(green) \scalebox{1.5}{\color{blue}$a$}
& *(blue) \scalebox{1.5}{\color{green}$b$}
& *(green) \scalebox{1.5}{\color{blue}$a$}
& *(blue) \scalebox{1.5}{\color{green}$b$}
& *(blue) \scalebox{1.5}{\color{green}$c$}\\
*(blue) \scalebox{1.5}{\color{green}$c$}
& *(green) \scalebox{1.5}{\color{blue}$b$}
& *(blue) \scalebox{1.5}{\color{green}$c$}
& *(green) \scalebox{1.5}{\color{blue}$b$}\\
*(green) \scalebox{1.5}{\color{blue}$a$}
& *(blue) \scalebox{1.5}{\color{green}$b$}
& *(green) \scalebox{1.5}{\color{blue}$a$} \\
*(blue) \scalebox{1.5}{\color{green}$c$}
& *(green) \scalebox{1.5}{\color{blue}$b$} \\
*(green) \scalebox{1.5}{\color{blue}$a$}
& *(blue) \scalebox{1.5}{\color{green}$b$} \\
*(blue) \scalebox{1.5}{\color{green}$c$}
& *(green) \scalebox{1.5}{\color{blue}$b$} \\
*(green) \scalebox{1.5}{\color{blue}$a$}
& *(blue) \scalebox{1.5}{\color{green}$b$} \\
*(blue) \scalebox{1.5}{\color{green}$c$} \\
*(green) \scalebox{1.5}{\color{blue}$a$}
\end{ytableau}}}_{S_2}
\]

\caption{Decomposing $\pi$ into $S_1$ and $S_2$}
\label{fig:grounded_S_1}
\end{figure}

We now define a map \(\delta\) from \(\mathcal{P}_{2,b}\) to pairs of partitions by decomposing the cells of \(\pi\in\mathcal{P}_{2,b}\) into disjoint sets \(S_1\) and \(S_2\). If the first row contains no \(a\)-cell, set \(S_1=\emptyset\) and put every cell in \(S_2\). Otherwise, start at the first \(a\)-cell in the first row and form a horizontal strip of alternating \(a\)-cells and blue \(b\)-cells. At each step, the next cell lies either immediately to the right or immediately below and to the right; continue as long as possible. Put in \(S_1\) all cells above or on this strip, together with all cells to the right of its last cell that lie weakly above that cell's row. Put the remaining cells in \(S_2\). Figure~\ref{fig:grounded_S_1} illustrates this decomposition.

\begin{figure}[h]
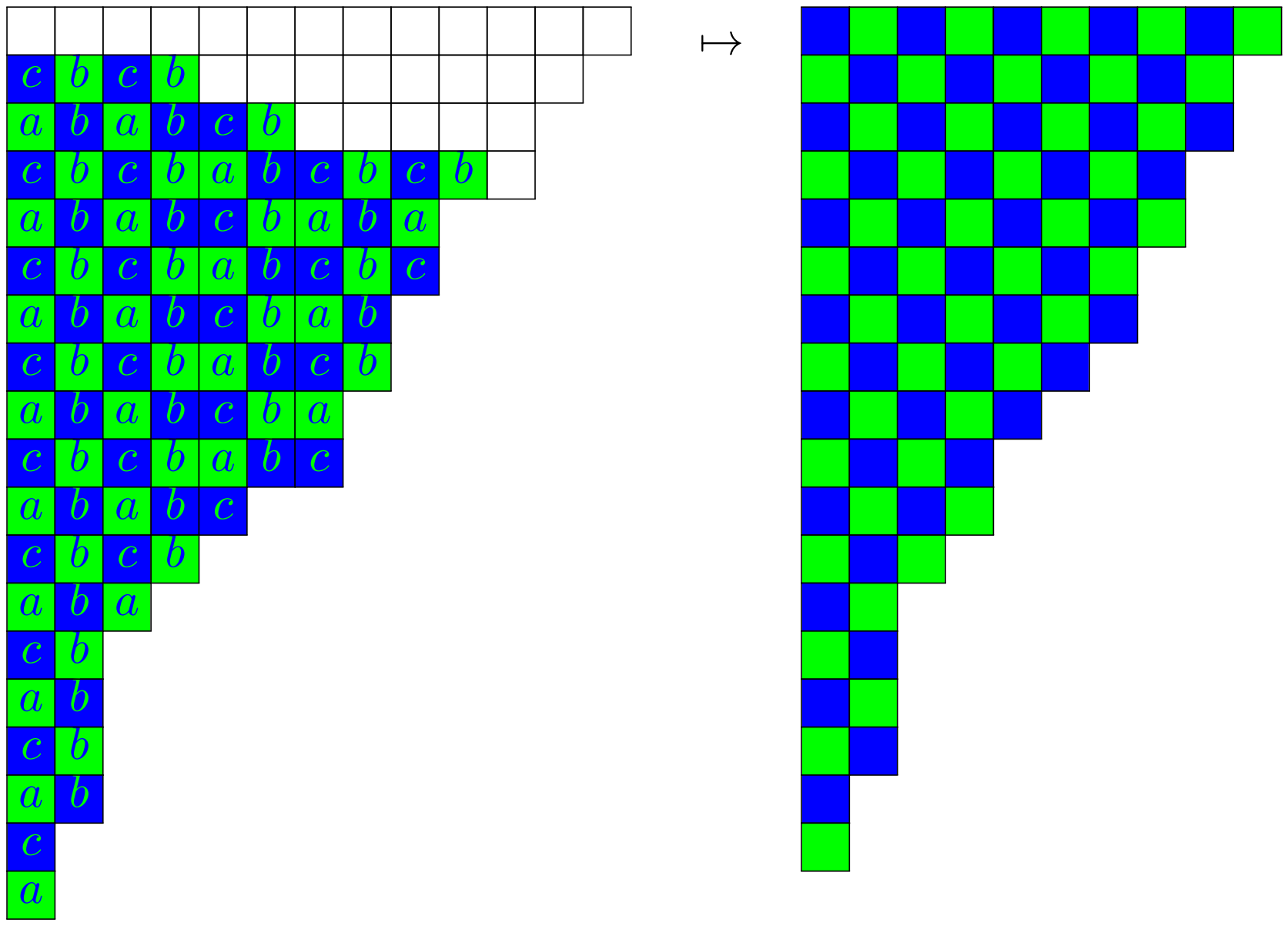

\centering
\[
\scalebox{0.7}{
\begin{ytableau}
\null
& \null
& \null
& \null
& \null
& \null
& \null
& \null
& \null
& \null
& \null
& \null
& \null \\
*(blue) \scalebox{1.5}{\color{green}$c$}
& *(green) \scalebox{1.5}{\color{blue}$b$}
& *(blue) \scalebox{1.5}{\color{green}$c$}
& *(green) \scalebox{1.5}{\color{blue}$b$}
& \null
& \null
& \null
& \null
& \null
& \null
& \null
& \null \\
*(green) \scalebox{1.5}{\color{blue}$a$}
& *(blue) \scalebox{1.5}{\color{green}$b$}
& *(green) \scalebox{1.5}{\color{blue}$a$}
& *(blue) \scalebox{1.5}{\color{green}$b$}
& *(blue) \scalebox{1.5}{\color{green}$c$}
& *(green) \scalebox{1.5}{\color{blue}$b$}
& \null
& \null
& \null
& \null
& \null \\
*(blue) \scalebox{1.5}{\color{green}$c$}
& *(green) \scalebox{1.5}{\color{blue}$b$}
& *(blue) \scalebox{1.5}{\color{green}$c$}
& *(green) \scalebox{1.5}{\color{blue}$b$}
& *(green) \scalebox{1.5}{\color{blue}$a$}
& *(blue) \scalebox{1.5}{\color{green}$b$}
& *(blue) \scalebox{1.5}{\color{green}$c$}
& *(green) \scalebox{1.5}{\color{blue}$b$}
& *(blue) \scalebox{1.5}{\color{green}$c$}
& *(green) \scalebox{1.5}{\color{blue}$b$}
& \null \\
*(green) \scalebox{1.5}{\color{blue}$a$}
& *(blue) \scalebox{1.5}{\color{green}$b$}
& *(green) \scalebox{1.5}{\color{blue}$a$}
& *(blue) \scalebox{1.5}{\color{green}$b$}
& *(blue) \scalebox{1.5}{\color{green}$c$}
& *(green) \scalebox{1.5}{\color{blue}$b$}
& *(green) \scalebox{1.5}{\color{blue}$a$}
& *(blue) \scalebox{1.5}{\color{green}$b$}
& *(green) \scalebox{1.5}{\color{blue}$a$} \\
*(blue) \scalebox{1.5}{\color{green}$c$}
& *(green) \scalebox{1.5}{\color{blue}$b$}
& *(blue) \scalebox{1.5}{\color{green}$c$}
& *(green) \scalebox{1.5}{\color{blue}$b$}
& *(green) \scalebox{1.5}{\color{blue}$a$}
& *(blue) \scalebox{1.5}{\color{green}$b$}
& *(blue) \scalebox{1.5}{\color{green}$c$}
& *(green) \scalebox{1.5}{\color{blue}$b$}
& *(blue) \scalebox{1.5}{\color{green}$c$} \\
*(green) \scalebox{1.5}{\color{blue}$a$}
& *(blue) \scalebox{1.5}{\color{green}$b$}
& *(green) \scalebox{1.5}{\color{blue}$a$}
& *(blue) \scalebox{1.5}{\color{green}$b$}
& *(blue) \scalebox{1.5}{\color{green}$c$}
& *(green) \scalebox{1.5}{\color{blue}$b$}
& *(green) \scalebox{1.5}{\color{blue}$a$}
& *(blue) \scalebox{1.5}{\color{green}$b$} \\
*(blue) \scalebox{1.5}{\color{green}$c$}
& *(green) \scalebox{1.5}{\color{blue}$b$}
& *(blue) \scalebox{1.5}{\color{green}$c$}
& *(green) \scalebox{1.5}{\color{blue}$b$}
& *(green) \scalebox{1.5}{\color{blue}$a$}
& *(blue) \scalebox{1.5}{\color{green}$b$}
& *(blue) \scalebox{1.5}{\color{green}$c$}
& *(green) \scalebox{1.5}{\color{blue}$b$} \\
*(green) \scalebox{1.5}{\color{blue}$a$}
& *(blue) \scalebox{1.5}{\color{green}$b$}
& *(green) \scalebox{1.5}{\color{blue}$a$}
& *(blue) \scalebox{1.5}{\color{green}$b$}
& *(blue) \scalebox{1.5}{\color{green}$c$}
& *(green) \scalebox{1.5}{\color{blue}$b$}
& *(green) \scalebox{1.5}{\color{blue}$a$} \\
*(blue) \scalebox{1.5}{\color{green}$c$}
& *(green) \scalebox{1.5}{\color{blue}$b$}
& *(blue) \scalebox{1.5}{\color{green}$c$}
& *(green) \scalebox{1.5}{\color{blue}$b$}
& *(green) \scalebox{1.5}{\color{blue}$a$}
& *(blue) \scalebox{1.5}{\color{green}$b$}
& *(blue) \scalebox{1.5}{\color{green}$c$} \\
*(green) \scalebox{1.5}{\color{blue}$a$}
& *(blue) \scalebox{1.5}{\color{green}$b$}
& *(green) \scalebox{1.5}{\color{blue}$a$}
& *(blue) \scalebox{1.5}{\color{green}$b$}
& *(blue) \scalebox{1.5}{\color{green}$c$}\\
*(blue) \scalebox{1.5}{\color{green}$c$}
& *(green) \scalebox{1.5}{\color{blue}$b$}
& *(blue) \scalebox{1.5}{\color{green}$c$}
& *(green) \scalebox{1.5}{\color{blue}$b$}\\
*(green) \scalebox{1.5}{\color{blue}$a$}
& *(blue) \scalebox{1.5}{\color{green}$b$}
& *(green) \scalebox{1.5}{\color{blue}$a$} \\
*(blue) \scalebox{1.5}{\color{green}$c$}
& *(green) \scalebox{1.5}{\color{blue}$b$} \\
*(green) \scalebox{1.5}{\color{blue}$a$}
& *(blue) \scalebox{1.5}{\color{green}$b$} \\
*(blue) \scalebox{1.5}{\color{green}$c$}
& *(green) \scalebox{1.5}{\color{blue}$b$} \\
*(green) \scalebox{1.5}{\color{blue}$a$}
& *(blue) \scalebox{1.5}{\color{green}$b$} \\
*(blue) \scalebox{1.5}{\color{green}$c$} \\
*(green) \scalebox{1.5}{\color{blue}$a$}
\end{ytableau}
}
\hspace{1em}\mapsto \hspace{1em}
\scalebox{0.7}{\begin{ytableau}
*(blue) & *(green) & *(blue) & *(green) & *(blue) & *(green) & *(blue) & *(green) & *(blue) & *(green) \\
*(green) & *(blue) & *(green) & *(blue) & *(green) & *(blue) & *(green) & *(blue) & *(green) \\
*(blue) & *(green) & *(blue) & *(green) & *(blue) & *(green) & *(blue) & *(green) & *(blue) \\
*(green) & *(blue) & *(green) & *(blue) & *(green) & *(blue) & *(green) & *(blue) \\
*(blue) & *(green) & *(blue) & *(green) & *(blue) & *(green) & *(blue) & *(green) \\
*(green) & *(blue) & *(green) & *(blue) & *(green) & *(blue) & *(green) \\
*(blue) & *(green) & *(blue) & *(green) & *(blue) & *(green) & *(blue) \\
*(green) & *(blue) & *(green) & *(blue) & *(green)& *(blue) \\
*(blue) & *(green) & *(blue) & *(green) & *(blue) \\
*(green) & *(blue) & *(green) & *(blue)  \\
*(blue) & *(green) & *(blue) & *(green) \\
*(green) & *(blue) & *(green) \\
*(blue) & *(green)  \\
*(green) & *(blue) \\
*(blue) & *(green)  \\
*(green) & *(blue) \\
*(blue) \\
*(green)
\end{ytableau}}
\]

\caption{Obtaining $\lambda^{(2)}$ from $S_2$}
\label{fig:lambda2_from_S2}
\end{figure}

The partition \( \lambda^{(2)} \) is obtained from \( S_2 \) by sliding its columns upward and forgetting the entries \(a, b, c\). This is illustrated in Figure~\ref{fig:lambda2_from_S2}.

\begin{figure}[h]
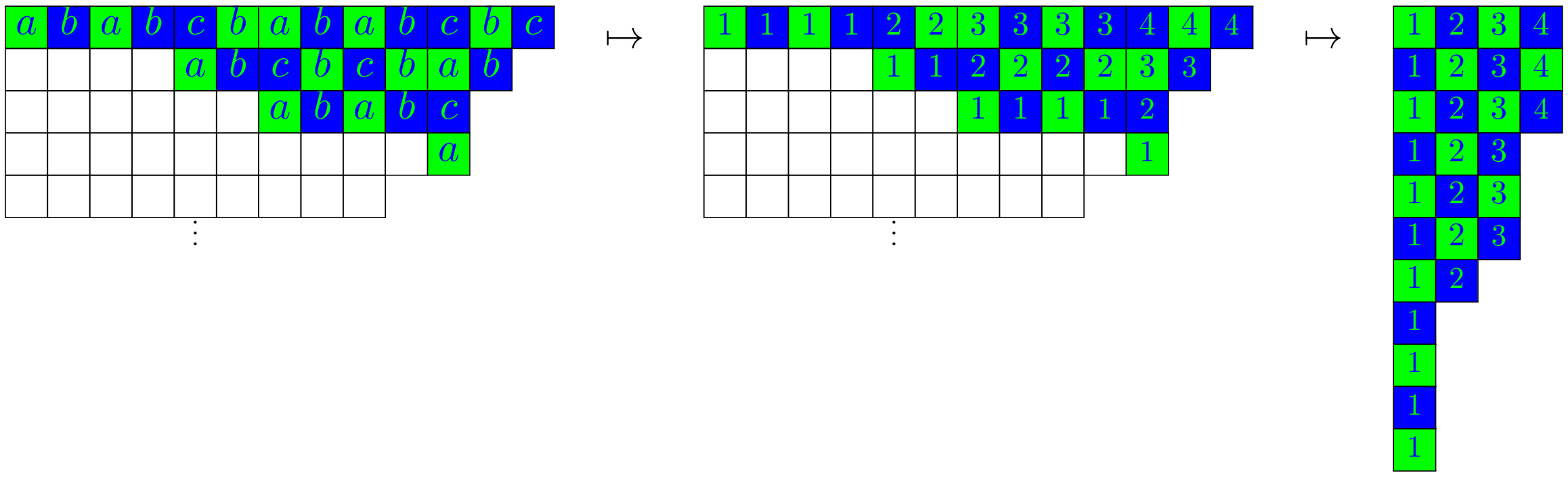

\centering
\[
\scalebox{0.7}{
\begin{ytableau}
*(green) \scalebox{1.5}{\color{blue}$a$}
& *(blue) \scalebox{1.5}{\color{green}$b$}
& *(green) \scalebox{1.5}{\color{blue}$a$}
& *(blue) \scalebox{1.5}{\color{green}$b$}
& *(blue) \scalebox{1.5}{\color{green}$c$}
& *(green) \scalebox{1.5}{\color{blue}$b$}
& *(green) \scalebox{1.5}{\color{blue}$a$}
& *(blue) \scalebox{1.5}{\color{green}$b$}
& *(green) \scalebox{1.5}{\color{blue}$a$}
& *(blue) \scalebox{1.5}{\color{green}$b$}
& *(blue) \scalebox{1.5}{\color{green}$c$}
& *(green) \scalebox{1.5}{\color{blue}$b$}
& *(blue) \scalebox{1.5}{\color{green}$c$} \\
\null
& \null
& \null
& \null
& *(green) \scalebox{1.5}{\color{blue}$a$}
& *(blue) \scalebox{1.5}{\color{green}$b$}
& *(blue) \scalebox{1.5}{\color{green}$c$}
& *(green) \scalebox{1.5}{\color{blue}$b$}
& *(blue) \scalebox{1.5}{\color{green}$c$}
& *(green) \scalebox{1.5}{\color{blue}$b$}
& *(green) \scalebox{1.5}{\color{blue}$a$}
& *(blue) \scalebox{1.5}{\color{green}$b$} \\
\null
& \null
& \null
& \null
& \null
& \null
& *(green) \scalebox{1.5}{\color{blue}$a$}
& *(blue) \scalebox{1.5}{\color{green}$b$}
& *(green) \scalebox{1.5}{\color{blue}$a$}
& *(blue) \scalebox{1.5}{\color{green}$b$}
& *(blue) \scalebox{1.5}{\color{green}$c$} \\
\null
& \null
& \null
& \null
& \null
& \null
& \null
& \null
& \null
& \null
& *(green) \scalebox{1.5}{\color{blue}$a$} \\
\null
& \null
& \null
& \null
& \null
& \null
& \null
& \null
& \null \\
\none & \none & \none & \none &   \none[\vdots]
\end{ytableau}}
\hspace{1em}\mapsto \hspace{1em}
\scalebox{0.7}{
\begin{ytableau}
*(green) \scalebox{1.2}{\color{blue}$1$}
& *(blue) \scalebox{1.2}{\color{green}$1$}
& *(green) \scalebox{1.2}{\color{blue}$1$}
& *(blue) \scalebox{1.2}{\color{green}$1$}
& *(blue) \scalebox{1.2}{\color{green}$2$}
& *(green) \scalebox{1.2}{\color{blue}$2$}
& *(green) \scalebox{1.2}{\color{blue}$3$}
& *(blue) \scalebox{1.2}{\color{green}$3$}
& *(green) \scalebox{1.2}{\color{blue}$3$}
& *(blue) \scalebox{1.2}{\color{green}$3$}
& *(blue) \scalebox{1.2}{\color{green}$4$}
& *(green) \scalebox{1.2}{\color{blue}$4$}
& *(blue) \scalebox{1.1}{\color{green}$4$} \\
\null
& \null
& \null
& \null
& *(green) \scalebox{1.2}{\color{blue}$1$}
& *(blue) \scalebox{1.2}{\color{green}$1$}
& *(blue) \scalebox{1.2}{\color{green}$2$}
& *(green) \scalebox{1.2}{\color{blue}$2$}
& *(blue) \scalebox{1.2}{\color{green}$2$}
& *(green) \scalebox{1.2}{\color{blue}$2$}
& *(green) \scalebox{1.2}{\color{blue}$3$}
& *(blue) \scalebox{1.1}{\color{green}$3$} \\
\null
& \null
& \null
& \null
& \null
& \null
& *(green) \scalebox{1.2}{\color{blue}$1$}
& *(blue) \scalebox{1.2}{\color{green}$1$}
& *(green) \scalebox{1.2}{\color{blue}$1$}
& *(blue) \scalebox{1.1}{\color{green}$1$}
& *(blue) \scalebox{1.1}{\color{green}$2$} \\
\null
& \null
& \null
& \null
& \null
& \null
& \null
& \null
& \null
& \null
& *(green) \scalebox{1.1}{\color{blue}$1$} \\
\null
& \null
& \null
& \null
& \null
& \null
& \null
& \null
& \null \\
\none & \none & \none & \none &   \none[\vdots]
\end{ytableau}}
\hspace{1em}\mapsto \hspace{1em}
\scalebox{0.7}{\begin{ytableau}
*(green) \scalebox{1.2}{\color{blue}$1$} & *(blue) \scalebox{1.2}{\color{green}$2$} & *(green) \scalebox{1.2}{\color{blue}$3$}& *(blue) \scalebox{1.2}{\color{green}$4$} \\
*(blue) \scalebox{1.2}{\color{green}$1$} & *(green) \scalebox{1.2}{\color{blue}$2$} & *(blue) \scalebox{1.2}{\color{green}$3$} & *(green) \scalebox{1.2}{\color{blue}$4$} \\
*(green) \scalebox{1.2}{\color{blue}$1$} & *(blue) \scalebox{1.2}{\color{green}$2$} & *(green) \scalebox{1.2}{\color{blue}$3$} & *(blue) \scalebox{1.1}{\color{green}$4$} \\
*(blue) \scalebox{1.2}{\color{green}$1$} & *(green) \scalebox{1.2}{\color{blue}$2$} & *(blue) \scalebox{1.2}{\color{green}$3$} \\
*(green) \scalebox{1.2}{\color{blue}$1$} & *(blue) \scalebox{1.2}{\color{green}$2$} & *(green) \scalebox{1.2}{\color{blue}$3$} \\
*(blue) \scalebox{1.2}{\color{green}$1$} & *(green) \scalebox{1.2}{\color{blue}$2$} & *(blue) \scalebox{1.1}{\color{green}$3$} \\
*(green) \scalebox{1.2}{\color{blue}$1$} & *(blue) \scalebox{1.1}{\color{green}$2$}\\
*(blue) \scalebox{1.1}{\color{green}$1$} \\
*(green) \scalebox{1.1}{\color{blue}$1$} \\
*(blue) \scalebox{1.1}{\color{green}$1$} \\
*(green) \scalebox{1.1}{\color{blue}$1$}
\end{ytableau}}
\]
\caption{Obtaining \( \lambda^{(1)} \) from $S_1$}
\label{fig:lambda1_from_S1}
\end{figure}

The partition $\lambda^{(1)}$ is constructed from $S_1$ as follows. We decompose $S_1$ into horizontal strips. The first strip is the strip constructed above. The odd-indexed strips have $a$- and blue $b$-cells alternatingly, and the even-indexed strips have $c$- and green $b$-cells alternatingly. Let $h_1, h_2, \ldots$ denote the lengths of these horizontal strips. Then $\lambda^{(1)}$ is defined as the partition whose $i$th column has length $h_i$. In our example, we label the cells of these strips with integers $1,\ldots,4$, as illustrated in the middle part of Figure~\ref{fig:lambda1_from_S1}. These horizontal strips have lengths $(h_1,h_2,h_3,h_4)=(11,7,6,3)$. Then we have $\lambda^{(1)}=(4,4,4,3,3,3,2,1,1,1,1)$.

\begin{lemma}\label{lemma:delta_well-def}
Let $\pi\in\mathcal{P}_{2,b}$ and $\delta(\pi)=(\lambda^{(1)},\lambda^{(2)})$. Then $\lambda^{(1)}$ and $\lambda^{(2)}$ are $2$-restricted partitions.
\end{lemma}
\begin{proof}
For \(\lambda^{(1)}\), consider consecutive strips of lengths \(h_i\) and \(h_{i+1}\). Each strip contains at most one cell in each column. By construction, at least two columns meet the \(i\)th strip but not the \((i+1)\)th strip. The difference conditions for \( \mathcal{P}_{2,b} \) allow at most one column meeting the \((i+1)\)th strip but not the \(i\)th strip. Hence \( h_i > h_{i+1} \), so the column lengths of \( \lambda^{(1)} \) are distinct and \( \lambda^{(1)} \) is $2$-restricted.

Now we consider $\lambda^{(2)}$. By the difference conditions of $\mathcal{P}_{2,b}$, there are at most two consecutive columns in $\pi$ of the same length. We have seen in the proof of Proposition~\ref{prop:vertices_P2_crystal} that cells in two consecutive columns of the same length have an alternating green-blue domino pattern. By construction, the horizontal strip at the south wall of $S_1$ does not contain two consecutive cells of the same colour. Thus, if two consecutive columns in $\pi$ have the same length, then their intersection with $S_2$ cannot have the same length. In other words, all columns of $\lambda^{(2)}$ are distinct.
\end{proof}

We say that two binary words \(r\) and \(s\) on the alphabet $\{0,1\}$ are \emph{equivalent},
written \(r\sim s\), if the subwords of $r$ and $s$ restricted to
unmatched letters are equal. We then show that \(\delta\) preserves \(\type\) up to this equivalence.

\begin{lemma}\label{lemma:delta_type}
We have $\type(\pi)\sim\type(\delta(\pi))$ for all $\pi \in \mathcal{P}_{2,b}$.
\end{lemma}
\begin{proof}
Write \(\delta(\pi)=(\lambda^{(1)},\lambda^{(2)})\), with \(S_1\) and \(S_2\) as in the construction of \(\delta\). If \(S_1=\emptyset\), then \(\delta(\pi)=(\emptyset,\pi)\) and the types agree. We may therefore assume that \(S_1\neq\emptyset\).

  If \(S_2=\emptyset\), the first column consists only of the initial \(a\)-cell of the strip defining \(S_1\). Thus \(\pi=1_a\). In this case \(\delta(\pi)=((1),\emptyset)\), and both types are \(000\). We may therefore also assume that \(S_2\neq\emptyset\).

  Let \(\type(S_2)\) be the factor of \(\type(\pi)\) formed by the letters corresponding to columns whose bottom cells lie in \(S_2\). Then \(\type(\lambda^{(2)})=\type(S_2)\) if the first row of \(\lambda^{(2)}\) ends in blue, and \(\type(\lambda^{(2)})=\type(S_2)0\) if it ends in green.

  Let \(w(\pi)\) and \(w(S_2)\) denote the numbers of columns containing cells of \(\pi\) and \(S_2\), respectively. We induct on \(w(\pi)-w(S_2)\geq0\). In the base case these widths are equal, so the last column of \(\pi\) intersects \(S_2\). We distinguish four cases according to the entry and colour of its top and bottom cells within \(S_1\). By construction, the bottom cell of \(S_1\) in a column that continues into \(S_2\) is either an \(a\)-cell or a blue \(b\)-cell.

\begin{enumerate}[label=(1\alph*)]
\item Assume that the last column of \(\pi\) has an \(a\)-cell both at
  the top and the bottom of \(S_1\). Since colours alternate along
  columns and \(a\)-cells are green, the first row of \(\lambda^{(2)}\)
  ends in a blue cell. By definition of \(\delta\), each strip starts
  with an \(a\)- or \(c\)-cell and no strip contains both. Thus, the
  topmost strip of \(\pi\) starts with an \(a\)-cell, and the first
  row of \(\lambda^{(1)}\) ends in green. Since the last column of \(\pi\)
  intersects \(S_2\), the strip endings and the cells of the last
  column within \(S_1\) are in bijection. It follows that
\[
\type(\lambda^{(1)},\lambda^{(2)})=\type(\lambda^{(2)})\type(\lambda^{(1)})=\type(S_2)\cdot0101\cdots010\,0\sim\type(S_2)00.
\]
Since the first row of \(\pi\) ends in an \(a\)-cell, we have
\(\type(\pi)=\type(S_2)00\). Thus,
\(\type(\pi)\sim\type(\lambda^{(1)},\lambda^{(2)})\).

\item Assume that the last column of \(\pi\) has a \(c\)-cell at the
  top and an \(a\)-cell at the bottom of \(S_1\). Then, the first row
  of \(\lambda^{(2)}\) ends in a blue cell. Since the topmost strip in
  \(\pi\) must start with a \(c\)-cell, the first row of \(\lambda^{(1)}\)
  ends in blue. Hence
  \[
    \type(\lambda^{(1)},\lambda^{(2)})=\type(S_2)\cdot
    0101\cdots01\sim\type(S_2).
  \]
  As the first row of \(\pi\) ends in a \(c\)-cell, we have
  \(\type(\pi)=\type(S_2)\). Therefore,
  \(\type(\pi)\sim\type(\lambda^{(1)},\lambda^{(2)})\).

\item Assume that the last column of \(\pi\) has a blue \(b\)-cell
  both at the top and the bottom of \(S_1\). Then, the first row of
  \(\lambda^{(2)}\) ends in a green cell. In any strip that contains a
  blue \(b\)-cell, \(a\)-cells and blue \(b\)-cells alternate. Hence
  the topmost strip in \(\pi\) must start with an \(a\)-cell, and the
  first row of \(\lambda^{(1)}\) ends in green. It follows that
  \[
    \type(\lambda^{(1)},\lambda^{(2)})=\type(S_2)0\cdot 1010\cdots101\,0\sim\type(S_2)0.
  \]
  Since the first row of \(\pi\) ends in a \(b\)-cell, we have
  \(\type(\pi)=\type(S_2)0\). Therefore,
  \(\type(\pi)\sim\type(\lambda^{(1)},\lambda^{(2)})\).

\item Assume that the last column of \(\pi\) has a green \(b\)-cell at
  the top and a blue \(b\)-cell at the bottom of \(S_1\). Then, the
  first row of \(\lambda^{(2)}\) ends in a green cell. In any strip that
  contains a green \(b\)-cell, \(c\)-cells and green \(b\)-cells
  alternate. Hence the topmost strip in \(\pi\) must start with a
  \(c\)-cell, and the first row of \(\lambda^{(1)}\) ends in blue. It
  follows that
  \[
    \type(\lambda^{(1)},\lambda^{(2)})=\type(S_2)0\cdot 1010\cdots10\sim\type(S_2)0.
  \]
  Since the first row of \(\pi\) ends in a \(b\)-cell, we have
  \(\type(\pi)=\type(S_2)0\). Therefore,
  \(\type(\pi)\sim\type(\lambda^{(1)},\lambda^{(2)})\).
\end{enumerate}

For the induction step, assume \(w(\pi)-w(S_2)>0\), and let \(\pi_\circ\) be the grounded partition obtained by deleting the last column of \(\pi\). Write \(\delta(\pi_\circ)=(\lambda_\circ^{(1)},\lambda_\circ^{(2)})\). By the case \(S_1=\emptyset\) or by induction, \(\type(\pi_\circ)\sim\type(\lambda_\circ^{(1)},\lambda_\circ^{(2)})\). Let \(i+1\) be the number of columns of \(\pi\), and denote the top and bottom cells of column \(j\) by \(T_j\) and \(B_j\), respectively. We distinguish twelve cases according to the entries and colours of \(T_i\), \(T_{i+1}\), \(B_{i+1}\).

\begin{enumerate}[label=(2\alph*)]

\item Assume that \( T_i \), \( T_{i+1} \), and \( B_{i+1} \) are a green \( b \)-cell, an \( a \)-cell, and an \( a \)-cell, respectively, illustrated in Figure~\ref{fig:delta_iso_2abcd}~(2a). Then the cell to the
  left of \( B_{i+1} \) is a green \( b \)-cell. We have
  \[
    \type(\pi_\circ)=s0, \qand \type(\pi)=s000,
  \]
  for some binary word \(s\).

  An even number of strips descend one row lower from the \(i\)th
  column to the \((i+1)\)th column; and a new strip starts in the top
  cell of the \((i+1)\)th column, which is a green cell. All strips
  continuing from \(\pi_\circ\) to \(\pi\) change the colour of the
  strip ending cell from green to blue, or from blue to green. All
  other strip endings in \(\pi_\circ\) are unchanged in \(\pi\). Note
  that the topmost strip in \(\pi\) starts with an \(a\)-cell, so the
  first row of \(\lambda^{(1)}\) ends in green; the topmost strip in
  \(\pi_\circ\), if \(\pi_\circ\) has a strip, starts with a
  \(c\)-cell, so \(\lambda_\circ^{(1)}\) is empty or its first row
  ends in blue. Moreover, if the green \(b\)-cell to the left of
  \( B_{i+1} \) belongs to \(S_1\), then the topmost strip ending that
  remains unchanged is green. If this green \(b\)-cell belongs to
  \(S_2\), then all strip endings in \(\pi_\circ\) are changed and the
  last letter of \(\type(\lambda^{(2)})\) is \(0\). These two cases
  are illustrated by the last two diagrams in
  Figure~\ref{fig:delta_iso_2abcd}~(2a). Therefore,
  \[
    \type(\lambda_\circ^{(1)},\lambda_\circ^{(2)})=r\,0\,1010\cdots10\sim r0, \qand \type(\lambda^{(1)},\lambda^{(2)})=r\,0\,0101\cdots01\,00\sim r000,
  \]
  for some binary word \(r\). 
  
\item Assume that \( T_i \), \( T_{i+1} \), and \( B_{i+1} \) are a blue \( b \)-cell, a \( c \)-cell, and a \( c \)-cell, respectively,
  illustrated in Figure~\ref{fig:delta_iso_2abcd}~(2b). Then the cell to the
  left of \( B_{i+1} \) is a blue \(b\)-cell. We have
  \[
    \type(\pi_\circ)=s0, \qand \type(\pi)=s1,
  \]
  for some binary word \(s\).

  An even number of strips descend one row lower from the \(i\)th
  column to the \((i+1)\)th column; and a new strip starts in the top
  cell of the \((i+1)\)th column, which is a blue cell. All other
  strip endings in \(\pi_\circ\) are unchanged. Note that the topmost
  strip in \(\pi\) starts with a \(c\)-cell, so the first row of
  \(\lambda^{(1)}\) ends in blue; the topmost strip in \(\pi_\circ\)
  starts with an \(a\)-cell, so the first row of \(\lambda_\circ^{(1)}\)
  ends in green. Therefore,
  \[
    \type(\lambda_\circ^{(1)},\lambda_\circ^{(2)})=r\,0101\cdots01\,0\sim r0, \qand \type(\lambda^{(1)},\lambda^{(2)})=r\,1010\cdots10\,1\sim r1,
  \]
  for some binary word \(r\).

\item Assume that \( T_i \), \( T_{i+1} \), and \( B_{i+1} \) are a green \( b \)-cell, an \( a \)-cell, and a \( c \)-cell, respectively,
  illustrated in Figure~\ref{fig:delta_iso_2abcd}~(2c). Then the cell to the
  left of \( B_{i+1} \) is a blue \(b\)-cell. We have
  \[
    \type(\pi_\circ)=s0, \qand \type(\pi)=s100,
  \]
  for some binary word \(s\).

  An odd number of strips descend one row lower from the \(i\)th
  column to the \((i+1)\)th column; and a new strip starts in the top
  cell of the \((i+1)\)th column, which is a green cell. All other
  strip endings in \(\pi_\circ\) are unchanged. Note that the topmost
  strip in \(\pi\) starts with an \(a\)-cell, so the first row of
  \(\lambda^{(1)}\) ends in green; the topmost strip in \(\pi_\circ\)
  starts with a \(c\)-cell, so the first row of \(\lambda_\circ^{(1)}\)
  ends in blue. Therefore,
  \[
    \type(\lambda_\circ^{(1)},\lambda_\circ^{(2)})=r\,0101\cdots010\sim r0, \qand \type(\lambda^{(1)},\lambda^{(2)})=r\,1010\cdots101\,00\sim r100,
  \]
  for some binary word \(r\). 

\item Assume that \( T_i \), \( T_{i+1} \), and \( B_{i+1} \) are a blue \( b \)-cell, a \( c \)-cell, and an \( a \)-cell, respectively,
  illustrated in Figure~\ref{fig:delta_iso_2abcd}~(2d). Then the cell to the
  left of \( B_{i+1} \) is a green \(b\)-cell. We have
  \[
    \type(\pi_\circ)=s0, \qand \type(\pi)=s0,
  \]
  for some binary word \(s\).

  An odd number of strips descend one row lower from the \(i\)th
  column to the \((i+1)\)th column; and a new strip starts in the top
  cell of the \((i+1)\)th column, which is a blue cell. All other
  strip endings in \(\pi_\circ\) are unchanged. Note that the topmost
  strip in \(\pi\) starts with a \(c\)-cell, so the first row of
  \(\lambda^{(1)}\) ends in blue; the topmost strip in \(\pi_\circ\)
  starts with an \(a\)-cell, so the first row of \(\lambda_\circ^{(1)}\)
  ends in green. Moreover, if the green \(b\)-cell to the left of
  \( B_{i+1} \) belongs to \(S_1\), then the topmost strip ending that
  remains unchanged is green. If this green \(b\)-cell
  belongs to \(S_2\), then all strip endings in \(\pi_\circ\) are
  changed and the last letter of \(\type(\lambda^{(2)})\) is \(0\). These
  two cases are illustrated by the last two diagrams in
  Figure~\ref{fig:delta_iso_2abcd}~(2d). Therefore,
  \[
    \type(\lambda_\circ^{(1)},\lambda_\circ^{(2)})=r\,0\,1010\cdots101\,0\sim
    r0, \qand \type(\lambda^{(1)},\lambda^{(2)})=r\,0\,0101\cdots010\,1\sim r0,
  \]
 for some binary word \(r\).

\begin{figure}[h]
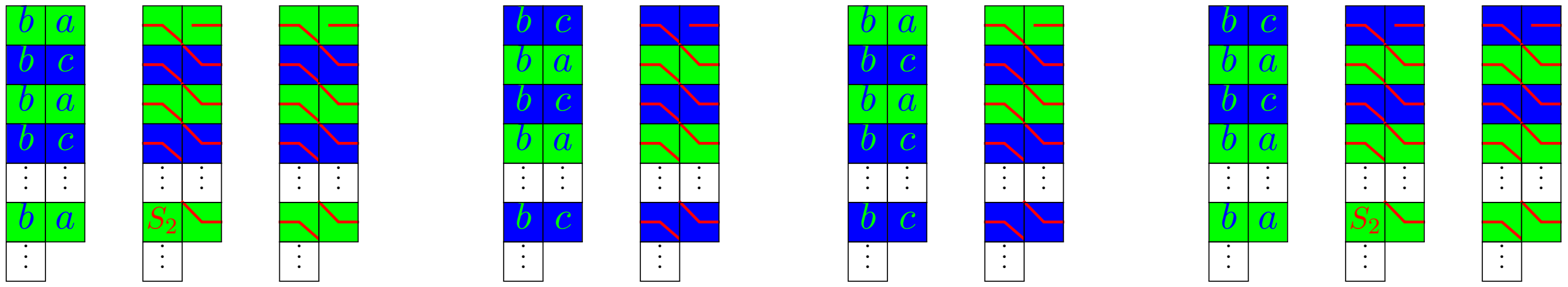

    \centering

    %%%%%%%%%%%%%%%%%%%%%%%%%%%%%%%
    % case (2a)
    %%%%%%%%%%%%%%%%%%%%%%%%%%%%%%%
    \begin{minipage}[t]{0.29\textwidth}
        \centering
        \scalebox{0.7}{%
        \begin{ytableau}
            *(green) \scalebox{1.5}{\color{blue}$b$} & *(green) \scalebox{1.5}{\color{blue}$a$} \\
            *(blue) \scalebox{1.5}{\color{green}$b$} & *(blue) \scalebox{1.5}{\color{green}$c$} \\
            *(green) \scalebox{1.5}{\color{blue}$b$} & *(green) \scalebox{1.5}{\color{blue}$a$}\\
            *(blue) \scalebox{1.5}{\color{green}$b$} & *(blue) \scalebox{1.5}{\color{green}$c$} \\
            \vdots & \vdots \\
            *(green) \scalebox{1.5}{\color{blue}$b$} & *(green) \scalebox{1.5}{\color{blue}$a$}\\
            \vdots
        \end{ytableau}
        }
        \quad
        \scalebox{0.7}{%
        \begin{ytableau}
            *(green) \tikz[baseline=(c.center)]{
                \node[minimum size=2em, inner sep=0pt] (c) {};
                \draw[red, very thick] ([xshift=-0.8em,yshift=0.6ex]c.center) -- ([xshift=-0em,yshift=0.6ex]c.center) -- ([xshift=0.8em,yshift=-1ex]c.center);
            } & *(green) \tikz[baseline=(c.center)]{
                \node[minimum size=2em, inner sep=0pt] (c) {};
                \draw[red, very thick] ([xshift=-0.4em,yshift=0.6ex]c.center) -- ([xshift=0.8em,yshift=0.6ex]c.center);
            } \\
            *(blue) \tikz[baseline=(c.center)]{
                \node[minimum size=2em, inner sep=0pt] (c) {};
                \draw[red, very thick] ([xshift=-0.8em,yshift=0.6ex]c.center) -- ([xshift=-0em,yshift=0.6ex]c.center) -- ([xshift=0.8em,yshift=-1ex]c.center);
            } & *(blue) \tikz[baseline=(c.center)]{
                \node[minimum size=2em, inner sep=0pt] (c) {};
                \draw[red, very thick] ([xshift=-0.8em,yshift=2.5ex]c.center) -- ([xshift=-0em,yshift=0.6ex]c.center) -- ([xshift=0.8em,yshift=0.6ex]c.center);
            } \\
            *(green) \tikz[baseline=(c.center)]{
                \node[minimum size=2em, inner sep=0pt] (c) {};
                \draw[red, very thick] ([xshift=-0.8em,yshift=0.6ex]c.center) -- ([xshift=-0em,yshift=0.6ex]c.center) -- ([xshift=0.8em,yshift=-1ex]c.center);
            } & *(green) \tikz[baseline=(c.center)]{
                \node[minimum size=2em, inner sep=0pt] (c) {};
                \draw[red, very thick] ([xshift=-0.8em,yshift=2.5ex]c.center) -- ([xshift=-0em,yshift=0.6ex]c.center) -- ([xshift=0.8em,yshift=0.6ex]c.center);
            }\\
            *(blue) \tikz[baseline=(c.center)]{
                \node[minimum size=2em, inner sep=0pt] (c) {};
                \draw[red, very thick] ([xshift=-0.8em,yshift=0.6ex]c.center) -- ([xshift=-0em,yshift=0.6ex]c.center) -- ([xshift=0.8em,yshift=-1ex]c.center);
            } & *(blue) \tikz[baseline=(c.center)]{
                \node[minimum size=2em, inner sep=0pt] (c) {};
                \draw[red, very thick] ([xshift=-0.8em,yshift=2.5ex]c.center) -- ([xshift=-0em,yshift=0.6ex]c.center) -- ([xshift=0.8em,yshift=0.6ex]c.center);
            } \\
            \vdots & \vdots \\
            *(green) \scalebox{1.2}{\color{red}$S_2$} & *(green) \tikz[baseline=(c.center)]{
                \node[minimum size=2em, inner sep=0pt] (c) {};
                \draw[red, very thick] ([xshift=-0.8em,yshift=2.5ex]c.center) -- ([xshift=-0em,yshift=0.6ex]c.center) -- ([xshift=0.8em,yshift=0.6ex]c.center);
            } \\
            \vdots
        \end{ytableau}
        }
        \quad
        \scalebox{0.7}{%
        \begin{ytableau}
            *(green) \tikz[baseline=(c.center)]{
                \node[minimum size=2em, inner sep=0pt] (c) {};
                \draw[red, very thick] ([xshift=-0.8em,yshift=0.6ex]c.center) -- ([xshift=-0em,yshift=0.6ex]c.center) -- ([xshift=0.8em,yshift=-1ex]c.center);
            } & *(green) \tikz[baseline=(c.center)]{
                \node[minimum size=2em, inner sep=0pt] (c) {};
                \draw[red, very thick] ([xshift=-0.4em,yshift=0.6ex]c.center) -- ([xshift=0.8em,yshift=0.6ex]c.center);
            } \\
            *(blue) \tikz[baseline=(c.center)]{
                \node[minimum size=2em, inner sep=0pt] (c) {};
                \draw[red, very thick] ([xshift=-0.8em,yshift=0.6ex]c.center) -- ([xshift=-0em,yshift=0.6ex]c.center) -- ([xshift=0.8em,yshift=-1ex]c.center);
            } & *(blue) \tikz[baseline=(c.center)]{
                \node[minimum size=2em, inner sep=0pt] (c) {};
                \draw[red, very thick] ([xshift=-0.8em,yshift=2.5ex]c.center) -- ([xshift=-0em,yshift=0.6ex]c.center) -- ([xshift=0.8em,yshift=0.6ex]c.center);
            } \\
            *(green) \tikz[baseline=(c.center)]{
                \node[minimum size=2em, inner sep=0pt] (c) {};
                \draw[red, very thick] ([xshift=-0.8em,yshift=0.6ex]c.center) -- ([xshift=-0em,yshift=0.6ex]c.center) -- ([xshift=0.8em,yshift=-1ex]c.center);
            } & *(green) \tikz[baseline=(c.center)]{
                \node[minimum size=2em, inner sep=0pt] (c) {};
                \draw[red, very thick] ([xshift=-0.8em,yshift=2.5ex]c.center) -- ([xshift=-0em,yshift=0.6ex]c.center) -- ([xshift=0.8em,yshift=0.6ex]c.center);
            }\\
            *(blue) \tikz[baseline=(c.center)]{
                \node[minimum size=2em, inner sep=0pt] (c) {};
                \draw[red, very thick] ([xshift=-0.8em,yshift=0.6ex]c.center) -- ([xshift=-0em,yshift=0.6ex]c.center) -- ([xshift=0.8em,yshift=-1ex]c.center);
            } & *(blue) \tikz[baseline=(c.center)]{
                \node[minimum size=2em, inner sep=0pt] (c) {};
                \draw[red, very thick] ([xshift=-0.8em,yshift=2.5ex]c.center) -- ([xshift=-0em,yshift=0.6ex]c.center) -- ([xshift=0.8em,yshift=0.6ex]c.center);
            } \\
            \vdots & \vdots \\
            *(green) \tikz[baseline=(c.center)]{
                \node[minimum size=2em, inner sep=0pt] (c) {};
                \draw[red, very thick] ([xshift=-0.8em,yshift=0.6ex]c.center) -- ([xshift=-0em,yshift=0.6ex]c.center) -- ([xshift=0.8em,yshift=-1ex]c.center);
            } & *(green) \tikz[baseline=(c.center)]{
                \node[minimum size=2em, inner sep=0pt] (c) {};
                \draw[red, very thick] ([xshift=-0.8em,yshift=2.5ex]c.center) -- ([xshift=-0em,yshift=0.6ex]c.center) -- ([xshift=0.8em,yshift=0.6ex]c.center);
            } \\
            \vdots
        \end{ytableau}
        }

        \medskip
        {\small (2a)}
    \end{minipage}
    \hfill
    %%%%%%%%%%%%%%%%%%%%%%%%%%%%%%%
    % case (2b)
    %%%%%%%%%%%%%%%%%%%%%%%%%%%%%%%
    \begin{minipage}[t]{0.19\textwidth}
        \centering
        \scalebox{0.7}{%
        \begin{ytableau}
            *(blue) \scalebox{1.5}{\color{green}$b$} & *(blue) \scalebox{1.5}{\color{green}$c$} \\
            *(green) \scalebox{1.5}{\color{blue}$b$} & *(green) \scalebox{1.5}{\color{blue}$a$}\\
            *(blue) \scalebox{1.5}{\color{green}$b$} & *(blue) \scalebox{1.5}{\color{green}$c$} \\
            *(green) \scalebox{1.5}{\color{blue}$b$} & *(green) \scalebox{1.5}{\color{blue}$a$} \\
            \vdots & \vdots \\
            *(blue) \scalebox{1.5}{\color{green}$b$} & *(blue) \scalebox{1.5}{\color{green}$c$} \\
            \vdots
        \end{ytableau}
        }
        \quad
        \scalebox{0.7}{%
        \begin{ytableau}
            *(blue) \tikz[baseline=(c.center)]{
                \node[minimum size=2em, inner sep=0pt] (c) {};
                \draw[red, very thick] ([xshift=-0.8em,yshift=0.6ex]c.center) -- ([xshift=-0em,yshift=0.6ex]c.center) -- ([xshift=0.8em,yshift=-1ex]c.center);
            } & *(blue) \tikz[baseline=(c.center)]{
                \node[minimum size=2em, inner sep=0pt] (c) {};
                \draw[red, very thick] ([xshift=-0.4em,yshift=0.6ex]c.center) -- ([xshift=0.8em,yshift=0.6ex]c.center);
            } \\
            *(green) \tikz[baseline=(c.center)]{
                \node[minimum size=2em, inner sep=0pt] (c) {};
                \draw[red, very thick] ([xshift=-0.8em,yshift=0.6ex]c.center) -- ([xshift=-0em,yshift=0.6ex]c.center) -- ([xshift=0.8em,yshift=-1ex]c.center);
            } & *(green) \tikz[baseline=(c.center)]{
                \node[minimum size=2em, inner sep=0pt] (c) {};
                \draw[red, very thick] ([xshift=-0.8em,yshift=2.5ex]c.center) -- ([xshift=-0em,yshift=0.6ex]c.center) -- ([xshift=0.8em,yshift=0.6ex]c.center);
            }\\
            *(blue) \tikz[baseline=(c.center)]{
                \node[minimum size=2em, inner sep=0pt] (c) {};
                \draw[red, very thick] ([xshift=-0.8em,yshift=0.6ex]c.center) -- ([xshift=-0em,yshift=0.6ex]c.center) -- ([xshift=0.8em,yshift=-1ex]c.center);
            } & *(blue) \tikz[baseline=(c.center)]{
                \node[minimum size=2em, inner sep=0pt] (c) {};
                \draw[red, very thick] ([xshift=-0.8em,yshift=2.5ex]c.center) -- ([xshift=-0em,yshift=0.6ex]c.center) -- ([xshift=0.8em,yshift=0.6ex]c.center);
            } \\
            *(green) \tikz[baseline=(c.center)]{
                \node[minimum size=2em, inner sep=0pt] (c) {};
                \draw[red, very thick] ([xshift=-0.8em,yshift=0.6ex]c.center) -- ([xshift=-0em,yshift=0.6ex]c.center) -- ([xshift=0.8em,yshift=-1ex]c.center);
            } & *(green) \tikz[baseline=(c.center)]{
                \node[minimum size=2em, inner sep=0pt] (c) {};
                \draw[red, very thick] ([xshift=-0.8em,yshift=2.5ex]c.center) -- ([xshift=-0em,yshift=0.6ex]c.center) -- ([xshift=0.8em,yshift=0.6ex]c.center);
            } \\
            \vdots & \vdots \\
            *(blue) \tikz[baseline=(c.center)]{
                \node[minimum size=2em, inner sep=0pt] (c) {};
                \draw[red, very thick] ([xshift=-0.8em,yshift=0.6ex]c.center) -- ([xshift=-0em,yshift=0.6ex]c.center) -- ([xshift=0.8em,yshift=-1ex]c.center);
            } & *(blue) \tikz[baseline=(c.center)]{
                \node[minimum size=2em, inner sep=0pt] (c) {};
                \draw[red, very thick] ([xshift=-0.8em,yshift=2.5ex]c.center) -- ([xshift=-0em,yshift=0.6ex]c.center) -- ([xshift=0.8em,yshift=0.6ex]c.center);
            } \\
            \vdots
        \end{ytableau}
        }

        \medskip
        {\small (2b)}
    \end{minipage}
    \hfill
    %%%%%%%%%%%%%%%%%%%%%%%%%%%%%%%
    % case (2c)
    %%%%%%%%%%%%%%%%%%%%%%%%%%%%%%%
    \begin{minipage}[t]{0.19\textwidth}
        \centering
        \scalebox{0.7}{%
        \begin{ytableau}
            *(green) \scalebox{1.5}{\color{blue}$b$} & *(green) \scalebox{1.5}{\color{blue}$a$} \\
            *(blue) \scalebox{1.5}{\color{green}$b$} & *(blue) \scalebox{1.5}{\color{green}$c$} \\
            *(green) \scalebox{1.5}{\color{blue}$b$} & *(green) \scalebox{1.5}{\color{blue}$a$}\\
            *(blue) \scalebox{1.5}{\color{green}$b$} & *(blue) \scalebox{1.5}{\color{green}$c$} \\
            \vdots & \vdots \\
            *(blue) \scalebox{1.5}{\color{green}$b$} & *(blue) \scalebox{1.5}{\color{green}$c$}\\
            \vdots
        \end{ytableau}
        }
        \quad
        \scalebox{0.7}{%
        \begin{ytableau}
            *(green) \tikz[baseline=(c.center)]{
                \node[minimum size=2em, inner sep=0pt] (c) {};
                \draw[red, very thick] ([xshift=-0.8em,yshift=0.6ex]c.center) -- ([xshift=-0em,yshift=0.6ex]c.center) -- ([xshift=0.8em,yshift=-1ex]c.center);
            } & *(green) \tikz[baseline=(c.center)]{
                \node[minimum size=2em, inner sep=0pt] (c) {};
                \draw[red, very thick] ([xshift=-0.4em,yshift=0.6ex]c.center) -- ([xshift=0.8em,yshift=0.6ex]c.center);
            } \\
            *(blue) \tikz[baseline=(c.center)]{
                \node[minimum size=2em, inner sep=0pt] (c) {};
                \draw[red, very thick] ([xshift=-0.8em,yshift=0.6ex]c.center) -- ([xshift=-0em,yshift=0.6ex]c.center) -- ([xshift=0.8em,yshift=-1ex]c.center);
            } & *(blue) \tikz[baseline=(c.center)]{
                \node[minimum size=2em, inner sep=0pt] (c) {};
                \draw[red, very thick] ([xshift=-0.8em,yshift=2.5ex]c.center) -- ([xshift=-0em,yshift=0.6ex]c.center) -- ([xshift=0.8em,yshift=0.6ex]c.center);
            } \\
            *(green) \tikz[baseline=(c.center)]{
                \node[minimum size=2em, inner sep=0pt] (c) {};
                \draw[red, very thick] ([xshift=-0.8em,yshift=0.6ex]c.center) -- ([xshift=-0em,yshift=0.6ex]c.center) -- ([xshift=0.8em,yshift=-1ex]c.center);
            } & *(green) \tikz[baseline=(c.center)]{
                \node[minimum size=2em, inner sep=0pt] (c) {};
                \draw[red, very thick] ([xshift=-0.8em,yshift=2.5ex]c.center) -- ([xshift=-0em,yshift=0.6ex]c.center) -- ([xshift=0.8em,yshift=0.6ex]c.center);
            }\\
            *(blue) \tikz[baseline=(c.center)]{
                \node[minimum size=2em, inner sep=0pt] (c) {};
                \draw[red, very thick] ([xshift=-0.8em,yshift=0.6ex]c.center) -- ([xshift=-0em,yshift=0.6ex]c.center) -- ([xshift=0.8em,yshift=-1ex]c.center);
            } & *(blue) \tikz[baseline=(c.center)]{
                \node[minimum size=2em, inner sep=0pt] (c) {};
                \draw[red, very thick] ([xshift=-0.8em,yshift=2.5ex]c.center) -- ([xshift=-0em,yshift=0.6ex]c.center) -- ([xshift=0.8em,yshift=0.6ex]c.center);
            } \\
            \vdots & \vdots \\
            *(blue) \tikz[baseline=(c.center)]{
                \node[minimum size=2em, inner sep=0pt] (c) {};
                \draw[red, very thick] ([xshift=-0.8em,yshift=0.6ex]c.center) -- ([xshift=-0em,yshift=0.6ex]c.center) -- ([xshift=0.8em,yshift=-1ex]c.center);
            } & *(blue) \tikz[baseline=(c.center)]{
                \node[minimum size=2em, inner sep=0pt] (c) {};
                \draw[red, very thick] ([xshift=-0.8em,yshift=2.5ex]c.center) -- ([xshift=-0em,yshift=0.6ex]c.center) -- ([xshift=0.8em,yshift=0.6ex]c.center);
            } \\
            \vdots
        \end{ytableau}
        }

        \medskip
        {\small (2c)}
    \end{minipage}
    \hfill
    %%%%%%%%%%%%%%%%%%%%%%%%%%%%%%%
    % case (2d)
    %%%%%%%%%%%%%%%%%%%%%%%%%%%%%%%
    \begin{minipage}[t]{0.29\textwidth}
        \centering
        \scalebox{0.7}{%
        \begin{ytableau}
            *(blue) \scalebox{1.5}{\color{green}$b$} & *(blue) \scalebox{1.5}{\color{green}$c$} \\
            *(green) \scalebox{1.5}{\color{blue}$b$} & *(green) \scalebox{1.5}{\color{blue}$a$}\\
            *(blue) \scalebox{1.5}{\color{green}$b$} & *(blue) \scalebox{1.5}{\color{green}$c$} \\
            *(green) \scalebox{1.5}{\color{blue}$b$} & *(green) \scalebox{1.5}{\color{blue}$a$} \\
            \vdots & \vdots \\
            *(green) \scalebox{1.5}{\color{blue}$b$} & *(green) \scalebox{1.5}{\color{blue}$a$} \\
            \vdots
        \end{ytableau}
        }
        \quad
        \scalebox{0.7}{%
        \begin{ytableau}
            *(blue) \tikz[baseline=(c.center)]{
                \node[minimum size=2em, inner sep=0pt] (c) {};
                \draw[red, very thick] ([xshift=-0.8em,yshift=0.6ex]c.center) -- ([xshift=-0em,yshift=0.6ex]c.center) -- ([xshift=0.8em,yshift=-1ex]c.center);
            } & *(blue) \tikz[baseline=(c.center)]{
                \node[minimum size=2em, inner sep=0pt] (c) {};
                \draw[red, very thick] ([xshift=-0.4em,yshift=0.6ex]c.center) -- ([xshift=0.8em,yshift=0.6ex]c.center);
            } \\
            *(green) \tikz[baseline=(c.center)]{
                \node[minimum size=2em, inner sep=0pt] (c) {};
                \draw[red, very thick] ([xshift=-0.8em,yshift=0.6ex]c.center) -- ([xshift=-0em,yshift=0.6ex]c.center) -- ([xshift=0.8em,yshift=-1ex]c.center);
            } & *(green) \tikz[baseline=(c.center)]{
                \node[minimum size=2em, inner sep=0pt] (c) {};
                \draw[red, very thick] ([xshift=-0.8em,yshift=2.5ex]c.center) -- ([xshift=-0em,yshift=0.6ex]c.center) -- ([xshift=0.8em,yshift=0.6ex]c.center);
            } \\
            *(blue) \tikz[baseline=(c.center)]{
                \node[minimum size=2em, inner sep=0pt] (c) {};
                \draw[red, very thick] ([xshift=-0.8em,yshift=0.6ex]c.center) -- ([xshift=-0em,yshift=0.6ex]c.center) -- ([xshift=0.8em,yshift=-1ex]c.center);
            } & *(blue) \tikz[baseline=(c.center)]{
                \node[minimum size=2em, inner sep=0pt] (c) {};
                \draw[red, very thick] ([xshift=-0.8em,yshift=2.5ex]c.center) -- ([xshift=-0em,yshift=0.6ex]c.center) -- ([xshift=0.8em,yshift=0.6ex]c.center);
            }\\
            *(green) \tikz[baseline=(c.center)]{
                \node[minimum size=2em, inner sep=0pt] (c) {};
                \draw[red, very thick] ([xshift=-0.8em,yshift=0.6ex]c.center) -- ([xshift=-0em,yshift=0.6ex]c.center) -- ([xshift=0.8em,yshift=-1ex]c.center);
            } & *(green) \tikz[baseline=(c.center)]{
                \node[minimum size=2em, inner sep=0pt] (c) {};
                \draw[red, very thick] ([xshift=-0.8em,yshift=2.5ex]c.center) -- ([xshift=-0em,yshift=0.6ex]c.center) -- ([xshift=0.8em,yshift=0.6ex]c.center);
            } \\
            \vdots & \vdots \\
            *(green) \scalebox{1.2}{\color{red}$S_2$} & *(green) \tikz[baseline=(c.center)]{
                \node[minimum size=2em, inner sep=0pt] (c) {};
                \draw[red, very thick] ([xshift=-0.8em,yshift=2.5ex]c.center) -- ([xshift=-0em,yshift=0.6ex]c.center) -- ([xshift=0.8em,yshift=0.6ex]c.center);
            } \\
            \vdots
        \end{ytableau}
        }
        \quad
        \scalebox{0.7}{%
        \begin{ytableau}
            *(blue) \tikz[baseline=(c.center)]{
                \node[minimum size=2em, inner sep=0pt] (c) {};
                \draw[red, very thick] ([xshift=-0.8em,yshift=0.6ex]c.center) -- ([xshift=-0em,yshift=0.6ex]c.center) -- ([xshift=0.8em,yshift=-1ex]c.center);
            } & *(blue) \tikz[baseline=(c.center)]{
                \node[minimum size=2em, inner sep=0pt] (c) {};
                \draw[red, very thick] ([xshift=-0.4em,yshift=0.6ex]c.center) -- ([xshift=0.8em,yshift=0.6ex]c.center);
            } \\
            *(green) \tikz[baseline=(c.center)]{
                \node[minimum size=2em, inner sep=0pt] (c) {};
                \draw[red, very thick] ([xshift=-0.8em,yshift=0.6ex]c.center) -- ([xshift=-0em,yshift=0.6ex]c.center) -- ([xshift=0.8em,yshift=-1ex]c.center);
            } & *(green) \tikz[baseline=(c.center)]{
                \node[minimum size=2em, inner sep=0pt] (c) {};
                \draw[red, very thick] ([xshift=-0.8em,yshift=2.5ex]c.center) -- ([xshift=-0em,yshift=0.6ex]c.center) -- ([xshift=0.8em,yshift=0.6ex]c.center);
            } \\
            *(blue) \tikz[baseline=(c.center)]{
                \node[minimum size=2em, inner sep=0pt] (c) {};
                \draw[red, very thick] ([xshift=-0.8em,yshift=0.6ex]c.center) -- ([xshift=-0em,yshift=0.6ex]c.center) -- ([xshift=0.8em,yshift=-1ex]c.center);
            } & *(blue) \tikz[baseline=(c.center)]{
                \node[minimum size=2em, inner sep=0pt] (c) {};
                \draw[red, very thick] ([xshift=-0.8em,yshift=2.5ex]c.center) -- ([xshift=-0em,yshift=0.6ex]c.center) -- ([xshift=0.8em,yshift=0.6ex]c.center);
            }\\
            *(green) \tikz[baseline=(c.center)]{
                \node[minimum size=2em, inner sep=0pt] (c) {};
                \draw[red, very thick] ([xshift=-0.8em,yshift=0.6ex]c.center) -- ([xshift=-0em,yshift=0.6ex]c.center) -- ([xshift=0.8em,yshift=-1ex]c.center);
            } & *(green) \tikz[baseline=(c.center)]{
                \node[minimum size=2em, inner sep=0pt] (c) {};
                \draw[red, very thick] ([xshift=-0.8em,yshift=2.5ex]c.center) -- ([xshift=-0em,yshift=0.6ex]c.center) -- ([xshift=0.8em,yshift=0.6ex]c.center);
            } \\
            \vdots & \vdots \\
            *(green) \tikz[baseline=(c.center)]{
                \node[minimum size=2em, inner sep=0pt] (c) {};
                \draw[red, very thick] ([xshift=-0.8em,yshift=0.6ex]c.center) -- ([xshift=-0em,yshift=0.6ex]c.center) -- ([xshift=0.8em,yshift=-1ex]c.center);
            } & *(green) \tikz[baseline=(c.center)]{
                \node[minimum size=2em, inner sep=0pt] (c) {};
                \draw[red, very thick] ([xshift=-0.8em,yshift=2.5ex]c.center) -- ([xshift=-0em,yshift=0.6ex]c.center) -- ([xshift=0.8em,yshift=0.6ex]c.center);
            } \\
            \vdots
        \end{ytableau}
        }

        \medskip
        {\small (2d)}
    \end{minipage}

    \caption{Inductive steps for cases (2a), (2b), (2c), and (2d).}
    \label{fig:delta_iso_2abcd}
\end{figure}

\item Assume that \( T_i \) and \( B_{i+1} \) are \( c \)-cells and a
  blue \( b \)-cell, respectively, illustrated in
  Figure~\ref{fig:delta_iso_2efg}~(2e). Then \( T_{i+1} \) is a green
  \(b\)-cell, and the cell to the left of \( B_{i+1} \) is an
  \(a\)-cell. We have
  \[
    \type(\pi_\circ)=s, \qand \type(\pi)=s10,
  \]
  for some binary word \(s\).

  An even number of strips traverse the \(i\)th and the \((i+1)\)th
  column horizontally. All other strip endings in \(\pi_\circ\) are
  unchanged. Note that the topmost strip in both \(\pi\) and
  \(\pi_\circ\) starts with a \(c\)-cell, so the first row of both
  \(\lambda^{(1)}\) and \(\lambda_\circ^{(1)}\) ends in blue. Therefore,
  \[
    \type(\lambda_\circ^{(1)},\lambda_\circ^{(2)})=r\,0101\cdots01\sim r, \qand \type(\lambda^{(1)},\lambda^{(2)})=r\,1010\cdots10\sim r10,
  \]
  for some binary word \(r\).

\item Assume that \( T_i \), \( T_{i+1} \), and \( B_{i+1} \) are a green \( b \)-cell, a \( c \)-cell, and a \( c \)-cell, respectively,
  illustrated in Figure~\ref{fig:delta_iso_2efg}~(2f). Then the cell to the
  left of \( B_{i+1} \) is a green \(b\)-cell. We have
  \[
    \type(\pi_\circ)=s0, \qand \type(\pi)=s1,
  \]
  for some binary word \(s\).

  An odd number of strips traverse the \(i\)th and the \((i+1)\)th
  column horizontally. All other strip endings in \(\pi_\circ\) are
  unchanged. Note that the topmost strip in both \(\pi\) and
  \(\pi_\circ\) starts with a \(c\)-cell, so the first row of both
  \(\lambda^{(1)}\) and \(\lambda_\circ^{(1)}\) ends in blue. Therefore,
  \[
    \type(\lambda_\circ^{(1)},\lambda_\circ^{(2)})=r\,0101\cdots010\sim r0, \qand \type(\lambda^{(1)},\lambda^{(2)})=r\,1010\cdots101\sim r1
  \]
  for some binary word \(r\).

\item Assume that \( T_i \), \( T_{i+1} \), and \( B_{i+1} \) are a green \( b \)-cell, a \( c \)-cell, and an \( a \)-cell, respectively,
  illustrated in Figure~\ref{fig:delta_iso_2efg}~(2g). Then the \(b\)-cell to
  the left of \( B_{i+1} \) is blue. We have
  \[
    \type(\pi_\circ)=s0, \qand \type(\pi)=s0,
  \]
  for some binary word \(s\).

  An even number of strips traverse the \(i\)th and the \((i+1)\)th
  column horizontally. All other strip endings in \(\pi_\circ\) are
  unchanged. Note that the topmost strip in both \(\pi\) and
  \(\pi_\circ\) starts with a \(c\)-cell, so the first row of both
  \(\lambda^{(1)}\) and \(\lambda_\circ^{(1)}\) ends in blue. Since \( T_i \)
  is a green \(b\)-cell, the cells in the \((i-1)\)th column are the
  same as those in the \((i+1)\)th column in each row. Thus, by the
  difference conditions of \(\mathcal{P}_{2,b}\), the \(i\)th and
  \((i+1)\)th columns cannot have the same length (see the first
  diagram in Figure~\ref{fig:delta_iso_2efg}~(2g)). Therefore, there is
  a green \(b\)-cell in the \(i\)th column, diagonally
  south-west of \( B_{i+1} \). If this green \(b\)-cell belongs to
  \(S_1\), then the topmost strip ending of \(\pi_\circ\) that remains
  unchanged contributes a letter \(0\). If this green \(b\)-cell
  belongs to \(S_2\), then all strip endings in \(\pi_\circ\) are
  changed and the last letter of \(\type(\lambda^{(2)})\) is \(0\). These
  two cases are illustrated by the last two diagrams in
  Figure~\ref{fig:delta_iso_2efg}~(2g). Therefore,
  \[
    \type(\lambda_\circ^{(1)},\lambda_\circ^{(2)})=r\,0\,1010\cdots10\sim r0, \qand \type(\lambda^{(1)},\lambda^{(2)})=r\,0\,0101\cdots01\sim r0,
  \]
  for some binary word \(r\).
  
\begin{figure}[h]
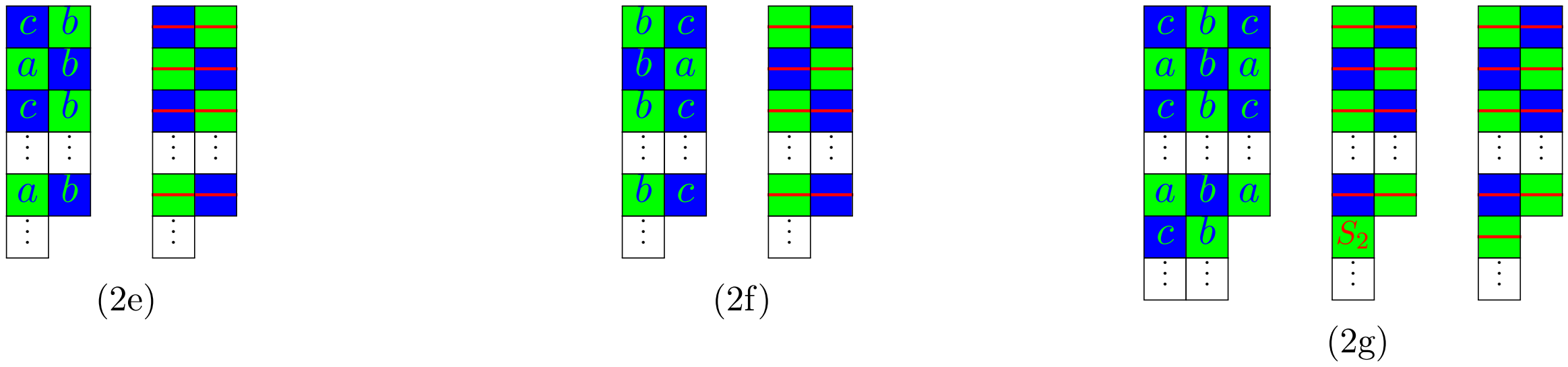

    \centering

    \begin{minipage}[t]{0.32\textwidth}
        \centering
        \scalebox{0.7}{%
        \begin{ytableau}
            *(blue) \scalebox{1.5}{\color{green}$c$} & *(green) \scalebox{1.5}{\color{blue}$b$} \\
            *(green) \scalebox{1.5}{\color{blue}$a$} & *(blue) \scalebox{1.5}{\color{green}$b$} \\
            *(blue) \scalebox{1.5}{\color{green}$c$} & *(green) \scalebox{1.5}{\color{blue}$b$} \\
            \vdots & \vdots \\
            *(green) \scalebox{1.5}{\color{blue}$a$} & *(blue) \scalebox{1.5}{\color{green}$b$} \\
            \vdots
        \end{ytableau}
        }
        \quad
        \scalebox{0.7}{%
        \begin{ytableau}
            *(blue) \tikz[baseline=(c.center)]{
                \node[minimum size=2em, inner sep=0pt] (c) {};
                \draw[red, very thick] ([xshift=-0.8em,yshift=0.6ex]c.center) -- ([xshift=0.8em,yshift=0.6ex]c.center);
            } & *(green) \tikz[baseline=(c.center)]{
                \node[minimum size=2em, inner sep=0pt] (c) {};
                \draw[red, very thick] ([xshift=-0.8em,yshift=0.6ex]c.center) -- ([xshift=0.8em,yshift=0.6ex]c.center);
            } \\
            *(green) \tikz[baseline=(c.center)]{
                \node[minimum size=2em, inner sep=0pt] (c) {};
                \draw[red, very thick] ([xshift=-0.8em,yshift=0.6ex]c.center) -- ([xshift=0.8em,yshift=0.6ex]c.center);
            } & *(blue) \tikz[baseline=(c.center)]{
                \node[minimum size=2em, inner sep=0pt] (c) {};
                \draw[red, very thick] ([xshift=-0.8em,yshift=0.6ex]c.center) -- ([xshift=0.8em,yshift=0.6ex]c.center);
            } \\
            *(blue) \tikz[baseline=(c.center)]{
                \node[minimum size=2em, inner sep=0pt] (c) {};
                \draw[red, very thick] ([xshift=-0.8em,yshift=0.6ex]c.center) -- ([xshift=0.8em,yshift=0.6ex]c.center);
            } & *(green) \tikz[baseline=(c.center)]{
                \node[minimum size=2em, inner sep=0pt] (c) {};
                \draw[red, very thick] ([xshift=-0.8em,yshift=0.6ex]c.center) -- ([xshift=0.8em,yshift=0.6ex]c.center);
            } \\
            \vdots & \vdots \\
            *(green) \tikz[baseline=(c.center)]{
                \node[minimum size=2em, inner sep=0pt] (c) {};
                \draw[red, very thick] ([xshift=-0.8em,yshift=0.6ex]c.center) -- ([xshift=0.8em,yshift=0.6ex]c.center);
            } & *(blue) \tikz[baseline=(c.center)]{
                \node[minimum size=2em, inner sep=0pt] (c) {};
                \draw[red, very thick] ([xshift=-0.8em,yshift=0.6ex]c.center) -- ([xshift=0.8em,yshift=0.6ex]c.center);
            } \\
            \vdots
        \end{ytableau}
        }

        \medskip
        {\small (2e)}
    \end{minipage}
    \hfill
    \begin{minipage}[t]{0.32\textwidth}
        \centering
        \scalebox{0.7}{%
        \begin{ytableau}
            *(green) \scalebox{1.5}{\color{blue}$b$} & *(blue) \scalebox{1.5}{\color{green}$c$} \\
            *(blue) \scalebox{1.5}{\color{green}$b$} & *(green) \scalebox{1.5}{\color{blue}$a$} \\
            *(green) \scalebox{1.5}{\color{blue}$b$} & *(blue) \scalebox{1.5}{\color{green}$c$}\\
            \vdots & \vdots \\
            *(green) \scalebox{1.5}{\color{blue}$b$} & *(blue) \scalebox{1.5}{\color{green}$c$}\\
            \vdots
        \end{ytableau}
        }
        \quad
        \scalebox{0.7}{%
        \begin{ytableau}
            *(green) \tikz[baseline=(c.center)]{
                \node[minimum size=2em, inner sep=0pt] (c) {};
                \draw[red, very thick] ([xshift=-0.8em,yshift=0.6ex]c.center) -- ([xshift=0.8em,yshift=0.6ex]c.center);
            } & *(blue) \tikz[baseline=(c.center)]{
                \node[minimum size=2em, inner sep=0pt] (c) {};
                \draw[red, very thick] ([xshift=-0.8em,yshift=0.6ex]c.center) -- ([xshift=0.8em,yshift=0.6ex]c.center);
            }\\
            *(blue) \tikz[baseline=(c.center)]{
                \node[minimum size=2em, inner sep=0pt] (c) {};
                \draw[red, very thick] ([xshift=-0.8em,yshift=0.6ex]c.center) -- ([xshift=0.8em,yshift=0.6ex]c.center);
            } & *(green) \tikz[baseline=(c.center)]{
                \node[minimum size=2em, inner sep=0pt] (c) {};
                \draw[red, very thick] ([xshift=-0.8em,yshift=0.6ex]c.center) -- ([xshift=0.8em,yshift=0.6ex]c.center);
            }\\
            *(green) \tikz[baseline=(c.center)]{
                \node[minimum size=2em, inner sep=0pt] (c) {};
                \draw[red, very thick] ([xshift=-0.8em,yshift=0.6ex]c.center) -- ([xshift=0.8em,yshift=0.6ex]c.center);
            } & *(blue) \tikz[baseline=(c.center)]{
                \node[minimum size=2em, inner sep=0pt] (c) {};
                \draw[red, very thick] ([xshift=-0.8em,yshift=0.6ex]c.center) -- ([xshift=0.8em,yshift=0.6ex]c.center);
            }\\
            \vdots & \vdots \\
            *(green) \tikz[baseline=(c.center)]{
                \node[minimum size=2em, inner sep=0pt] (c) {};
                \draw[red, very thick] ([xshift=-0.8em,yshift=0.6ex]c.center) -- ([xshift=0.8em,yshift=0.6ex]c.center);
            } & *(blue) \tikz[baseline=(c.center)]{
                \node[minimum size=2em, inner sep=0pt] (c) {};
                \draw[red, very thick] ([xshift=-0.8em,yshift=0.6ex]c.center) -- ([xshift=0.8em,yshift=0.6ex]c.center);
            } \\
            \vdots
        \end{ytableau}
        }

        \medskip
        {\small (2f)}
    \end{minipage}
    \hfill
    \begin{minipage}[t]{0.32\textwidth}
        \centering
        \scalebox{0.7}{%
        \begin{ytableau}
            *(blue) \scalebox{1.5}{\color{green}$c$} & *(green) \scalebox{1.5}{\color{blue}$b$} & *(blue) \scalebox{1.5}{\color{green}$c$} \\
            *(green) \scalebox{1.5}{\color{blue}$a$} & *(blue) \scalebox{1.5}{\color{green}$b$} & *(green) \scalebox{1.5}{\color{blue}$a$} \\
            *(blue) \scalebox{1.5}{\color{green}$c$} & *(green) \scalebox{1.5}{\color{blue}$b$} & *(blue) \scalebox{1.5}{\color{green}$c$}\\
            \vdots & \vdots & \vdots \\
            *(green) \scalebox{1.5}{\color{blue}$a$} & *(blue) \scalebox{1.5}{\color{green}$b$} & *(green) \scalebox{1.5}{\color{blue}$a$}\\
            *(blue) \scalebox{1.5}{\color{green}$c$} & *(green) \scalebox{1.5}{\color{blue}$b$} \\
            \vdots & \vdots
        \end{ytableau}
        }
        \quad
        \scalebox{0.7}{%
        \begin{ytableau}
            *(green) \tikz[baseline=(c.center)]{
                \node[minimum size=2em, inner sep=0pt] (c) {};
                \draw[red, very thick] ([xshift=-0.8em,yshift=0.6ex]c.center) -- ([xshift=0.8em,yshift=0.6ex]c.center);
            } & *(blue) \tikz[baseline=(c.center)]{
                \node[minimum size=2em, inner sep=0pt] (c) {};
                \draw[red, very thick] ([xshift=-0.8em,yshift=0.6ex]c.center) -- ([xshift=0.8em,yshift=0.6ex]c.center);
            }\\
            *(blue) \tikz[baseline=(c.center)]{
                \node[minimum size=2em, inner sep=0pt] (c) {};
                \draw[red, very thick] ([xshift=-0.8em,yshift=0.6ex]c.center) -- ([xshift=0.8em,yshift=0.6ex]c.center);
            } & *(green) \tikz[baseline=(c.center)]{
                \node[minimum size=2em, inner sep=0pt] (c) {};
                \draw[red, very thick] ([xshift=-0.8em,yshift=0.6ex]c.center) -- ([xshift=0.8em,yshift=0.6ex]c.center);
            }\\
            *(green) \tikz[baseline=(c.center)]{
                \node[minimum size=2em, inner sep=0pt] (c) {};
                \draw[red, very thick] ([xshift=-0.8em,yshift=0.6ex]c.center) -- ([xshift=0.8em,yshift=0.6ex]c.center);
            } & *(blue) \tikz[baseline=(c.center)]{
                \node[minimum size=2em, inner sep=0pt] (c) {};
                \draw[red, very thick] ([xshift=-0.8em,yshift=0.6ex]c.center) -- ([xshift=0.8em,yshift=0.6ex]c.center);
            }\\
            \vdots & \vdots \\
            *(blue) \tikz[baseline=(c.center)]{
                \node[minimum size=2em, inner sep=0pt] (c) {};
                \draw[red, very thick] ([xshift=-0.8em,yshift=0.6ex]c.center) -- ([xshift=0.8em,yshift=0.6ex]c.center);
            } & *(green) \tikz[baseline=(c.center)]{
                \node[minimum size=2em, inner sep=0pt] (c) {};
                \draw[red, very thick] ([xshift=-0.8em,yshift=0.6ex]c.center) -- ([xshift=0.8em,yshift=0.6ex]c.center);
            } \\
            *(green) \scalebox{1.2}{\color{red}$S_2$} \\
            \vdots
        \end{ytableau}
        }
        \quad
        \scalebox{0.7}{%
        \begin{ytableau}
            *(green) \tikz[baseline=(c.center)]{
                \node[minimum size=2em, inner sep=0pt] (c) {};
                \draw[red, very thick] ([xshift=-0.8em,yshift=0.6ex]c.center) -- ([xshift=0.8em,yshift=0.6ex]c.center);
            } & *(blue) \tikz[baseline=(c.center)]{
                \node[minimum size=2em, inner sep=0pt] (c) {};
                \draw[red, very thick] ([xshift=-0.8em,yshift=0.6ex]c.center) -- ([xshift=0.8em,yshift=0.6ex]c.center);
            }\\
            *(blue) \tikz[baseline=(c.center)]{
                \node[minimum size=2em, inner sep=0pt] (c) {};
                \draw[red, very thick] ([xshift=-0.8em,yshift=0.6ex]c.center) -- ([xshift=0.8em,yshift=0.6ex]c.center);
            } & *(green) \tikz[baseline=(c.center)]{
                \node[minimum size=2em, inner sep=0pt] (c) {};
                \draw[red, very thick] ([xshift=-0.8em,yshift=0.6ex]c.center) -- ([xshift=0.8em,yshift=0.6ex]c.center);
            }\\
            *(green) \tikz[baseline=(c.center)]{
                \node[minimum size=2em, inner sep=0pt] (c) {};
                \draw[red, very thick] ([xshift=-0.8em,yshift=0.6ex]c.center) -- ([xshift=0.8em,yshift=0.6ex]c.center);
            } & *(blue) \tikz[baseline=(c.center)]{
                \node[minimum size=2em, inner sep=0pt] (c) {};
                \draw[red, very thick] ([xshift=-0.8em,yshift=0.6ex]c.center) -- ([xshift=0.8em,yshift=0.6ex]c.center);
            }\\
            \vdots & \vdots \\
            *(blue) \tikz[baseline=(c.center)]{
                \node[minimum size=2em, inner sep=0pt] (c) {};
                \draw[red, very thick] ([xshift=-0.8em,yshift=0.6ex]c.center) -- ([xshift=0.8em,yshift=0.6ex]c.center);
            } & *(green) \tikz[baseline=(c.center)]{
                \node[minimum size=2em, inner sep=0pt] (c) {};
                \draw[red, very thick] ([xshift=-0.8em,yshift=0.6ex]c.center) -- ([xshift=0.8em,yshift=0.6ex]c.center);
            } \\
            *(green) \tikz[baseline=(c.center)]{
                \node[minimum size=2em, inner sep=0pt] (c) {};
                \draw[red, very thick] ([xshift=-0.8em,yshift=0.6ex]c.center) -- ([xshift=0.8em,yshift=0.6ex]c.center);
            } \\
            \vdots
        \end{ytableau}
        }

        \medskip
        {\small (2g)}
    \end{minipage}

    \caption{Inductive steps for cases (2e), (2f), and (2g).}
    \label{fig:delta_iso_2efg}
\end{figure}

\item Assume that \( T_i \) and \( B_{i+1} \) are an \( a \)-cell and
  a green \( b \)-cell, respectively, illustrated in
  Figure~\ref{fig:delta_iso_2hij}~(2h). Then \( T_{i+1} \) is a blue
  \(b\)-cell, and the cell to the left of \( B_{i+1} \) is a
  \(c\)-cell. We have
  \[
    \type(\pi_\circ)=s00, \qand \type(\pi)=s00,
  \]
  for some binary word \(s\).

  An even number of strips traverse the \(i\)th and the \((i+1)\)th
  column horizontally. All other strip endings in \(\pi_\circ\) are
  unchanged. Note that the topmost strip in both \(\pi\) and
  \(\pi_\circ\) starts with an \(a\)-cell, so the first row of both
  \(\lambda^{(1)}\) and \(\lambda_\circ^{(1)}\) ends in green. By the
  difference conditions in \(\mathcal{P}_{2,b}\), there must be an
  \(a\)-cell diagonally south-west of the \(b\)-cell \( B_{i+1} \).
  Moreover, by definition of the map \(\delta\), this cell must belong
  to \(S_1\). Thus, the topmost strip ending of \(\pi_\circ\) that
  remains unchanged contributes a letter \(0\). Therefore,
  \[
    \type(\lambda_\circ^{(1)},\lambda_\circ^{(2)})=r\,0\,1010\cdots10\,0\sim
    r00, \qand \type(\lambda^{(1)},\lambda^{(2)})=r\,0\,0101\cdots01\,0\sim r00,
  \]
 for some binary word \(r\). 

\item Assume that \( T_i \), \( T_{i+1} \), and \( B_{i+1} \) are a
  blue \( b \)-cell, an \( a \)-cell, and a \( c \)-cell,
  respectively, illustrated in
  Figure~\ref{fig:delta_iso_2hij}~(2i). Then the cell to the left
  of \( B_{i+1} \) is a green \(b\)-cell. We have
  \[
    \type(\pi_\circ)=s0, \qand \type(\pi)=s100,
  \]
  for some binary word \(s\).

  An even number of strips traverse the \(i\)th and the \((i+1)\)th
  column horizontally. All other strip endings in \(\pi_\circ\) are
  unchanged. Note that the topmost strip in both \(\pi\) and
  \(\pi_\circ\) starts with an \(a\)-cell, so the first row of both
  \(\lambda^{(1)}\) and \(\lambda_\circ^{(1)}\) ends in green. Therefore,
  \[
    \type(\lambda_\circ^{(1)},\lambda_\circ^{(2)})=r\,0101\cdots01\,0\sim r0, \qand \type(\lambda^{(1)},\lambda^{(2)})=r\,1010\cdots10\,0\sim r100,
  \]
  for some binary word \(r\).

  \item Assume that \( T_i \), \( T_{i+1} \), and \( B_{i+1} \) are a
blue \( b \)-cell, an \( a \)-cell, and an \( a \)-cell, respectively,
  illustrated in Figure~\ref{fig:delta_iso_2hij}~(2j). Then the cell to the
  left of \( B_{i+1} \) is a blue \( b \)-cell. We have
  \[
    \type(\pi_\circ)=s0, \qand \type(\pi)=s000,
  \]
  for some binary word \(s\).

  An odd number of strips traverse the \(i\)th and the \((i+1)\)th
  column horizontally. All other strip endings in \(\pi_\circ\) are
  unchanged. Note that the topmost strip in both \(\pi\) and
  \(\pi_\circ\) starts with an \(a\)-cell, so the first row of both
  \(\lambda^{(1)}\) and \(\lambda_\circ^{(1)}\) ends in green. Since
  \( T_i \) is a blue \(b\)-cell, the cells in the \((i-1)\)th column
  are the same as those in the \((i+1)\)th column in each row. Thus,
  by the difference conditions of \(\mathcal{P}_{2,b}\), the \(i\)th
  and \((i+1)\)th columns cannot have the same length (see the first
  diagram in Figure~\ref{fig:delta_iso_2hij}~(2j)). Therefore, there is a
  green \(b\)-cell in the \(i\)th column, diagonally south-west of
  \( B_{i+1} \). If this green \(b\)-cell belongs to \(S_1\), then the
  topmost strip ending of \(\pi_\circ\) that remains unchanged
  contributes a letter \(0\). If this green \(b\)-cell belongs to
  \(S_2\), then all strip endings in \(\pi_\circ\) are changed and the
  last letter of \(\type(\lambda^{(2)})\) is \(0\). These two cases
  are illustrated by the last two diagrams in
  Figure~\ref{fig:delta_iso_2hij}~(2j). Therefore,
  \[
    \type(\lambda_\circ^{(1)},\lambda_\circ^{(2)})=r\,0\,1010\cdots101\,0\sim
    r0, \qand \type(\lambda^{(1)},\lambda^{(2)})=r\,0\,0101\cdots010\,0\sim r000,
  \]
 for some binary word \(r\).
\begin{figure}[h]
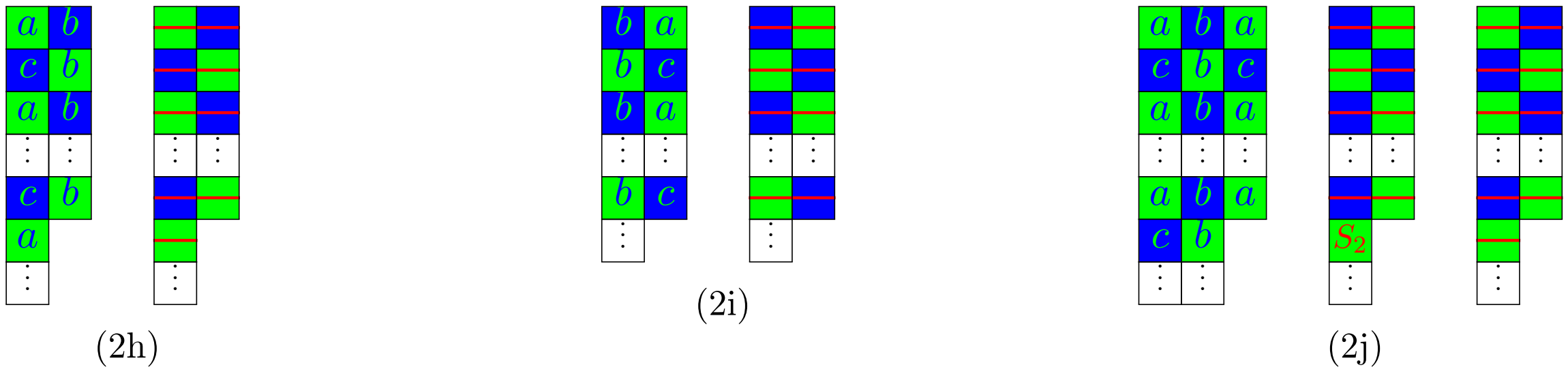

    \centering

    \begin{minipage}[t]{0.31\textwidth}
        \centering
        \scalebox{0.7}{%
        \begin{ytableau}
            *(green) \scalebox{1.5}{\color{blue}$a$}
            & *(blue) \scalebox{1.5}{\color{green}$b$} \\
            *(blue) \scalebox{1.5}{\color{green}$c$} & *(green) \scalebox{1.5}{\color{blue}$b$} \\
            *(green) \scalebox{1.5}{\color{blue}$a$} & *(blue) \scalebox{1.5}{\color{green}$b$} \\
            \vdots & \vdots \\
            *(blue) \scalebox{1.5}{\color{green}$c$} & *(green) \scalebox{1.5}{\color{blue}$b$} \\
            *(green) \scalebox{1.5}{\color{blue}$a$} \\
            \vdots
        \end{ytableau}
        }
        \quad
        \scalebox{0.7}{%
        \begin{ytableau}
            *(green) \tikz[baseline=(c.center)]{
                \node[minimum size=2em, inner sep=0pt] (c) {};
                \draw[red, very thick] ([xshift=-0.8em,yshift=0.6ex]c.center) -- ([xshift=0.8em,yshift=0.6ex]c.center);
            }
            & *(blue) \tikz[baseline=(c.center)]{
                \node[minimum size=2em, inner sep=0pt] (c) {};
                \draw[red, very thick] ([xshift=-0.8em,yshift=0.6ex]c.center) -- ([xshift=0.8em,yshift=0.6ex]c.center);
            } \\
            *(blue) \tikz[baseline=(c.center)]{
                \node[minimum size=2em, inner sep=0pt] (c) {};
                \draw[red, very thick] ([xshift=-0.8em,yshift=0.6ex]c.center) -- ([xshift=0.8em,yshift=0.6ex]c.center);
            } & *(green) \tikz[baseline=(c.center)]{
                \node[minimum size=2em, inner sep=0pt] (c) {};
                \draw[red, very thick] ([xshift=-0.8em,yshift=0.6ex]c.center) -- ([xshift=0.8em,yshift=0.6ex]c.center);
            } \\
             *(green) \tikz[baseline=(c.center)]{
                \node[minimum size=2em, inner sep=0pt] (c) {};
                \draw[red, very thick] ([xshift=-0.8em,yshift=0.6ex]c.center) -- ([xshift=0.8em,yshift=0.6ex]c.center);
            } & *(blue) \tikz[baseline=(c.center)]{
                \node[minimum size=2em, inner sep=0pt] (c) {};
                \draw[red, very thick] ([xshift=-0.8em,yshift=0.6ex]c.center) -- ([xshift=0.8em,yshift=0.6ex]c.center);
            } \\
            \vdots & \vdots \\
            *(blue) \tikz[baseline=(c.center)]{
                \node[minimum size=2em, inner sep=0pt] (c) {};
                \draw[red, very thick] ([xshift=-0.8em,yshift=0.6ex]c.center) -- ([xshift=0.8em,yshift=0.6ex]c.center);
            } & *(green) \tikz[baseline=(c.center)]{
                \node[minimum size=2em, inner sep=0pt] (c) {};
                \draw[red, very thick] ([xshift=-0.8em,yshift=0.6ex]c.center) -- ([xshift=0.8em,yshift=0.6ex]c.center);
            } \\
            *(green) \tikz[baseline=(c.center)]{
                \node[minimum size=2em, inner sep=0pt] (c) {};
                \draw[red, very thick] ([xshift=-0.8em,yshift=0.6ex]c.center) -- ([xshift=0.8em,yshift=0.6ex]c.center);
            } \\
            \vdots
        \end{ytableau}
        }

        \medskip
        {\small (2h)}
    \end{minipage}
    \hfill
    \begin{minipage}[t]{0.31\textwidth}
        \centering
        \scalebox{0.7}{%
        \begin{ytableau}
            *(blue) \scalebox{1.5}{\color{green}$b$} & *(green) \scalebox{1.5}{\color{blue}$a$} \\
            *(green) \scalebox{1.5}{\color{blue}$b$} & *(blue) \scalebox{1.5}{\color{green}$c$}\\
            *(blue) \scalebox{1.5}{\color{green}$b$} & *(green) \scalebox{1.5}{\color{blue}$a$} \\
            \vdots & \vdots \\
            *(green) \scalebox{1.5}{\color{blue}$b$} & *(blue) \scalebox{1.5}{\color{green}$c$}\\
            \vdots
        \end{ytableau}
        }
        \quad
        \scalebox{0.7}{%
        \begin{ytableau}
            *(blue) \tikz[baseline=(c.center)]{
                \node[minimum size=2em, inner sep=0pt] (c) {};
                \draw[red, very thick] ([xshift=-0.8em,yshift=0.6ex]c.center) -- ([xshift=0.8em,yshift=0.6ex]c.center);
            } & *(green) \tikz[baseline=(c.center)]{
                \node[minimum size=2em, inner sep=0pt] (c) {};
                \draw[red, very thick] ([xshift=-0.8em,yshift=0.6ex]c.center) -- ([xshift=0.8em,yshift=0.6ex]c.center);
            }\\
            *(green) \tikz[baseline=(c.center)]{
                \node[minimum size=2em, inner sep=0pt] (c) {};
                \draw[red, very thick] ([xshift=-0.8em,yshift=0.6ex]c.center) -- ([xshift=0.8em,yshift=0.6ex]c.center);
            } & *(blue) \tikz[baseline=(c.center)]{
                \node[minimum size=2em, inner sep=0pt] (c) {};
                \draw[red, very thick] ([xshift=-0.8em,yshift=0.6ex]c.center) -- ([xshift=0.8em,yshift=0.6ex]c.center);
            }\\
            *(blue) \tikz[baseline=(c.center)]{
                \node[minimum size=2em, inner sep=0pt] (c) {};
                \draw[red, very thick] ([xshift=-0.8em,yshift=0.6ex]c.center) -- ([xshift=0.8em,yshift=0.6ex]c.center);
            } & *(green) \tikz[baseline=(c.center)]{
                \node[minimum size=2em, inner sep=0pt] (c) {};
                \draw[red, very thick] ([xshift=-0.8em,yshift=0.6ex]c.center) -- ([xshift=0.8em,yshift=0.6ex]c.center);
            }\\
            \vdots & \vdots \\
            *(green) \tikz[baseline=(c.center)]{
                \node[minimum size=2em, inner sep=0pt] (c) {};
                \draw[red, very thick] ([xshift=-0.8em,yshift=0.6ex]c.center) -- ([xshift=0.8em,yshift=0.6ex]c.center);
            } & *(blue) \tikz[baseline=(c.center)]{
                \node[minimum size=2em, inner sep=0pt] (c) {};
                \draw[red, very thick] ([xshift=-0.8em,yshift=0.6ex]c.center) -- ([xshift=0.8em,yshift=0.6ex]c.center);
            } \\
            \vdots
        \end{ytableau}
        }

        \medskip
        {\small (2i)}
    \end{minipage}
    \hfill
    \begin{minipage}[t]{0.35\textwidth}
        \centering
        \scalebox{0.7}{%
        \begin{ytableau}
            *(green) \scalebox{1.5}{\color{blue}$a$} & *(blue) \scalebox{1.5}{\color{green}$b$} & *(green) \scalebox{1.5}{\color{blue}$a$} \\
            *(blue) \scalebox{1.5}{\color{green}$c$} & *(green) \scalebox{1.5}{\color{blue}$b$} & *(blue) \scalebox{1.5}{\color{green}$c$}\\
            *(green) \scalebox{1.5}{\color{blue}$a$} & *(blue) \scalebox{1.5}{\color{green}$b$} & *(green) \scalebox{1.5}{\color{blue}$a$} \\
            \vdots & \vdots & \vdots \\
            *(green) \scalebox{1.5}{\color{blue}$a$} & *(blue) \scalebox{1.5}{\color{green}$b$} & *(green) \scalebox{1.5}{\color{blue}$a$}\\
            *(blue) \scalebox{1.5}{\color{green}$c$} & *(green) \scalebox{1.5}{\color{blue}$b$} \\
            \vdots & \vdots
        \end{ytableau}
        }
        \quad
        \scalebox{0.7}{%
        \begin{ytableau}
            *(blue) \tikz[baseline=(c.center)]{
                \node[minimum size=2em, inner sep=0pt] (c) {};
                \draw[red, very thick] ([xshift=-0.8em,yshift=0.6ex]c.center) -- ([xshift=0.8em,yshift=0.6ex]c.center);
            } & *(green) \tikz[baseline=(c.center)]{
                \node[minimum size=2em, inner sep=0pt] (c) {};
                \draw[red, very thick] ([xshift=-0.8em,yshift=0.6ex]c.center) -- ([xshift=0.8em,yshift=0.6ex]c.center);
            }\\
            *(green) \tikz[baseline=(c.center)]{
                \node[minimum size=2em, inner sep=0pt] (c) {};
                \draw[red, very thick] ([xshift=-0.8em,yshift=0.6ex]c.center) -- ([xshift=0.8em,yshift=0.6ex]c.center);
            } & *(blue) \tikz[baseline=(c.center)]{
                \node[minimum size=2em, inner sep=0pt] (c) {};
                \draw[red, very thick] ([xshift=-0.8em,yshift=0.6ex]c.center) -- ([xshift=0.8em,yshift=0.6ex]c.center);
            }\\
            *(blue) \tikz[baseline=(c.center)]{
                \node[minimum size=2em, inner sep=0pt] (c) {};
                \draw[red, very thick] ([xshift=-0.8em,yshift=0.6ex]c.center) -- ([xshift=0.8em,yshift=0.6ex]c.center);
            } & *(green) \tikz[baseline=(c.center)]{
                \node[minimum size=2em, inner sep=0pt] (c) {};
                \draw[red, very thick] ([xshift=-0.8em,yshift=0.6ex]c.center) -- ([xshift=0.8em,yshift=0.6ex]c.center);
            }\\
            \vdots & \vdots \\
            *(blue) \tikz[baseline=(c.center)]{
                \node[minimum size=2em, inner sep=0pt] (c) {};
                \draw[red, very thick] ([xshift=-0.8em,yshift=0.6ex]c.center) -- ([xshift=0.8em,yshift=0.6ex]c.center);
            } & *(green) \tikz[baseline=(c.center)]{
                \node[minimum size=2em, inner sep=0pt] (c) {};
                \draw[red, very thick] ([xshift=-0.8em,yshift=0.6ex]c.center) -- ([xshift=0.8em,yshift=0.6ex]c.center);
            } \\
            *(green) \scalebox{1.2}{\color{red}$S_2$} \\
            \vdots
        \end{ytableau}
        }
        \quad
        \scalebox{0.7}{%
        \begin{ytableau}
            *(green) \tikz[baseline=(c.center)]{
                \node[minimum size=2em, inner sep=0pt] (c) {};
                \draw[red, very thick] ([xshift=-0.8em,yshift=0.6ex]c.center) -- ([xshift=0.8em,yshift=0.6ex]c.center);
            } & *(blue) \tikz[baseline=(c.center)]{
                \node[minimum size=2em, inner sep=0pt] (c) {};
                \draw[red, very thick] ([xshift=-0.8em,yshift=0.6ex]c.center) -- ([xshift=0.8em,yshift=0.6ex]c.center);
            }\\
            *(blue) \tikz[baseline=(c.center)]{
                \node[minimum size=2em, inner sep=0pt] (c) {};
                \draw[red, very thick] ([xshift=-0.8em,yshift=0.6ex]c.center) -- ([xshift=0.8em,yshift=0.6ex]c.center);
            } & *(green) \tikz[baseline=(c.center)]{
                \node[minimum size=2em, inner sep=0pt] (c) {};
                \draw[red, very thick] ([xshift=-0.8em,yshift=0.6ex]c.center) -- ([xshift=0.8em,yshift=0.6ex]c.center);
            }\\
            *(green) \tikz[baseline=(c.center)]{
                \node[minimum size=2em, inner sep=0pt] (c) {};
                \draw[red, very thick] ([xshift=-0.8em,yshift=0.6ex]c.center) -- ([xshift=0.8em,yshift=0.6ex]c.center);
            } & *(blue) \tikz[baseline=(c.center)]{
                \node[minimum size=2em, inner sep=0pt] (c) {};
                \draw[red, very thick] ([xshift=-0.8em,yshift=0.6ex]c.center) -- ([xshift=0.8em,yshift=0.6ex]c.center);
            }\\
            \vdots & \vdots \\
            *(blue) \tikz[baseline=(c.center)]{
                \node[minimum size=2em, inner sep=0pt] (c) {};
                \draw[red, very thick] ([xshift=-0.8em,yshift=0.6ex]c.center) -- ([xshift=0.8em,yshift=0.6ex]c.center);
            } & *(green) \tikz[baseline=(c.center)]{
                \node[minimum size=2em, inner sep=0pt] (c) {};
                \draw[red, very thick] ([xshift=-0.8em,yshift=0.6ex]c.center) -- ([xshift=0.8em,yshift=0.6ex]c.center);
            } \\
            *(green) \tikz[baseline=(c.center)]{
                \node[minimum size=2em, inner sep=0pt] (c) {};
                \draw[red, very thick] ([xshift=-0.8em,yshift=0.6ex]c.center) -- ([xshift=0.8em,yshift=0.6ex]c.center);
            } \\
            \vdots
        \end{ytableau}
        }

        \medskip
        {\small (2j)}
    \end{minipage}

    \caption{Inductive steps for cases (2h), (2i), and (2j).}
    \label{fig:delta_iso_2hij}
\end{figure}

\item Assume that \( T_i \) and \( B_{i+1} \) are an \( a \)-cell and
  a blue \( b \)-cell, respectively, illustrated in
  Figure~\ref{fig:delta_iso_2kl}~(2k). Then \( T_{i+1} \) is a blue
  \(b\)-cell, and the cell to the left of \( B_{i+1} \) is an
  \(a\)-cell. We have
  \[
    \type(\pi_\circ)=s00, \qand \type(\pi)=s10,
  \]
  for some binary word \(s\).

  An odd number of strips traverse the \(i\)th and the \((i+1)\)th
  column horizontally. All other strip endings in \(\pi_\circ\) are
  unchanged. Note that the topmost strip in both \(\pi\) and
  \(\pi_\circ\) starts with an \(a\)-cell, so the first row of both
  \(\lambda^{(1)}\) and \(\lambda_\circ^{(1)}\) ends in green. Therefore,
  \[
    \type(\lambda_\circ^{(1)},\lambda_\circ^{(2)})=r\,0101\cdots010\,0\sim r00, \qand \type(\lambda^{(1)},\lambda^{(2)})=r\,1010\cdots101\,0\sim r10,
  \]
  for some binary word \(r\).

\item Assume that \( T_i \) and \( B_{i+1} \) are a \( c \)-cell and
  a green \( b \)-cell, respectively, illustrated in
  Figure~\ref{fig:delta_iso_2kl}~(2l). Then \( T_{i+1} \) is a green
  \(b\)-cell, and the cell to the left of \( B_{i+1} \) is a
  \(c\)-cell. We have
  \[
    \type(\pi_\circ)=s, \qand \type(\pi)=s00,
  \]
  for some binary word \(s\).

  An odd number of strips traverse the \(i\)th and the \((i+1)\)th
  column horizontally. All other strip endings in \(\pi_\circ\) are
  unchanged. Note that the topmost strip in both \(\pi\) and
  \(\pi_\circ\) starts with a \(c\)-cell, so the first row of both
  \(\lambda^{(1)}\) and \(\lambda_\circ^{(1)}\) ends in blue. By the
  difference conditions in \(\mathcal{P}_{2,b}\), there must be an
  \(a\)-cell diagonally south-west of the \(b\)-cell \( B_{i+1} \).
  Moreover, by definition of the map \(\delta\), this cell must belong
  to \(S_1\). Thus, the topmost strip ending of \(\pi_\circ\) that
  remains unchanged contributes a letter \(0\). Therefore,
  \[
    \type(\lambda_\circ^{(1)},\lambda_\circ^{(2)})=r\,0\,1010\cdots101\sim
    r, \qand \type(\lambda^{(1)},\lambda^{(2)})=r\,0\,0101\cdots010\sim r00,
  \]
 for some binary word \(r\).
\begin{figure}[h]
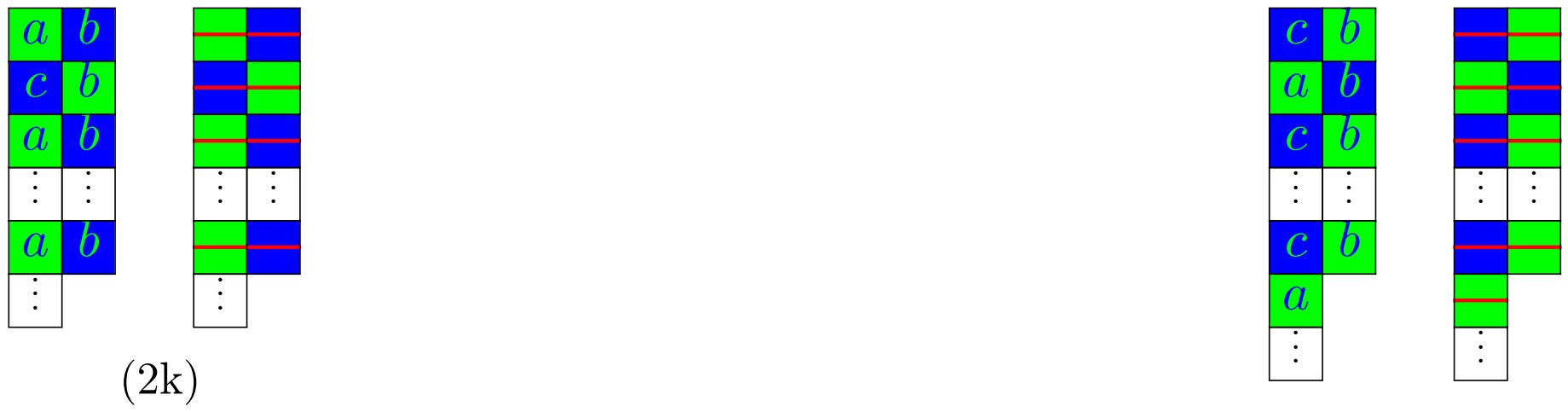

    \centering

    \begin{minipage}[t]{0.45\textwidth}
        \centering
        \scalebox{0.7}{%
        \begin{ytableau}
            *(green) \scalebox{1.5}{\color{blue}$a$} & *(blue) \scalebox{1.5}{\color{green}$b$} \\
            *(blue) \scalebox{1.5}{\color{green}$c$} & *(green) \scalebox{1.5}{\color{blue}$b$} \\
            *(green) \scalebox{1.5}{\color{blue}$a$} & *(blue) \scalebox{1.5}{\color{green}$b$} \\
            \vdots & \vdots \\
            *(green) \scalebox{1.5}{\color{blue}$a$} & *(blue) \scalebox{1.5}{\color{green}$b$} \\
            \vdots
        \end{ytableau}
        }
        \quad
        \scalebox{0.7}{%
        \begin{ytableau}
            *(green) \tikz[baseline=(c.center)]{
                \node[minimum size=2em, inner sep=0pt] (c) {};
                \draw[red, very thick] ([xshift=-0.8em,yshift=0.6ex]c.center) -- ([xshift=0.8em,yshift=0.6ex]c.center);
            } & *(blue) \tikz[baseline=(c.center)]{
                \node[minimum size=2em, inner sep=0pt] (c) {};
                \draw[red, very thick] ([xshift=-0.8em,yshift=0.6ex]c.center) -- ([xshift=0.8em,yshift=0.6ex]c.center);
            } \\
            *(blue) \tikz[baseline=(c.center)]{
                \node[minimum size=2em, inner sep=0pt] (c) {};
                \draw[red, very thick] ([xshift=-0.8em,yshift=0.6ex]c.center) -- ([xshift=0.8em,yshift=0.6ex]c.center);
            } & *(green) \tikz[baseline=(c.center)]{
                \node[minimum size=2em, inner sep=0pt] (c) {};
                \draw[red, very thick] ([xshift=-0.8em,yshift=0.6ex]c.center) -- ([xshift=0.8em,yshift=0.6ex]c.center);
            } \\
            *(green) \tikz[baseline=(c.center)]{
                \node[minimum size=2em, inner sep=0pt] (c) {};
                \draw[red, very thick] ([xshift=-0.8em,yshift=0.6ex]c.center) -- ([xshift=0.8em,yshift=0.6ex]c.center);
            } & *(blue) \tikz[baseline=(c.center)]{
                \node[minimum size=2em, inner sep=0pt] (c) {};
                \draw[red, very thick] ([xshift=-0.8em,yshift=0.6ex]c.center) -- ([xshift=0.8em,yshift=0.6ex]c.center);
            } \\
            \vdots & \vdots \\
            *(green) \tikz[baseline=(c.center)]{
                \node[minimum size=2em, inner sep=0pt] (c) {};
                \draw[red, very thick] ([xshift=-0.8em,yshift=0.6ex]c.center) -- ([xshift=0.8em,yshift=0.6ex]c.center);
            } & *(blue) \tikz[baseline=(c.center)]{
                \node[minimum size=2em, inner sep=0pt] (c) {};
                \draw[red, very thick] ([xshift=-0.8em,yshift=0.6ex]c.center) -- ([xshift=0.8em,yshift=0.6ex]c.center);
            } \\
            \vdots
        \end{ytableau}
        }

        \medskip
        {\small (2k)}
    \end{minipage}
    \hfill
    \begin{minipage}[t]{0.45\textwidth}
        \centering
        \scalebox{0.7}{%
        \begin{ytableau}
            *(blue) \scalebox{1.5}{\color{green}$c$} & *(green) \scalebox{1.5}{\color{blue}$b$} \\
            *(green) \scalebox{1.5}{\color{blue}$a$} & *(blue) \scalebox{1.5}{\color{green}$b$} \\
            *(blue) \scalebox{1.5}{\color{green}$c$} & *(green) \scalebox{1.5}{\color{blue}$b$} \\
            \vdots & \vdots \\
            *(blue) \scalebox{1.5}{\color{green}$c$} & *(green) \scalebox{1.5}{\color{blue}$b$} \\
            *(green) \scalebox{1.5}{\color{blue}$a$} \\
            \vdots
        \end{ytableau}
        }
        \quad
        \scalebox{0.7}{%
        \begin{ytableau}
            *(blue) \tikz[baseline=(c.center)]{
                \node[minimum size=2em, inner sep=0pt] (c) {};
                \draw[red, very thick] ([xshift=-0.8em,yshift=0.6ex]c.center) -- ([xshift=0.8em,yshift=0.6ex]c.center);
            } & *(green) \tikz[baseline=(c.center)]{
                \node[minimum size=2em, inner sep=0pt] (c) {};
                \draw[red, very thick] ([xshift=-0.8em,yshift=0.6ex]c.center) -- ([xshift=0.8em,yshift=0.6ex]c.center);
            } \\
            *(green) \tikz[baseline=(c.center)]{
                \node[minimum size=2em, inner sep=0pt] (c) {};
                \draw[red, very thick] ([xshift=-0.8em,yshift=0.6ex]c.center) -- ([xshift=0.8em,yshift=0.6ex]c.center);
            } & *(blue) \tikz[baseline=(c.center)]{
                \node[minimum size=2em, inner sep=0pt] (c) {};
                \draw[red, very thick] ([xshift=-0.8em,yshift=0.6ex]c.center) -- ([xshift=0.8em,yshift=0.6ex]c.center);
            } \\
            *(blue) \tikz[baseline=(c.center)]{
                \node[minimum size=2em, inner sep=0pt] (c) {};
                \draw[red, very thick] ([xshift=-0.8em,yshift=0.6ex]c.center) -- ([xshift=0.8em,yshift=0.6ex]c.center);
            } & *(green) \tikz[baseline=(c.center)]{
                \node[minimum size=2em, inner sep=0pt] (c) {};
                \draw[red, very thick] ([xshift=-0.8em,yshift=0.6ex]c.center) -- ([xshift=0.8em,yshift=0.6ex]c.center);
            } \\
            \vdots & \vdots \\
            *(blue) \tikz[baseline=(c.center)]{
                \node[minimum size=2em, inner sep=0pt] (c) {};
                \draw[red, very thick] ([xshift=-0.8em,yshift=0.6ex]c.center) -- ([xshift=0.8em,yshift=0.6ex]c.center);
            } & *(green) \tikz[baseline=(c.center)]{
                \node[minimum size=2em, inner sep=0pt] (c) {};
                \draw[red, very thick] ([xshift=-0.8em,yshift=0.6ex]c.center) -- ([xshift=0.8em,yshift=0.6ex]c.center);
            } \\
            *(green) \tikz[baseline=(c.center)]{
                \node[minimum size=2em, inner sep=0pt] (c) {};
                \draw[red, very thick] ([xshift=-0.8em,yshift=0.6ex]c.center) -- ([xshift=0.8em,yshift=0.6ex]c.center);
            } \\
            \vdots
        \end{ytableau}
        }

        \medskip
        {\small (2l)}
    \end{minipage}

    \caption{Inductive steps for cases (2k) and (2l).}
    \label{fig:delta_iso_2kl}
\end{figure} \qedhere
\end{enumerate}
\end{proof}

\begin{theorem}\label{thm:delta_iso}
  The graphs \(\mathcal{B}_{2,b}^{\mathcal{P}}\) and
  \(\mathcal{B}_{2,b}^{\mathcal{A}}\) are isomorphic under the map \( \delta \).
\end{theorem}

\begin{proof}
\emph{$\delta$ is a morphism of edge-coloured directed graphs.} Let $\pi\in\mathcal{P}_{2,b}$. By Lemma~\ref{lemma:delta_type}, $\type(\pi)\sim\type(\delta(\pi))$, so a blue arrow starts at $\pi$ if and only if a blue arrow starts at $\delta(\pi)$. In this case, let $\pi\to\pi'$ be the blue arrow from $\pi$, and let $u$ be the cell that it adds. We have to show that the blue arrow from $\delta(\pi)$ ends at $\delta(\pi')$. For green arrows, the proof is analogous.

Adding $u$ does not change the membership of any existing cell in $S_1$ or $S_2$. If $u\in S_1$, then $u$ is added at the end of a strip, or forms a new strip of length $1$; the other strips are unchanged. Thus $\delta(\pi')$ is obtained from $\delta(\pi)$ by adding one cell. Write $\type(\pi)=x0y$, where the $0$ is the letter corresponding to $u$, and $\type(\delta(\pi))=z0w$, where the $0$ is the letter corresponding to this cell. We will show that $x\sim z$ and $y\sim w$. This implies that if the $0$ of $x0y$ is its leftmost unmatched $0$, then the $0$ of $z0w$ is its leftmost unmatched $0$.

As the blue arrow from $\pi$ adds $u$, the $0$ of $x0y$ is its leftmost unmatched $0$. Thus $x$ has no unmatched $0$. If $u\in S_2$, then $x=z$. Assume that $u\in S_1$. The words $x$ and $z$ do not depend on the columns to the right of $u$. Removing these columns from $\pi$ and $\pi'$ gives two grounded partitions. Their types are $x0$ and $x1$ followed only by appended $0$s. In the second one, $u$ is in the last column, so the strips above the strip of $u$ end at the cells above $u$. Thus the types of their images under $\delta$ are $z0$ and $z1$ followed only by $0101\cdots$ and possibly an appended $0$. By Lemma~\ref{lemma:delta_type},
\begin{align}
x\,0\,0\cdots0&\sim z\,0\,0101\cdots,\label{eq:delta_iso_0}\\
x\,1\,0\cdots0&\sim z\,1\,0101\cdots.\label{eq:delta_iso_1}
\end{align}
The $1$s of $0101\cdots$ are matched. Thus, by \eqref{eq:delta_iso_0}, $x$ and $z$ have the same number of unmatched $1$s. By \eqref{eq:delta_iso_1}, $x1$ and $z1$ have the same number of unmatched $1$s. As $x$ has no unmatched $0$, $x1$ has one more unmatched $1$ than $x$. Thus $z1$ has one more unmatched $1$ than $z$, so $z$ has no unmatched $0$. As $x$ and $z$ have no unmatched $0$ and $x$ and $z$ have the same number of unmatched $1$s, we have $x\sim z$. Therefore
\begin{align}
x1y=\type(\pi')
&\sim\type(\delta(\pi'))=z1w
\qquad{\scriptstyle\textcolor{gray!60}{\text{by Lemma~\ref{lemma:delta_type}}}}\\[0.5em]
&\sim x1w.
\qquad{\scriptstyle\textcolor{gray!60}{\text{as $z\sim x$}}}
\end{align}
As $x1$ has no unmatched $0$, this gives $y\sim w$.

By Proposition~\ref{prop:vertices_P2_crystal} and induction along the paths from $\emptyset$, $\delta$ maps $\mathcal{P}_{2,b}$ to the vertices of $\mathcal{B}_{2,b}^{\mathcal{A}}$.

\emph{$\delta$ is surjective.} Let $y$ be a vertex of $\mathcal{B}_{2,b}^{\mathcal{A}}$. By Definition~\ref{def:pairs_crystal}, there is a path $(\emptyset,\emptyset)=y_0\to y_1\to\cdots\to y_r=y$. We show by induction on $j$ that $y_j=\delta(v_j)$ for a vertex $v_j$ of $\mathcal{B}_{2,b}^{\mathcal{P}}$; for $j=0$ take $v_0=\emptyset$. Let $y_j=\delta(v_j)$, and let $y_j\to y_{j+1}$ be a blue arrow. Then $\type(\delta(v_j))$ has an unmatched $0$, so $\type(v_j)$ has an unmatched $0$ by Lemma~\ref{lemma:delta_type}, and $v_j$ has a blue arrow $v_j\to v_{j+1}$. As $\delta$ is a morphism, $\delta(v_j)\to\delta(v_{j+1})$ is a blue arrow. As $\delta(v_j)$ has at most one outgoing blue arrow, $\delta(v_{j+1})=y_{j+1}$. For green arrows, the proof is analogous.

\emph{$\delta$ is injective.} Let $\pi,\tau\in\mathcal{P}_{2,b}$ with $\delta(\pi)=\delta(\tau)$. We show that $\pi=\tau$ by induction on the size, which $\delta$ preserves. If $\pi$ is empty, so is $\tau$. Otherwise, $\delta(\pi)$ is the end of an arrow, say a blue one, so $\type(\delta(\pi))$ has an unmatched $1$. We have
\begin{align}
\type(\pi)\sim\type(\delta(\pi))=\type(\delta(\tau))\sim\type(\tau).
\qquad{\scriptstyle\textcolor{gray!60}{\text{by Lemma~\ref{lemma:delta_type}}}}
\end{align}
Thus $\type(\pi)$ and $\type(\tau)$ have an unmatched $1$. As shown in the proof of Proposition~\ref{prop:vertices_P2_crystal}, $\pi$ and $\tau$ are then the ends of blue arrows $\pi_0\to\pi$ and $\tau_0\to\tau$. As $\delta$ is a morphism of edge-coloured directed graphs, $\delta(\pi_0)\to\delta(\pi)$ and $\delta(\tau_0)\to\delta(\tau)$ are blue arrows. As only one blue arrow ends at $\delta(\pi)=\delta(\tau)$, we have $\delta(\pi_0)=\delta(\tau_0)$. By induction, $\pi_0=\tau_0$, so $\pi=\tau$.
\end{proof}

As a direct consequence of Theorem~\ref{thm:delta_iso}, we have the main result of this section.

\begin{theorem}\label{thm:grounded_crystal_b}
  The graph \(\mathcal{B}_{2,b}^{\mathcal{P}}\) is the affine crystal graph
  of type \(A_1^{(1)}\) and highest weight \(\Lambda_0+\Lambda_1\).
\end{theorem}

We now construct the analogous grounded partition models for the affine crystal graphs of highest weights \( 2\Lambda_0 \) and \( 2\Lambda_1 \). We interchange the roles of even-indexed columns, including the \( 0 \)th column, and odd-indexed columns.

\begin{definition}\label{def:P2_crystal2}
Define \(\mathcal{B}_{2,ac}^{\mathcal{P}}\) and \(\mathcal{B}_{2,ca}^{\mathcal{P}}\) by the rules for the arrows in Definition~\ref{def:P2_crystal}. Keep the colouring and type rules of Definition~\ref{def:P2_type}, but replace its condition on the entries by the following:
\begin{enumerate}
\item[(i\( ' \))] All cells in columns with even indices have entry
  \( a \) or \( c \), while all cells in columns with odd indices have
  entry \( b \). As an initial condition, the cells in the \(0\)th
  column from top to bottom are \(a,c,a,c,\ldots\) (respectively,
  \(c,a,c,a,\ldots\)).
\end{enumerate}

\end{definition}

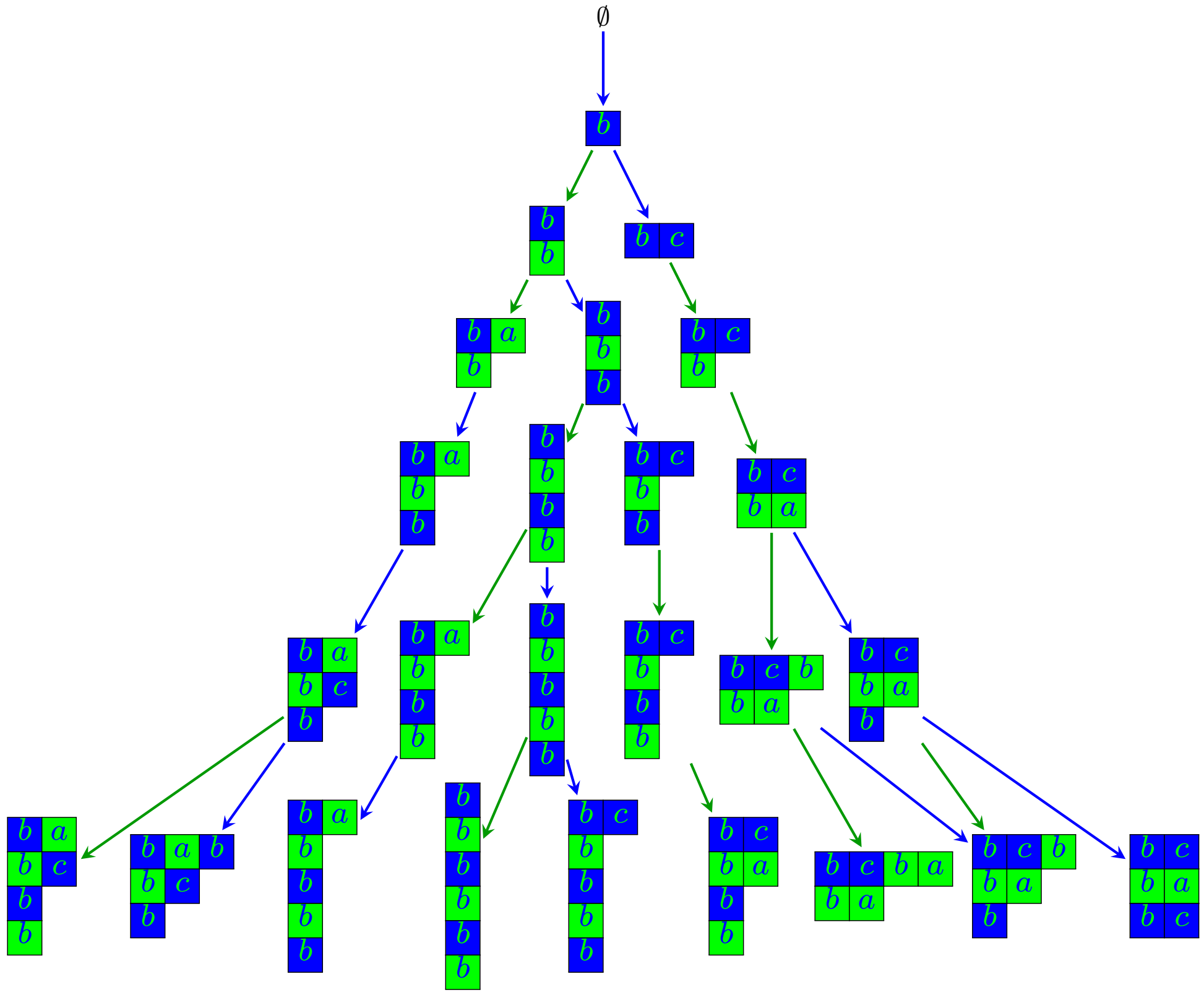
\begin{figure}
\centering
{\color{black}\begin{tikzpicture}[>=stealth, thick, scale=0.6]
  %%%%%%%%%%%%%%
  % Coordinates (unchanged)
  %%%%%%%%%%%%%%
  % Row 1: 1 vertex
  \coordinate (R1C1) at (0,0);
  
  % Row 2: 1 vertex
  \coordinate (R2C1) at (0,-2);
  
  % Row 3: 2 vertex
  \coordinate (R3C1) at (-1,-4);
  \coordinate (R3C2) at (1,-4);
  
  % Row 4: 3 vertices
  \coordinate (R4C1) at (-2,-6);
  \coordinate (R4C2) at (0,-6);
  \coordinate (R4C3) at (2,-6);
  
  % Row 5: 4 vertices
  \coordinate (R5C1) at (-3,-8.5);
  \coordinate (R5C2) at (-1,-8.5);
  \coordinate (R5C3) at (1,-8.5);
  \coordinate (R5C4) at (3,-8.5);

  % Row 6: 6 vertices
  \coordinate (R6C1) at (-5,-12);
  \coordinate (R6C2) at (-3,-12);
  \coordinate (R6C3) at (-1,-12);
  \coordinate (R6C4) at (1,-12);
  \coordinate (R6C5) at (3,-12);
  \coordinate (R6C6) at (5,-12);

  % Row 6: 9 vertices
  \coordinate (R7C1) at (-10,-15.5);
  \coordinate (R7C2) at (-7.5,-15.5);
  \coordinate (R7C3) at (-5,-15.5);
  \coordinate (R7C4) at (-2.5,-15.5);
  \coordinate (R7C5) at (0,-15.5);
  \coordinate (R7C6) at (2.5,-15.5);
  \coordinate (R7C7) at (5,-15.5);
  \coordinate (R7C8) at (7.5,-15.5);
  \coordinate (R7C9) at (10,-15.5);
  
  %%%%%%%%%%%%%%
  % Nodes (transposed shapes)
  %%%%%%%%%%%%%%
  \tikzset{mynode/.style={inner sep=0pt, font=\scriptsize}}
  
  \node[mynode] (N1C1) at (R1C1) {\(\emptyset\)};
  
  \node[mynode] (N2C1) at (R2C1) {\(\scalebox{0.8}{\begin{ytableau}
*(blue) \scalebox{1.5}{\color{green}$b$}
\end{ytableau}}\)};
  
  \node[mynode] (N3C1) at (R3C1) {\(\scalebox{0.8}{\begin{ytableau}
*(blue) \scalebox{1.5}{\color{green}$b$} \\ *(green) \scalebox{1.5}{\color{blue}$b$}
\end{ytableau}}\)};

  \node[mynode] (N3C2) at (R3C2) {\(\scalebox{0.8}{\begin{ytableau}
*(blue) \scalebox{1.5}{\color{green}$b$} &
*(blue) \scalebox{1.5}{\color{green}$c$}
\end{ytableau}}\)};
  
  \node[mynode] (N4C1) at (R4C1) {\(\scalebox{0.8}{\begin{ytableau}
*(blue) \scalebox{1.5}{\color{green}$b$} & *(green) \scalebox{1.5}{\color{blue}$a$} \\
*(green) \scalebox{1.5}{\color{blue}$b$}
\end{ytableau}}\)};

  \node[mynode] (N4C2) at (R4C2) {\(\scalebox{0.8}{\begin{ytableau}
*(blue) \scalebox{1.5}{\color{green}$b$} \\
*(green) \scalebox{1.5}{\color{blue}$b$} \\
*(blue) \scalebox{1.5}{\color{green}$b$}
\end{ytableau}}\)};

  \node[mynode] (N4C3) at (R4C3) {\(\scalebox{0.8}{\begin{ytableau}
*(blue) \scalebox{1.5}{\color{green}$b$} & *(blue) \scalebox{1.5}{\color{green}$c$} \\
*(green) \scalebox{1.5}{\color{blue}$b$}
\end{ytableau}}\)};

  \node[mynode] (N5C1) at (R5C1) {\(\scalebox{0.8}{\begin{ytableau}
*(blue) \scalebox{1.5}{\color{green}$b$} & *(green) \scalebox{1.5}{\color{blue}$a$} \\
*(green) \scalebox{1.5}{\color{blue}$b$} \\
*(blue) \scalebox{1.5}{\color{green}$b$}
\end{ytableau}}\)};

  \node[mynode] (N5C2) at (R5C2) {\(\scalebox{0.8}{\begin{ytableau}
*(blue) \scalebox{1.5}{\color{green}$b$} \\
*(green) \scalebox{1.5}{\color{blue}$b$} \\
*(blue) \scalebox{1.5}{\color{green}$b$} \\
*(green) \scalebox{1.5}{\color{blue}$b$}
\end{ytableau}}\)};

  \node[mynode] (N5C3) at (R5C3) {\(\scalebox{0.8}{\begin{ytableau}
*(blue) \scalebox{1.5}{\color{green}$b$} & *(blue) \scalebox{1.5}{\color{green}$c$} \\
*(green) \scalebox{1.5}{\color{blue}$b$} \\
*(blue) \scalebox{1.5}{\color{green}$b$}
\end{ytableau}}\)};

  \node[mynode] (N5C4) at (R5C4) {\(\scalebox{0.8}{\begin{ytableau}
*(blue) \scalebox{1.5}{\color{green}$b$} & *(blue) \scalebox{1.5}{\color{green}$c$} \\
*(green) \scalebox{1.5}{\color{blue}$b$} & *(green) \scalebox{1.5}{\color{blue}$a$}
\end{ytableau}}\)};

  \node[mynode] (N6C1) at (R6C1) {\(\scalebox{0.8}{\begin{ytableau}
*(blue) \scalebox{1.5}{\color{green}$b$} & *(green) \scalebox{1.5}{\color{blue}$a$} \\
*(green) \scalebox{1.5}{\color{blue}$b$} & *(blue) \scalebox{1.5}{\color{green}$c$} \\
*(blue) \scalebox{1.5}{\color{green}$b$}
\end{ytableau}}\)};

  \node[mynode] (N6C2) at (R6C2) {\(\scalebox{0.8}{\begin{ytableau}
*(blue) \scalebox{1.5}{\color{green}$b$} & *(green) \scalebox{1.5}{\color{blue}$a$} \\
*(green) \scalebox{1.5}{\color{blue}$b$} \\
*(blue) \scalebox{1.5}{\color{green}$b$} \\
*(green) \scalebox{1.5}{\color{blue}$b$}
\end{ytableau}}\)};

  \node[mynode] (N6C3) at (R6C3) {\(\scalebox{0.8}{\begin{ytableau}
*(blue) \scalebox{1.5}{\color{green}$b$} \\
*(green) \scalebox{1.5}{\color{blue}$b$} \\
*(blue) \scalebox{1.5}{\color{green}$b$} \\
*(green) \scalebox{1.5}{\color{blue}$b$} \\
*(blue) \scalebox{1.5}{\color{green}$b$}
\end{ytableau}}\)};

  \node[mynode] (N6C4) at (R6C4) {\(\scalebox{0.8}{\begin{ytableau}
*(blue) \scalebox{1.5}{\color{green}$b$} & *(blue) \scalebox{1.5}{\color{green}$c$} \\
*(green) \scalebox{1.5}{\color{blue}$b$} \\
*(blue) \scalebox{1.5}{\color{green}$b$} \\
*(green) \scalebox{1.5}{\color{blue}$b$}
\end{ytableau}}\)};

  \node[mynode] (N6C5) at (R6C5) {\(\scalebox{0.8}{\begin{ytableau}
*(blue) \scalebox{1.5}{\color{green}$b$} & *(blue) \scalebox{1.5}{\color{green}$c$} & *(green) \scalebox{1.5}{\color{blue}$b$} \\
*(green) \scalebox{1.5}{\color{blue}$b$} & *(green) \scalebox{1.5}{\color{blue}$a$}
\end{ytableau}}\)};

  \node[mynode] (N6C6) at (R6C6) {\(\scalebox{0.8}{\begin{ytableau}
*(blue) \scalebox{1.5}{\color{green}$b$} & *(blue) \scalebox{1.5}{\color{green}$c$} \\
*(green) \scalebox{1.5}{\color{blue}$b$} & *(green) \scalebox{1.5}{\color{blue}$a$} \\
*(blue) \scalebox{1.5}{\color{green}$b$}
\end{ytableau}}\)};

  \node[mynode] (N7C1) at (R7C1) {\(\scalebox{0.8}{\begin{ytableau}
*(blue) \scalebox{1.5}{\color{green}$b$} & *(green) \scalebox{1.5}{\color{blue}$a$} \\
*(green) \scalebox{1.5}{\color{blue}$b$} & *(blue) \scalebox{1.5}{\color{green}$c$} \\
*(blue) \scalebox{1.5}{\color{green}$b$} \\
*(green) \scalebox{1.5}{\color{blue}$b$}
\end{ytableau}}\)};

  \node[mynode] (N7C2) at (R7C2) {\(\scalebox{0.8}{\begin{ytableau}
*(blue) \scalebox{1.5}{\color{green}$b$} & *(green) \scalebox{1.5}{\color{blue}$a$} & *(blue) \scalebox{1.5}{\color{green}$b$} \\
*(green) \scalebox{1.5}{\color{blue}$b$} & *(blue) \scalebox{1.5}{\color{green}$c$} \\
*(blue) \scalebox{1.5}{\color{green}$b$}
\end{ytableau}}\)};

  \node[mynode] (N7C3) at (R7C3) {\(\scalebox{0.8}{\begin{ytableau}
*(blue) \scalebox{1.5}{\color{green}$b$} & *(green) \scalebox{1.5}{\color{blue}$a$} \\
*(green) \scalebox{1.5}{\color{blue}$b$} \\
*(blue) \scalebox{1.5}{\color{green}$b$} \\
*(green) \scalebox{1.5}{\color{blue}$b$} \\
*(blue) \scalebox{1.5}{\color{green}$b$}
\end{ytableau}}\)};

  \node[mynode] (N7C4) at (R7C4) {\(\scalebox{0.8}{\begin{ytableau}
*(blue) \scalebox{1.5}{\color{green}$b$} \\
*(green) \scalebox{1.5}{\color{blue}$b$} \\
*(blue) \scalebox{1.5}{\color{green}$b$} \\
*(green) \scalebox{1.5}{\color{blue}$b$} \\
*(blue) \scalebox{1.5}{\color{green}$b$} \\
*(green) \scalebox{1.5}{\color{blue}$b$}
\end{ytableau}}\)};

  \node[mynode] (N7C5) at (R7C5) {\(\scalebox{0.8}{\begin{ytableau}
*(blue) \scalebox{1.5}{\color{green}$b$} & *(blue) \scalebox{1.5}{\color{green}$c$}\\
*(green) \scalebox{1.5}{\color{blue}$b$} \\
*(blue) \scalebox{1.5}{\color{green}$b$} \\
*(green) \scalebox{1.5}{\color{blue}$b$} \\
*(blue) \scalebox{1.5}{\color{green}$b$}
\end{ytableau}}\)};

  \node[mynode] (N7C6) at (R7C6) {\(\scalebox{0.8}{\begin{ytableau}
*(blue) \scalebox{1.5}{\color{green}$b$} & *(blue) \scalebox{1.5}{\color{green}$c$} \\
*(green) \scalebox{1.5}{\color{blue}$b$} & *(green) \scalebox{1.5}{\color{blue}$a$} \\
*(blue) \scalebox{1.5}{\color{green}$b$} \\
*(green) \scalebox{1.5}{\color{blue}$b$}
\end{ytableau}}\)};

  \node[mynode] (N7C7) at (R7C7) {\(\scalebox{0.8}{\begin{ytableau}
*(blue) \scalebox{1.5}{\color{green}$b$} & *(blue) \scalebox{1.5}{\color{green}$c$} & *(green) \scalebox{1.5}{\color{blue}$b$} & *(green) \scalebox{1.5}{\color{blue}$a$} \\
*(green) \scalebox{1.5}{\color{blue}$b$} & *(green) \scalebox{1.5}{\color{blue}$a$}
\end{ytableau}}\)};

  \node[mynode] (N7C8) at (R7C8) {\(\scalebox{0.8}{\begin{ytableau}
*(blue) \scalebox{1.5}{\color{green}$b$} & *(blue) \scalebox{1.5}{\color{green}$c$} & *(green) \scalebox{1.5}{\color{blue}$b$} \\
*(green) \scalebox{1.5}{\color{blue}$b$} & *(green) \scalebox{1.5}{\color{blue}$a$} \\
*(blue) \scalebox{1.5}{\color{green}$b$}
\end{ytableau}}\)};

  \node[mynode] (N7C9) at (R7C9) {\(\scalebox{0.8}{\begin{ytableau}
*(blue) \scalebox{1.5}{\color{green}$b$} & *(blue) \scalebox{1.5}{\color{green}$c$} \\
*(green) \scalebox{1.5}{\color{blue}$b$} & *(green) \scalebox{1.5}{\color{blue}$a$} \\
*(blue) \scalebox{1.5}{\color{green}$b$} & *(blue) \scalebox{1.5}{\color{green}$c$}
\end{ytableau}}\)};
  
  %%%%%%%%%%%%%%
  % Arrows (unchanged)
  %%%%%%%%%%%%%%
  \tikzset{arrow style/.style={->, shorten <=1pt, shorten >=1pt, line width=0.8pt}}
  \draw[arrow style, blue] (N1C1) -- (N2C1);
  
  \draw[arrow style, green!60!black] (N2C1) -- (N3C1);
  \draw[arrow style, blue] (N2C1) -- (N3C2);
  
  \draw[arrow style, green!60!black] (N3C1) -- (N4C1);
  \draw[arrow style, blue] (N3C1) -- (N4C2);
  \draw[arrow style, green!60!black] (N3C2) -- (N4C3);
  
  \draw[arrow style, blue] (N4C1) -- (N5C1);
  \draw[arrow style, green!60!black] (N4C2) -- (N5C2);
  \draw[arrow style, blue] (N4C2) -- (N5C3);
  \draw[arrow style, green!60!black] (N4C3) -- (N5C4);

  \draw[arrow style, blue] (N5C1) -- (N6C1);
  \draw[arrow style, green!60!black] (N5C2) -- (N6C2);
  \draw[arrow style, blue] (N5C2) -- (N6C3);
  \draw[arrow style, green!60!black] (N5C3) -- (N6C4);
  \draw[arrow style, green!60!black] (N5C4) -- (N6C5);
  \draw[arrow style, blue] (N5C4) -- (N6C6);

  \draw[arrow style, green!60!black] (N6C1) -- (N7C1);
  \draw[arrow style, blue] (N6C1) -- (N7C2);
  \draw[arrow style, blue] (N6C2) -- (N7C3);
  \draw[arrow style, green!60!black] (N6C3) -- (N7C4);
  \draw[arrow style, blue] (N6C3) -- (N7C5);
  \draw[arrow style, green!60!black] (N6C4) -- (N7C6);
  \draw[arrow style, green!60!black] (N6C5) -- (N7C7);
  \draw[arrow style, blue] (N6C5) -- (N7C8);
  \draw[arrow style, green!60!black] (N6C6) -- (N7C8);
  \draw[arrow style, blue] (N6C6) -- (N7C9);
\end{tikzpicture}}
\caption{The graph $\mathcal{B}_{2,ac}^{\mathcal{P}}$}
\label{fig:level2ac_grounded_crystal_diags}
\end{figure}

Figure~\ref{fig:level2ac_grounded_crystal_diags} shows \(\mathcal{B}_{2,ac}^{\mathcal{P}}\). The argument of Proposition~\ref{prop:vertices_P2_crystal} gives the following characterisations of the vertex sets \(\mathcal{P}_{2,ac}\) and \(\mathcal{P}_{2,ca}\) of \(\mathcal{B}_{2,ac}^{\mathcal{P}}\) and \(\mathcal{B}_{2,ca}^{\mathcal{P}}\):
\[
  \mathcal{P}_{2,ac} = \{ \pi \in \mathcal{P}_{2,a} \mid \text{\( \ell(\pi) \) is even}\} \sqcup \{ \pi \in \mathcal{P}_{2,c} \mid \text{\( \ell(\pi) \) is odd}\}, 
\]
and
\[
  \mathcal{P}_{2,ca} = \{ \pi \in \mathcal{P}_{2,c} \mid \text{\( \ell(\pi) \) is even}\} \sqcup \{ \pi \in \mathcal{P}_{2,a} \mid \text{\( \ell(\pi) \) is odd}\},
\]
respectively.

\begin{remark}
There is a size-preserving bijection from $\mathcal{P}_{2,a}$ to $\mathcal{P}_{2,ac}$, given by swapping the labels $a$ and $c$ when the length is odd, and by the identity when the length is even. The same rule gives a bijection from $\mathcal{P}_{2,c}$ to $\mathcal{P}_{2,ca}$. However, for the graph isomorphism, it is essential to have alternating entries $a$ and $c$ in the $0$th column.
\end{remark}

The asymptotic model $\mathcal{B}_{2,ac}^{\mathcal{A}}$ (respectively,
$\mathcal{B}_{2,ca}^{\mathcal{A}}$) of highest weight $2\Lambda_1$
(respectively, $2\Lambda_0$) is defined in the same way as
$\mathcal{B}_{2,b}^{\mathcal{A}}$, except that both partitions
$\lambda^{(1)}$ and $\lambda^{(2)}$ have a blue (respectively, green) cell in
the top-left corner. The graph $\mathcal{B}_{2,ca}^{\mathcal{A}}$ is
obtained from $\mathcal{B}_{2,ac}^{\mathcal{A}}$ by interchanging
green and blue cells as well as green and blue arrows.

Define $\delta$ on the vertices of $\mathcal{B}_{2,ac}^{\mathcal{P}}$ and $\mathcal{B}_{2,ca}^{\mathcal{P}}$ as follows. Starting in the first row, put cells in $S_2$ as long as blue and green alternate. When a cell repeats the colour of its left neighbour, put it and the cells to its right in $S_1$, and continue the alternating strip at the cell below it. Repeat until the strip reaches the edge of the diagram. The cells on or below this strip form $S_2$, and those above it form $S_1$. Construct $\lambda^{(1)}$ and $\lambda^{(2)}$ from these sets as in the case of $\mathcal{B}_{2,b}^{\mathcal{P}}$.

\begin{figure}[h]
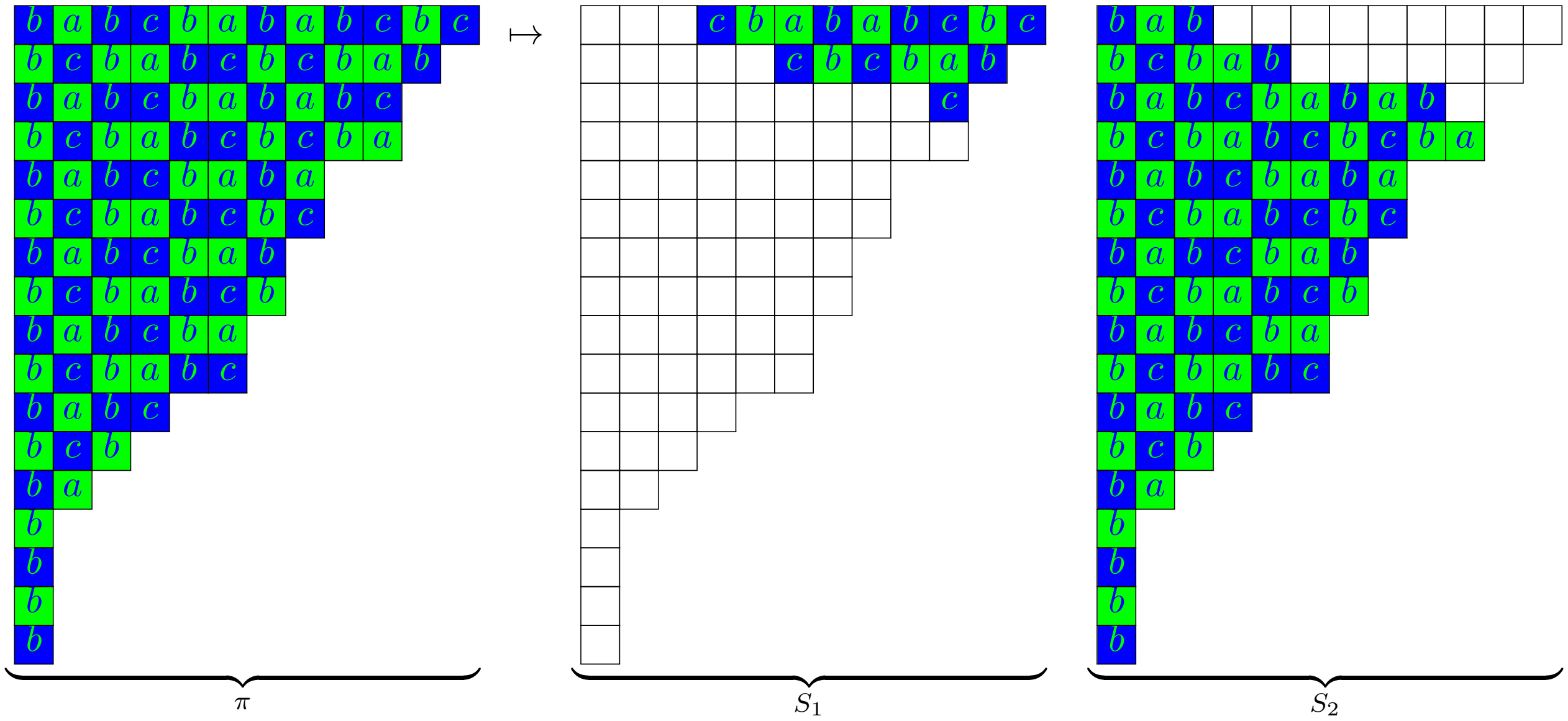

\centering
\[
\underbrace{\scalebox{0.7}{
\begin{ytableau}
*(blue) \scalebox{1.5}{\color{green}$b$}
& *(green) \scalebox{1.5}{\color{blue}$a$}
& *(blue) \scalebox{1.5}{\color{green}$b$}
& *(blue) \scalebox{1.5}{\color{green}$c$}
& *(green) \scalebox{1.5}{\color{blue}$b$}
& *(green) \scalebox{1.5}{\color{blue}$a$}
& *(blue) \scalebox{1.5}{\color{green}$b$}
& *(green) \scalebox{1.5}{\color{blue}$a$}
& *(blue) \scalebox{1.5}{\color{green}$b$}
& *(blue) \scalebox{1.5}{\color{green}$c$}
& *(green) \scalebox{1.5}{\color{blue}$b$}
& *(blue) \scalebox{1.5}{\color{green}$c$} \\
*(green) \scalebox{1.5}{\color{blue}$b$}
& *(blue) \scalebox{1.5}{\color{green}$c$}
& *(green) \scalebox{1.5}{\color{blue}$b$}
& *(green) \scalebox{1.5}{\color{blue}$a$}
& *(blue) \scalebox{1.5}{\color{green}$b$}
& *(blue) \scalebox{1.5}{\color{green}$c$}
& *(green) \scalebox{1.5}{\color{blue}$b$}
& *(blue) \scalebox{1.5}{\color{green}$c$}
& *(green) \scalebox{1.5}{\color{blue}$b$}
& *(green) \scalebox{1.5}{\color{blue}$a$}
& *(blue) \scalebox{1.5}{\color{green}$b$} \\
*(blue) \scalebox{1.5}{\color{green}$b$}
& *(green) \scalebox{1.5}{\color{blue}$a$}
& *(blue) \scalebox{1.5}{\color{green}$b$}
& *(blue) \scalebox{1.5}{\color{green}$c$}
& *(green) \scalebox{1.5}{\color{blue}$b$}
& *(green) \scalebox{1.5}{\color{blue}$a$}
& *(blue) \scalebox{1.5}{\color{green}$b$}
& *(green) \scalebox{1.5}{\color{blue}$a$}
& *(blue) \scalebox{1.5}{\color{green}$b$}
& *(blue) \scalebox{1.5}{\color{green}$c$} \\
*(green) \scalebox{1.5}{\color{blue}$b$}
& *(blue) \scalebox{1.5}{\color{green}$c$}
& *(green) \scalebox{1.5}{\color{blue}$b$}
& *(green) \scalebox{1.5}{\color{blue}$a$}
& *(blue) \scalebox{1.5}{\color{green}$b$}
& *(blue) \scalebox{1.5}{\color{green}$c$}
& *(green) \scalebox{1.5}{\color{blue}$b$}
& *(blue) \scalebox{1.5}{\color{green}$c$}
& *(green) \scalebox{1.5}{\color{blue}$b$}
& *(green) \scalebox{1.5}{\color{blue}$a$} \\
*(blue) \scalebox{1.5}{\color{green}$b$}
& *(green) \scalebox{1.5}{\color{blue}$a$}
& *(blue) \scalebox{1.5}{\color{green}$b$}
& *(blue) \scalebox{1.5}{\color{green}$c$}
& *(green) \scalebox{1.5}{\color{blue}$b$}
& *(green) \scalebox{1.5}{\color{blue}$a$}
& *(blue) \scalebox{1.5}{\color{green}$b$}
& *(green) \scalebox{1.5}{\color{blue}$a$} \\
*(green) \scalebox{1.5}{\color{blue}$b$}
& *(blue) \scalebox{1.5}{\color{green}$c$}
& *(green) \scalebox{1.5}{\color{blue}$b$}
& *(green) \scalebox{1.5}{\color{blue}$a$}
& *(blue) \scalebox{1.5}{\color{green}$b$}
& *(blue) \scalebox{1.5}{\color{green}$c$}
& *(green) \scalebox{1.5}{\color{blue}$b$}
& *(blue) \scalebox{1.5}{\color{green}$c$} \\
*(blue) \scalebox{1.5}{\color{green}$b$}
& *(green) \scalebox{1.5}{\color{blue}$a$}
& *(blue) \scalebox{1.5}{\color{green}$b$}
& *(blue) \scalebox{1.5}{\color{green}$c$}
& *(green) \scalebox{1.5}{\color{blue}$b$}
& *(green) \scalebox{1.5}{\color{blue}$a$}
& *(blue) \scalebox{1.5}{\color{green}$b$} \\
*(green) \scalebox{1.5}{\color{blue}$b$}
& *(blue) \scalebox{1.5}{\color{green}$c$}
& *(green) \scalebox{1.5}{\color{blue}$b$}
& *(green) \scalebox{1.5}{\color{blue}$a$}
& *(blue) \scalebox{1.5}{\color{green}$b$}
& *(blue) \scalebox{1.5}{\color{green}$c$}
& *(green) \scalebox{1.5}{\color{blue}$b$} \\
*(blue) \scalebox{1.5}{\color{green}$b$}
& *(green) \scalebox{1.5}{\color{blue}$a$}
& *(blue) \scalebox{1.5}{\color{green}$b$}
& *(blue) \scalebox{1.5}{\color{green}$c$}
& *(green) \scalebox{1.5}{\color{blue}$b$}
& *(green) \scalebox{1.5}{\color{blue}$a$} \\
*(green) \scalebox{1.5}{\color{blue}$b$}
& *(blue) \scalebox{1.5}{\color{green}$c$}
& *(green) \scalebox{1.5}{\color{blue}$b$}
& *(green) \scalebox{1.5}{\color{blue}$a$}
& *(blue) \scalebox{1.5}{\color{green}$b$}
& *(blue) \scalebox{1.5}{\color{green}$c$} \\
*(blue) \scalebox{1.5}{\color{green}$b$}
& *(green) \scalebox{1.5}{\color{blue}$a$}
& *(blue) \scalebox{1.5}{\color{green}$b$}
& *(blue) \scalebox{1.5}{\color{green}$c$}\\
*(green) \scalebox{1.5}{\color{blue}$b$}
& *(blue) \scalebox{1.5}{\color{green}$c$}
& *(green) \scalebox{1.5}{\color{blue}$b$}\\
*(blue) \scalebox{1.5}{\color{green}$b$}
& *(green) \scalebox{1.5}{\color{blue}$a$} \\
*(green) \scalebox{1.5}{\color{blue}$b$} \\
*(blue) \scalebox{1.5}{\color{green}$b$} \\
*(green) \scalebox{1.5}{\color{blue}$b$} \\
*(blue) \scalebox{1.5}{\color{green}$b$}
\end{ytableau}}}_{\pi} \hspace{0.5em}\mapsto \hspace{0.5em}
\underbrace{\scalebox{0.7}{
\begin{ytableau}
\null
& \null
& \null
& *(blue) \scalebox{1.5}{\color{green}$c$}
& *(green) \scalebox{1.5}{\color{blue}$b$}
& *(green) \scalebox{1.5}{\color{blue}$a$}
& *(blue) \scalebox{1.5}{\color{green}$b$}
& *(green) \scalebox{1.5}{\color{blue}$a$}
& *(blue) \scalebox{1.5}{\color{green}$b$}
& *(blue) \scalebox{1.5}{\color{green}$c$}
& *(green) \scalebox{1.5}{\color{blue}$b$}
& *(blue) \scalebox{1.5}{\color{green}$c$} \\
\null
& \null
& \null
& \null
& \null
& *(blue) \scalebox{1.5}{\color{green}$c$}
& *(green) \scalebox{1.5}{\color{blue}$b$}
& *(blue) \scalebox{1.5}{\color{green}$c$}
& *(green) \scalebox{1.5}{\color{blue}$b$}
& *(green) \scalebox{1.5}{\color{blue}$a$}
& *(blue) \scalebox{1.5}{\color{green}$b$} \\
\null
& \null
& \null
& \null
& \null
& \null
& \null
& \null
& \null
& *(blue) \scalebox{1.5}{\color{green}$c$} \\
\null
& \null
& \null
& \null
& \null
& \null
& \null
& \null
& \null
& \null\\
\null
& \null
& \null
& \null
& \null
& \null
& \null
& \null \\
\null
& \null
& \null
& \null
& \null
& \null
& \null
& \null \\
\null
& \null
& \null
& \null
& \null
& \null
& \null \\
\null
& \null
& \null
& \null
& \null
& \null
& \null \\
\null
& \null
& \null
& \null
& \null
& \null \\
\null
& \null
& \null
& \null
& \null
& \null \\
\null
& \null
& \null
& \null \\
\null
& \null
& \null \\
\null
& \null \\
\null \\
\null \\
\null \\
\null
\end{ytableau}}}_{S_1}\quad
\underbrace{\scalebox{0.7}{
\begin{ytableau}
*(blue) \scalebox{1.5}{\color{green}$b$}
& *(green) \scalebox{1.5}{\color{blue}$a$}
& *(blue) \scalebox{1.5}{\color{green}$b$}
& \null
& \null
& \null
& \null
& \null
& \null
& \null
& \null
& \null \\
*(green) \scalebox{1.5}{\color{blue}$b$}
& *(blue) \scalebox{1.5}{\color{green}$c$}
& *(green) \scalebox{1.5}{\color{blue}$b$}
& *(green) \scalebox{1.5}{\color{blue}$a$}
& *(blue) \scalebox{1.5}{\color{green}$b$}
& \null
& \null
& \null
& \null
& \null
& \null \\
*(blue) \scalebox{1.5}{\color{green}$b$}
& *(green) \scalebox{1.5}{\color{blue}$a$}
& *(blue) \scalebox{1.5}{\color{green}$b$}
& *(blue) \scalebox{1.5}{\color{green}$c$}
& *(green) \scalebox{1.5}{\color{blue}$b$}
& *(green) \scalebox{1.5}{\color{blue}$a$}
& *(blue) \scalebox{1.5}{\color{green}$b$}
& *(green) \scalebox{1.5}{\color{blue}$a$}
& *(blue) \scalebox{1.5}{\color{green}$b$}
& \null \\
*(green) \scalebox{1.5}{\color{blue}$b$}
& *(blue) \scalebox{1.5}{\color{green}$c$}
& *(green) \scalebox{1.5}{\color{blue}$b$}
& *(green) \scalebox{1.5}{\color{blue}$a$}
& *(blue) \scalebox{1.5}{\color{green}$b$}
& *(blue) \scalebox{1.5}{\color{green}$c$}
& *(green) \scalebox{1.5}{\color{blue}$b$}
& *(blue) \scalebox{1.5}{\color{green}$c$}
& *(green) \scalebox{1.5}{\color{blue}$b$}
& *(green) \scalebox{1.5}{\color{blue}$a$} \\
*(blue) \scalebox{1.5}{\color{green}$b$}
& *(green) \scalebox{1.5}{\color{blue}$a$}
& *(blue) \scalebox{1.5}{\color{green}$b$}
& *(blue) \scalebox{1.5}{\color{green}$c$}
& *(green) \scalebox{1.5}{\color{blue}$b$}
& *(green) \scalebox{1.5}{\color{blue}$a$}
& *(blue) \scalebox{1.5}{\color{green}$b$}
& *(green) \scalebox{1.5}{\color{blue}$a$} \\
*(green) \scalebox{1.5}{\color{blue}$b$}
& *(blue) \scalebox{1.5}{\color{green}$c$}
& *(green) \scalebox{1.5}{\color{blue}$b$}
& *(green) \scalebox{1.5}{\color{blue}$a$}
& *(blue) \scalebox{1.5}{\color{green}$b$}
& *(blue) \scalebox{1.5}{\color{green}$c$}
& *(green) \scalebox{1.5}{\color{blue}$b$}
& *(blue) \scalebox{1.5}{\color{green}$c$} \\
*(blue) \scalebox{1.5}{\color{green}$b$}
& *(green) \scalebox{1.5}{\color{blue}$a$}
& *(blue) \scalebox{1.5}{\color{green}$b$}
& *(blue) \scalebox{1.5}{\color{green}$c$}
& *(green) \scalebox{1.5}{\color{blue}$b$}
& *(green) \scalebox{1.5}{\color{blue}$a$}
& *(blue) \scalebox{1.5}{\color{green}$b$} \\
*(green) \scalebox{1.5}{\color{blue}$b$}
& *(blue) \scalebox{1.5}{\color{green}$c$}
& *(green) \scalebox{1.5}{\color{blue}$b$}
& *(green) \scalebox{1.5}{\color{blue}$a$}
& *(blue) \scalebox{1.5}{\color{green}$b$}
& *(blue) \scalebox{1.5}{\color{green}$c$}
& *(green) \scalebox{1.5}{\color{blue}$b$} \\
*(blue) \scalebox{1.5}{\color{green}$b$}
& *(green) \scalebox{1.5}{\color{blue}$a$}
& *(blue) \scalebox{1.5}{\color{green}$b$}
& *(blue) \scalebox{1.5}{\color{green}$c$}
& *(green) \scalebox{1.5}{\color{blue}$b$}
& *(green) \scalebox{1.5}{\color{blue}$a$} \\
*(green) \scalebox{1.5}{\color{blue}$b$}
& *(blue) \scalebox{1.5}{\color{green}$c$}
& *(green) \scalebox{1.5}{\color{blue}$b$}
& *(green) \scalebox{1.5}{\color{blue}$a$}
& *(blue) \scalebox{1.5}{\color{green}$b$}
& *(blue) \scalebox{1.5}{\color{green}$c$} \\
*(blue) \scalebox{1.5}{\color{green}$b$}
& *(green) \scalebox{1.5}{\color{blue}$a$}
& *(blue) \scalebox{1.5}{\color{green}$b$}
& *(blue) \scalebox{1.5}{\color{green}$c$}\\
*(green) \scalebox{1.5}{\color{blue}$b$}
& *(blue) \scalebox{1.5}{\color{green}$c$}
& *(green) \scalebox{1.5}{\color{blue}$b$}\\
*(blue) \scalebox{1.5}{\color{green}$b$}
& *(green) \scalebox{1.5}{\color{blue}$a$} \\
*(green) \scalebox{1.5}{\color{blue}$b$} \\
*(blue) \scalebox{1.5}{\color{green}$b$} \\
*(green) \scalebox{1.5}{\color{blue}$b$} \\
*(blue) \scalebox{1.5}{\color{green}$b$}
\end{ytableau}}}_{S_2}
\]

\caption{Decomposing $\pi$ into $S_1$ and $S_2$ in $\mathcal{P}_{2,ac}$}
\label{fig:delta_ac}
\end{figure}

For a vertex $\pi$ of $\mathcal{B}_{2,ac}^{\mathcal{P}}$, this strip consists of alternating $a$-cells and blue $b$-cells, as shown in Figure~\ref{fig:delta_ac}. For a vertex of $\mathcal{B}_{2,ca}^{\mathcal{P}}$, it consists of alternating $c$-cells and green $b$-cells.

The proofs of Lemmas~\ref{lemma:delta_well-def} and~\ref{lemma:delta_type} and of Theorem~\ref{thm:delta_iso} carry over, so $\delta$ is a graph isomorphism in both cases. We therefore obtain the following result.

\begin{theorem}\label{thm:grounded_crystal_ac}
The graph $\mathcal{B}_{2,ac}^{\mathcal{P}}$ (respectively, $\mathcal{B}_{2,ca}^{\mathcal{P}}$) is the affine crystal graph of type $A_1^{(1)}$ and highest weight $2\Lambda_1$ (respectively, $2\Lambda_0$).
\end{theorem}

\section{Decomposition formulas at levels $1$ and $2$}\label{section:decomp}

We erase the green arrows from the affine crystal graphs in
Definitions~\ref{def:partition_model_1} and \ref{def:P2_crystal}. The
resulting graphs are the crystal graphs of the restrictions to the
Levi subalgebra, and their connected
components are finite blue strings. In this section, we prove our main
results, Theorems~\ref{thm:decomp1} and \ref{thm:decomp2}.

\subsection{Major index on Yamanouchi words}\label{subsection:major index}

In this subsection, we recall the \emph{major index} on words,
introduced by MacMahon~\cite[Art.~2]{Mac13}, and define a modification
$\widetilde{\mathrm{maj}}$ whose generating function on Yamanouchi
words is, up to a power of $q$, a $q$-Yamanouchi number. This
modification is essentially the statistic introduced by
Krattenthaler~\cite[(2.2)]{Kra89}. In
Subsection~\ref{subsection:decomp1}, we express the size of the
smallest starting point of a given type in
$\mathcal{B}_{1,b}^{\mathcal{P}}$ in terms of this statistic.

\begin{definition}
Let $w = w_1\,w_2\,\cdots\,w_k$ be a word of length \(k\) on $\{0,1,2,\ldots\}$. The \emph{descent set} of \(w\) is
    \[
      \operatorname{Des}(w)
      \coloneq
      \{\,i \mid 1\le i < k,\; w_i > w_{i+1}\}.
    \]
Similarly, the \emph{ascent set} of \(w\) is
    \[
      \operatorname{Asc}(w)
      \coloneq
      \{\,i \mid 1\le i < k,\; w_i < w_{i+1}\}.
    \]
The \emph{major index} and the \emph{ascent major index} are defined as
      \[
        \mathrm{maj}(w)
        \coloneq
        \sum_{i\in\operatorname{Des}(w)} i,\hspace{2em} \mathrm{comaj}(w)
        \coloneq
        \sum_{i\in\operatorname{Asc}(w)} i.
      \]
\end{definition}

\begin{definition}\label{def:Yam_word}
A word on the alphabet $\{0,1\}$ is called a \emph{Yamanouchi word} if the number of $0$s is at least as large as the number of $1$s in each prefix. We denote the set of Yamanouchi words with $n$ many $0$s and $k$ many $1$s by $\mathcal{Y}(n,k)$. 
\end{definition}

The number of Yamanouchi words in $\mathcal{Y}(n,k)$ is
\[\frac{n-k+1}{n+1}\binom{n+k}{n}.\]
This count is due to Whitworth~\cite[Proposition~II, p.~109]{Whi79}; see also Krattenthaler~\cite[Corollary~10.3.2]{Kra15}. We now
define a $q$-analogue of this binomial expression.

\begin{definition}\label{def:q-Yam}
The \emph{$q$-Yamanouchi number} $Y_q(n,k)$ is defined as
\[
Y_q(n,k) \coloneq \frac{[n-k+1]_q}{[n+1]_q}\begin{bmatrix} n+k \\ n \end{bmatrix}_q.
\]
\end{definition}

The generating function of the major index on $\mathcal{Y}(n,k)$ is
not equal to $Y_q(n,k)$ in general, up to a power of $q$. Its
coefficients need not be palindromic, whereas those of $Y_q(n,k)$ are.

\begin{definition}\label{def:mod_maj}
  For a finite binary word \( w \), append an additional letter $0$ at
  the end of $w$. Then the \emph{modified major index}
  $\widetilde{\operatorname{maj}}(w)$ of $w$ is defined as the major
  index of this new word with the dummy $0$. In other words, whenever
  $w$ ends with a letter $1$, we add its length to its major index;
  otherwise the major index remains unchanged.
\end{definition}

\begin{example}
The modified major index statistic on $\mathcal{Y}(5,3)$ is as follows.
\[
\begin{aligned}
\overset{\textcolor{blue}{8}}{\texttt{00000111}}\textcolor{orange}{\texttt{0}},\quad 
&\overset{\textcolor{blue}{8+5}}{\texttt{00001011}}\textcolor{orange}{\texttt{0}},\quad 
\overset{\textcolor{blue}{8+4}}{\texttt{00010011}}\textcolor{orange}{\texttt{0}},\quad 
\overset{\textcolor{blue}{8+3}}{\texttt{00100011}}\textcolor{orange}{\texttt{0}},\quad 
\overset{\textcolor{blue}{8+2}}{\texttt{01000011}}\textcolor{orange}{\texttt{0}},\quad 
\overset{\textcolor{blue}{8+6}}{\texttt{00001101}}\textcolor{orange}{\texttt{0}},\quad 
\overset{\textcolor{blue}{8+6+4}}{\texttt{00010101}}\textcolor{orange}{\texttt{0}}, \\[1ex]
\overset{\textcolor{blue}{8+6+3}}{\texttt{00100101}}\textcolor{orange}{\texttt{0}},\quad 
&\overset{\textcolor{blue}{8+6+2}}{\texttt{01000101}}\textcolor{orange}{\texttt{0}},\quad 
\overset{\textcolor{blue}{8+5}}{\texttt{00011001}}\textcolor{orange}{\texttt{0}},\quad 
\overset{\textcolor{blue}{8+5+3}}{\texttt{00101001}}\textcolor{orange}{\texttt{0}},\quad 
\overset{\textcolor{blue}{8+5+2}}{\texttt{01001001}}\textcolor{orange}{\texttt{0}},\quad 
\overset{\textcolor{blue}{8+4}}{\texttt{00110001}}\textcolor{orange}{\texttt{0}},\quad 
\overset{\textcolor{blue}{8+4+2}}{\texttt{01010001}}\textcolor{orange}{\texttt{0}}, \\[1ex]
\overset{\textcolor{blue}{7}}{\texttt{00001110}}\textcolor{orange}{\texttt{0}},\quad 
&\overset{\textcolor{blue}{7+4}}{\texttt{00010110}}\textcolor{orange}{\texttt{0}},\quad 
\overset{\textcolor{blue}{7+3}}{\texttt{00100110}}\textcolor{orange}{\texttt{0}},\quad 
\overset{\textcolor{blue}{7+2}}{\texttt{01000110}}\textcolor{orange}{\texttt{0}},\quad 
\overset{\textcolor{blue}{7+5}}{\texttt{00011010}}\textcolor{orange}{\texttt{0}},\quad 
\overset{\textcolor{blue}{7+5+3}}{\texttt{00101010}}\textcolor{orange}{\texttt{0}},\quad 
\overset{\textcolor{blue}{7+5+2}}{\texttt{01001010}}\textcolor{orange}{\texttt{0}}, \\[1ex]
\overset{\textcolor{blue}{7+4}}{\texttt{00110010}}\textcolor{orange}{\texttt{0}},\quad 
&\overset{\textcolor{blue}{7+4+2}}{\texttt{01010010}}\textcolor{orange}{\texttt{0}},\quad 
\overset{\textcolor{blue}{6}}{\texttt{00011100}}\textcolor{orange}{\texttt{0}},\quad 
\overset{\textcolor{blue}{6+3}}{\texttt{00101100}}\textcolor{orange}{\texttt{0}},\quad 
\overset{\textcolor{blue}{6+2}}{\texttt{01001100}}\textcolor{orange}{\texttt{0}},\quad 
\overset{\textcolor{blue}{6+4}}{\texttt{00110100}}\textcolor{orange}{\texttt{0}},\quad 
\overset{\textcolor{blue}{6+4+2}}{\texttt{01010100}}\textcolor{orange}{\texttt{0}}.
\end{aligned}
\]
The generating function is the \(q\)-analogue \(Y_q(n,k)\), multiplied
with the factor \(q^6\):
\[
\sum_{w\in\mathcal{Y}(5,3)}q^{\widetilde{\operatorname{maj}}(w)} =  q^6\frac{[3]_q}{[6]_q}\begin{bmatrix} 8 \\ 3 \end{bmatrix}_q.
\]
\end{example}

The following theorem will be used in
Subsection~\ref{subsection:decomp1}. This result is essentially a
special case of \cite[Theorem~5]{Kra89b}. We use a slightly different
modification of the major index.

\begin{theorem}\label{thm:Yam}
We have
\begin{align}\label{eq:q-Yam}
\sum_{w\in\mathcal{Y}(n,k)}q^{\widetilde{\operatorname{maj}}(w)} = q^{2k}\,Y_q(n,k).
\end{align}
\end{theorem}

In Subsection~\ref{subsection:decomp2}, we make use of the following
statement.

\begin{proposition}\label{prop:rev_maj}
\begin{enumerate}
\item Let $w$ be a binary word with $k$ many $1$s. Then $\mathrm{comaj}(w)=\widetilde{\mathrm{maj}}(w)-k$.
\item We have
\begin{align*}
\sum_{w\in\mathcal{Y}(n,k)}q^{\operatorname{comaj}(w)} = q^{k}\,Y_q(n,k).
\end{align*}
\item Let $\overline{w}$ denote the reverse word of $w$, that is, $\overline{w}\coloneq w_\ell\cdots w_1$ if $w=w_1\cdots w_\ell$. Then we have
\begin{align*}
\sum_{w\in\mathcal{Y}(n,k)}q^{\operatorname{maj}(\overline{w})} = 
\sum_{w\in\mathcal{Y}(n,k)}q^{\operatorname{comaj}(w)}.
\end{align*}
\end{enumerate}
\end{proposition}
\begin{proof}
  Each maximal string of $1$s in $w$, starting at position $s$ and of
  length $l$, contributes $s-1$ to $\operatorname{comaj}(w)$ and
  $s+l-1$ to $\widetilde{\operatorname{maj}}(w)$; the dummy letter $0$
  ensures the latter also for a string at the end of $w$. Summing the
  difference $l$ over all strings proves part~(i).

  Part~(ii) follows directly from Theorem~\ref{thm:Yam} and part~(i)
  of this proposition.

  For part~(iii), let $w\in\mathcal{Y}(n,k)$. The $n-k$ unmatched $0$s
  of $w$ give a unique decomposition
  \[
    w = D_0\,0\,D_1\,0\cdots0\,D_{n-k},
  \]
  where each $D_j$ is a \emph{Dyck word}, that is, a Yamanouchi word
  with equally many $0$s and $1$s, possibly empty. Let
  $\mathrm{rc}(D_j)$ denote the word obtained from $D_j$ by reversing
  it and interchanging $0$ and $1$. It is again a Dyck word. Define
  \[
    E(w) =
    \mathrm{rc}(D_{n-k})\,0\,\mathrm{rc}(D_{n-k-1})\,0\cdots0\,\mathrm{rc}(D_0).
  \]
  Then $E(w)\in\mathcal{Y}(n,k)$. Applying $E$ twice recovers $w$, so
  $E$ is a bijection from $\mathcal{Y}(n,k)$ to $\mathcal{Y}(n,k)$. We
  have
  $\operatorname{Asc}(E(w))=\{n+k-i : i\in\operatorname{Asc}(w)\}$.
  Thus
  \[
    \operatorname{comaj}(E(w)) = \sum_{i\in\operatorname{Asc}(w)}(n+k-i) = \operatorname{maj}(\overline{w}).
  \]
  Summing over all $w\in\mathcal{Y}(n,k)$ gives part~(iii).
\end{proof}

\subsection{Decomposition formula for $\mathcal{P}_{1,b}$}\label{subsection:decomp1}

In this subsection, we restrict \(\mathcal{B}_{1,b}^{\mathcal{P}}\) to its blue arrows. For each type, we compute the generating function of the starting points with that type and the length of their strings. Summing over all types gives the size generating function of \(\mathcal{P}_{1,b}\), which is \((-q;q)_\infty\), and proves the first identity in Theorem~\ref{thm:decomp1}. Applying the same procedure to the green arrows gives the second identity.

Recall that a starting point is a vertex with no incoming blue arrow. By Definition~\ref{def:level1_type_crystal}, \( \lambda \) is a starting point if and only if \( \type(\lambda) \) has no unmatched \( 1 \), that is, if and only if \( \type(\lambda) \) is a Yamanouchi word.

As long as the sequence of parts remains non-increasing, and we do not change the number of columns, adding or removing an even number of cells from a column does not change the type. In particular, there are infinitely many starting points with the same type.

\begin{definition}\label{def:pivot_P1}
  For a Yamanouchi word \( w \), we define the \emph{pivot} to be the
  smallest $2$-restricted partition \( \lambda \) with \( \type(\lambda) = w \), with
  respect to inclusion of diagrams.
\end{definition}

\begin{lemma}\label{lemma:piv_evenYam}
Pivots are in bijection with Yamanouchi words of even length.
\end{lemma}
\begin{proof}
If $\lambda$ is a $2$-restricted partition with $2k-1$ columns, then the first row of $\lambda$ ends in a green cell. Thus, $\type(\lambda)=r0$, where $r$ is the binary word of length $2k-1$ obtained by reading the bottom cells in the columns of $\lambda$ from left to right. If $\lambda$ has $2k$ columns, then the first row of $\lambda$ ends in a blue cell. Thus, $\type(\lambda)=r$, where $r$ is the binary word of length $2k$ obtained by reading the bottom cells in the columns of $\lambda$ from left to right. Therefore, $\type(\lambda)$ has even length for any partition $\lambda$.

If $\type(\lambda)$ has length $2k$, then $\lambda$ has $2k-1$ or $2k$ columns. Moreover, $\type(\lambda)$ determines the parity of each column height. The $2$-restricted condition requires the column heights to be strictly decreasing. Choosing the smallest allowed heights from right to left therefore gives the unique minimal $2$-restricted partition with the prescribed type.

This correspondence is illustrated below:
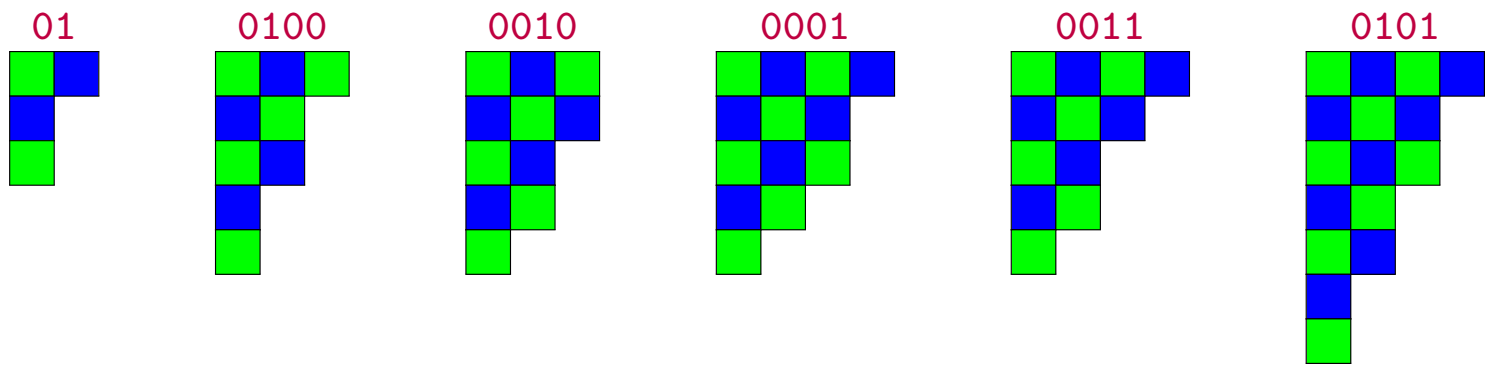
\begin{align*}
\overset{\textcolor{purple}{\texttt{01}}}{\scalebox{0.5}{\begin{ytableau}
  *(green)&*(blue)\\
  *(blue)\\
  *(green)
\end{ytableau}}}
\qquad
\overset{\textcolor{purple}{\texttt{0100}}}{\scalebox{0.5}{\begin{ytableau}
  *(green)&*(blue)&*(green)\\
  *(blue)&*(green)\\
  *(green)&*(blue)\\
  *(blue)\\
  *(green)
\end{ytableau}}}
\qquad
\overset{\textcolor{purple}{\texttt{0010}}}{\scalebox{0.5}{\begin{ytableau}
  *(green)&*(blue)&*(green)\\
  *(blue)&*(green)&*(blue)\\
  *(green)&*(blue)\\
  *(blue)&*(green)\\
  *(green)
\end{ytableau}}}
\qquad
\overset{\textcolor{purple}{\texttt{0001}}}{\scalebox{0.5}{\begin{ytableau}
  *(green)&*(blue)&*(green)&*(blue)\\
  *(blue)&*(green)&*(blue)\\
  *(green)&*(blue)&*(green)\\
  *(blue)&*(green)\\
  *(green)
\end{ytableau}}}
\qquad
\overset{\textcolor{purple}{\texttt{0011}}}{\scalebox{0.5}{\begin{ytableau}
  *(green)&*(blue)&*(green)&*(blue)\\
  *(blue)&*(green)&*(blue)\\
  *(green)&*(blue)\\
  *(blue)&*(green)\\
  *(green)
\end{ytableau}}}
\qquad
\overset{\textcolor{purple}{\texttt{0101}}}{\scalebox{0.5}{\begin{ytableau}
  *(green)&*(blue)&*(green)&*(blue)\\
  *(blue)&*(green)&*(blue)\\
  *(green)&*(blue)&*(green)\\
  *(blue)&*(green)\\
  *(green)&*(blue)\\
  *(blue)\\
  *(green)
\end{ytableau}}} 
\end{align*} \qedhere
\end{proof}

The following lemma shows that the size of a pivot can be expressed
in terms of the binary word \( \type(\lambda) \).

\begin{lemma}\label{lem:wt_pivot}
  Let \( \lambda \) be a pivot with
  \( w = \type(\lambda) \in \mathcal{Y}(n,k) \), where \( n+k \) is
  even and \( n \geq k \geq 0 \). Then we have
  \[
    |\lambda| = \binom{n+k}{2} + 2\, \widetilde{\operatorname{maj}}(w) -k.
  \]
\end{lemma}

\begin{proof}
  We construct \( \lambda \) from the staircase
  \( (n+k-1, n+k-2, \ldots, 1) \) by adding diagonals, as illustrated
  in Figure~\ref{fig:pivot_from_Yam}. Adding a diagonal of length
  \( d \) means adding one cell at the bottom of each of the first
  \( d \) columns. Let \( \Psi(w) \) be the partition whose parts are
  the lengths of the added diagonals.

  The largest part of \( \Psi(w) \), i.e. the number of cells along
  the first (blue) diagonal, is the largest descent of the word $w0$,
  i.e. the word obtained from $w$ by appending a letter $0$ at the
  end. The second part of $\Psi(w)$, i.e. the number of cells along
  the second (green) diagonal, is the largest ascent of $w$. In
  general, the $(2i-1)$th part of $\Psi(w)$ is the $i$th largest
  descent of $w0$, and the $(2i)$th part is the $i$th largest ascent
  of $w$. Then \( \lambda \) is the partition obtained by adding these
  diagonals to the staircase.

  For example, the Yamanouchi word $w = \texttt{00010011}$ corresponds
  to the pivot illustrated in Figure~\ref{fig:pivot_from_Yam}. The
  associated partition of diagonals is $\Psi(w) = (8,6,4,3)$.

\begin{figure}[h]
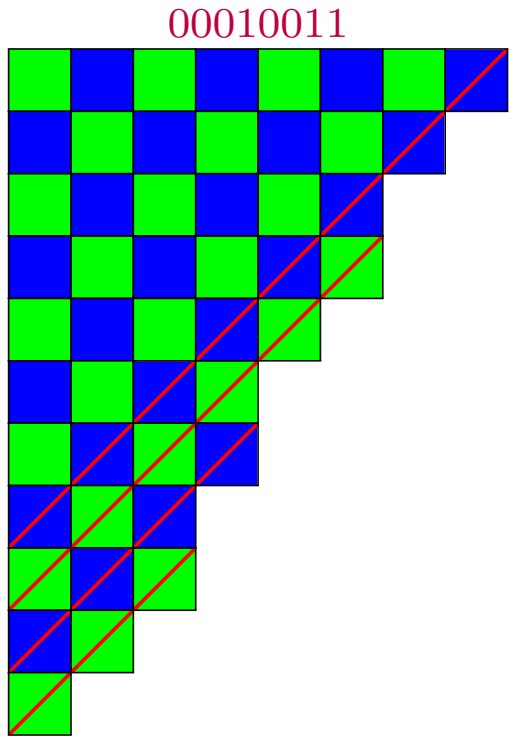

    \centering
\begin{align*}
\overset{\textcolor{purple}{00010011}}{
\scalebox{0.7}{$
\begin{ytableau}
  *(green) & *(blue)  & *(green) & *(blue)  & *(green)  & *(blue)  & *(green) & \DiagCellb \\
  *(blue)  & *(green) & *(blue)  & *(green) & *(blue)   & *(green) & \DiagCellb  \\
  *(green) & *(blue)  & *(green) & *(blue)  & *(green)  & \DiagCellb \\
  *(blue)  & *(green) & *(blue)  & *(green) & \DiagCellb & \DiagCellg \\
  *(green) & *(blue)  & *(green) & \DiagCellb & \DiagCellg \\
  *(blue)  & *(green) & \DiagCellb  & \DiagCellg \\
  *(green) & \DiagCellb & \DiagCellg & \DiagCellb \\
  \DiagCellb & \DiagCellg & \DiagCellb \\
  \DiagCellg  & \DiagCellb & \DiagCellg \\
  \DiagCellb & \DiagCellg \\
  \DiagCellg
\end{ytableau}
$}
}
\end{align*}
    \caption{Obtaining the pivot from the staircase}
    \label{fig:pivot_from_Yam}
\end{figure}

To see that any partition constructed in this way is a pivot, note that in the staircase consecutive parts differ by at most one, and we add at most one more cell to column $i$ than to column $i+1$ along the diagonals.

By definition, the contribution of the blue parts is
$\widetilde{\operatorname{maj}}(w)$ and the contribution of the green parts is $\operatorname{comaj}(w)$.
Therefore, by Proposition~\ref{prop:rev_maj}~(i), the size $|\Psi(w)|$ is given by
$2\,\widetilde{\operatorname{maj}}(w)-k$.
\end{proof}

Let \( \mu \) be a starting point with
\( \type(\mu) = w \in \mathcal{Y}(n,k) \), where \( n+k \) is even,
and \( n \geq k \geq 0 \). Then there exists a unique pivot \( \lambda \)
with the same type \( w \), and \( \mu \) can be obtained from
\( \lambda \) in a unique way by adding even numbers \( b_1 \geq \dots \geq b_{n+k} \)
to the columns of its Young diagram, beginning with the first column.
Therefore, we have the size generating function for starting points
for a fixed pivot.

\begin{lemma}\label{lem:gf_fixed_pivot}
  The size generating function of starting points corresponding to a
  pivot \( \lambda \) with \( \type(\lambda) \in \mathcal{Y}(n,k) \) is
  \[
    q^{|\lambda|}\,(q^2;q^2)_{n+k}^{-1}.
  \]
\end{lemma}

We are now ready to show the following identity, which relates the
infinite product \( (-q;q)_\infty \) to \( q \)-Yamanouchi numbers.
This identity is also stated in Theorem~\ref{thm:decomp1}.

\begin{theorem}\label{thm:decomp1_blue}
We have
\begin{align*}
(-q;q)_\infty = \sum_{\substack{n\geq k\geq 0}} q^{n(2n-1) + 3k}\, Y_{q^2}(2n-k,k)\,(q^2;q^2)_{2n}^{-1}\,[2(n-k)+1]_q.
\end{align*}
\end{theorem}
\begin{proof}
  Let \( \lambda \) be a pivot with
  \( w = \type(\lambda) \in \mathcal{Y}(n,k) \), where \( n+k \) is
  even and \( n \geq k \geq 0 \). By Lemma~\ref{lem:wt_pivot} and
  Lemma~\ref{lem:gf_fixed_pivot}, the generating function of starting
  points whose pivot is \( \lambda \) is
  \[
    q^{\binom{n+k}{2} + 2\, \widetilde{\operatorname{maj}}(w) -k} (q^2; q^2)_{n+k}^{-1}.
  \]
  All strings whose starting point has pivot $\lambda$ have length
  $n-k$, which is the number of unmatched $0$s in $\type(\lambda)$.
  Each arrow increases the size by one. Hence we multiply the
  generating function by $[n-k+1]_q$.

  To simplify the condition for \( n \) and \( k \), replace \( n \)
  by \( 2n-k \). The infinite summation formula, obtained by summing
  the product of the size generating function of starting points and
  the length of its blue string, is equal to
  \[
    \sum_{ n \geq k \geq 0 } \sum_{ w \in \mathcal{Y}(2n-k,k) } q^{\binom{2n}{2} + 2\, \widetilde{\operatorname{maj}}(w) - k} (q^2;q^2)_{2n}^{-1}\,[2(n-k)+1]_q.
  \]
  By Theorem~\ref{thm:Yam}, we obtain the desired result.
\end{proof}

We could also consider the restriction formula obtained by erasing the
blue arrows in the affine crystal graph $\mathcal{B}_{1,b}^{\mathcal{P}}$
(see Figure~\ref{fig:checkerboard}), instead of the green
arrows. We obtain the second formula in Theorem~\ref{thm:decomp1}.
\begin{theorem}\label{thm:decomp1_green}
We have
\begin{align*}
(-q;q)_\infty = \sum_{\substack{n\geq k\geq 0}} q^{(2n+1)n + 3k}\, Y_{q^2}(2n+1-k,k)\,(q^2;q^2)_{2n+1}^{-1}\,[2(n-k)+2]_q.
\end{align*}
\end{theorem}
\begin{proof}
The same argument goes through as in the proof of Theorem~\ref{thm:decomp1_blue}, with the only difference that we start building from an even staircase,
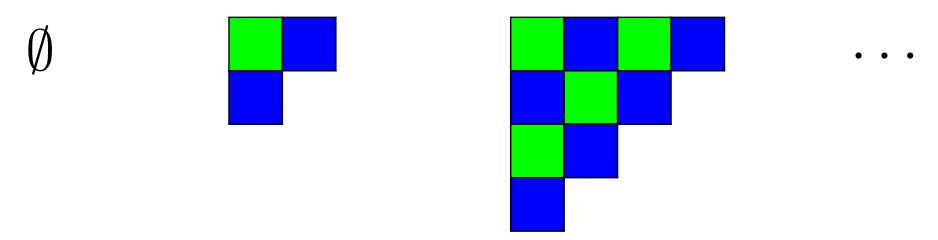
\[
\emptyset
\hspace{3em}
\scalebox{0.6}{\begin{ytableau}
*(green)&*(blue)\\
*(blue)
\end{ytableau}}\hspace{3em}
\scalebox{0.6}{\begin{ytableau}
  *(green)&*(blue)&*(green)&*(blue)\\
  *(blue)&*(green)&*(blue)\\
  *(green)&*(blue)\\
  *(blue)
\end{ytableau}}\hspace{2em}
\cdots\]
and the length of each green string is odd.
\end{proof}

\begin{remark}
Note that this is equivalent to restricting the affine crystal graph of highest weight $\Lambda_1$ and type $A_1^{(1)}$ to blue arrows. The partition model for this graph is the same as Definition~\ref{def:level1_type_crystal}, but each partition has a blue cell in the top-left corner instead of a green cell.
\end{remark}

\subsection{Decomposition formula for $\mathcal{P}_{2,b}$}\label{subsection:decomp2}

In this subsection, we consider the restriction of the affine crystal graph
\(\mathcal{B}_{2,b}^{\mathcal{P}}\) to the blue arrows. Once again,
the connected infinite crystal graph
\(\mathcal{B}_{2,b}^{\mathcal{P}}\), illustrated in
Figure~\ref{fig:level2_grounded_crystal}, breaks into a disjoint union
of infinitely many blue strings of finite length. Based on this
decomposition, we express the infinite product
\((-q;q^2)_\infty / (q;q^2)_\infty\) as an infinite sum of
\(q\)-binomial expressions.

In Subsection~\ref{subsection:decomp1}, we described starting points
in terms of the type. Then, we computed the generating function for pivots, the smallest partitions with given type. Similarly, in the
graph \(\mathcal{B}_{2,b}^{\mathcal{P}}\), a vertex \( \pi \) is a
starting point if and only if \( \type(\pi) \) is a Yamanouchi word.

Here, it is convenient to provide an alternative description that uses the sequence of labels of
the parts of \( \pi \) instead of the type. In this approach, we
define the pivot in a different way, again as the smallest partition
in an appropriate sense, and derive the decomposition formula. Note
that the type of \( \pi \) is obtained from the
bottommost cells in each column of its Young diagram, whereas the sequence of labels of the parts is obtained from the rightmost cells in each row.

\begin{definition}\label{def:pivot2}
  A \emph{starting point} in \(\mathcal{B}_{2,b}^{\mathcal{P}}\) is a
  grounded partition in \( \mathcal{P}_{2,b} \) that is not the target
  of a blue arrow. A starting point in \(\mathcal{B}_{2,b}^{\mathcal{P}}\) is called a \emph{pivot} if it has
  no \(b\)-parts. We denote by \(\mathrm{PIV}(n,k)\) the set of
  pivots with \(n\) many \(a\)-parts and \(k\) many \(c\)-parts.
\end{definition}

\begin{lemma}\label{lemma:level2decomp_step1}
  Let \(\pi\in\mathcal{P}_{2,b}\) be a grounded partition with \(n\)
  many \(a\)-parts, \(k\) many \(c\)-parts and no \(b\)-parts. Then
\begin{enumerate}
\item \(\pi\) is a starting point if and only if its sequence of
  labels, read from the smallest nonzero part upward, forms a
  Yamanouchi word, where \( a \) is interpreted as \( 0 \), and
  \( c \) as \( 1 \).
\item If \(\pi\) is a starting point with corresponding Yamanouchi
  word \(w\in\mathcal{Y}(n,k)\), then the length (i.e. the number of
  blue arrows) of the blue string starting from \(\pi\) is
  \(2(n - k) + 1\).
\end{enumerate}
\end{lemma}
\begin{proof}
  For a starting point \(\pi\), the word \(\type(\pi)\) starts with
  \(0\). Hence, the first part of \(\pi\) is \(1_a\). An \(a\)-part
  that follows an \(a\)-part contributes \(00\) to
  \(\type(\pi)\). An \(a\)-part that follows a \(c\)-part lies in the
  same column, and therefore it does not contribute to \(\type(\pi)\).
  Similarly, a \(c\)-part that follows an \(a\)-part does not
  contribute to \(\type(\pi)\). A \(c\)-part that follows a \(c\)-part
  contributes \(11\) to \(\type(\pi)\).

  Let \( w \) be the binary word obtained from the sequence of labels of
  the parts of \( \pi \). The above argument shows that
  \( \type(\pi) \) is a Yamanouchi word if and only if \( w \) is a
  Yamanouchi word. This proves the first statement.

  Moreover, after the first letter \(0\), we have \(2(n - k)\) many
  unmatched \(0\)s, so the length of the blue string starting from
  \(\pi\) is \(2(n - k)+1\).
\end{proof}

Next, we consider arbitrary starting points, not only pivots.

\begin{lemma}\label{lemma:level2decomp_step2}
  A grounded partition \( \pi \in\mathcal{P}_{2,b}\) is a starting
  point if and only if the subsequence of labels, starting from the
  smallest nonzero part and restricted to the letters \( a \) and
  \( c \), forms a Yamanouchi word, with \( a \) interpreted as
  \( 0 \) and \( c \) as \( 1 \). Moreover, the number of
  \( b \)-parts does not affect the length of the corresponding blue
  string.
\end{lemma}
\begin{proof}
  A \emph{string of \(b\)-parts} in \(\pi\) is a maximal sequence of
  consecutive \(b\)-parts \(k_b k_b\cdots k_b\), that is, all the
  \(b\)-parts of the same size. We proceed by induction on the number
  of strings of \(b\)-parts. The base case is
  Lemma~\ref{lemma:level2decomp_step1}, where \(\pi\) has no
  \(b\)-parts. Assume the statement holds for all \(\pi\) with at most
  \(i\) strings of \(b\)-parts. Let \(\pi\) have \(i+1\) strings of \(b\)-parts, and let \(k_b k_b \cdots k_b\) denote the
  \((i+1)\)th string. Let \(\tilde{\pi} \in \mathcal{P}_{2,b}\) be
  the grounded partition that has the same sequence of labels as \(\pi\),
  except that the \((i+1)\)th string of \(b\)’s is removed from the
  sequence of labels (note that there is a unique grounded partition in
  \(\mathcal{P}_{2,b}\) with a given sequence of labels). For example,
  for \(\pi=1_a1_c1_a2_b2_b3_a\) we have \(\tilde{\pi}=1_a1_c1_a3_a\);
  for \(\pi=1_a1_c2_b3_a\) we have \(\tilde{\pi}=1_a1_c1_a\).

  First, assume that \(\pi\) also contains a part of size \(k+1\),
  i.e. \(k_b\) is not the largest part. We then distinguish four
  cases, according to the labels of the odd parts that precede and
  follow the \((i+1)\)th string of \(b\)-parts.

\begin{enumerate}[label=(1\alph*)]
\item Assume the $(i+1)$th string of \( b \)-parts in $\pi$ lies between two \( a \)-parts, $(k-1)_a$ and $(k+1)_a$. This is illustrated in Figure~\ref{fig:sec4_1abcd}~(1a). Since the difference between two consecutive \(a\)-parts is \(2\), the number of columns in \(\pi\) and \(\tilde{\pi}\) is the same.  
Because \(a\) and \(c\) alternate in the odd columns and the part preceding the first $b$-part is an $a$-part, the bottom cell in the $k$th column is a green \(b\)-cell in both cases. Therefore, \(\type(\pi) = \type(\tilde{\pi})\). Hence $\pi$ is a starting point if and only if $\tilde{\pi}$ is a starting point, and in this case, the length of the corresponding strings is the same.

\begin{figure}[h]
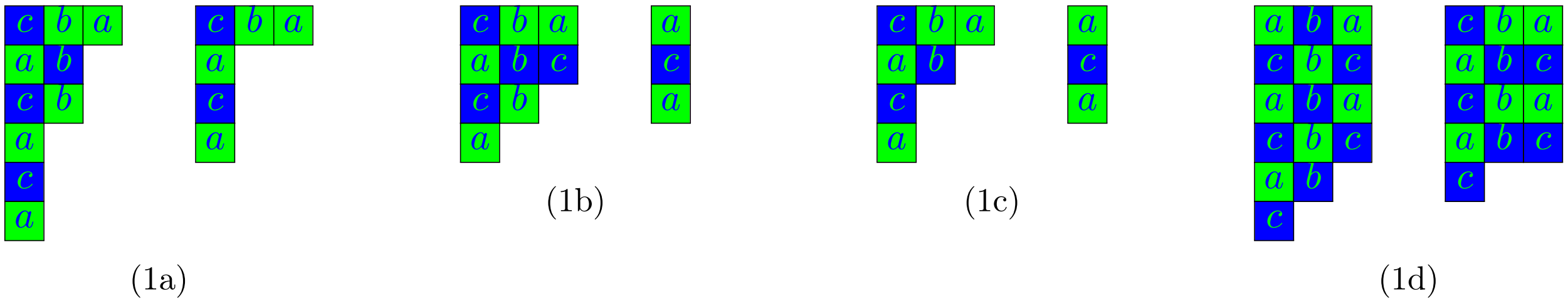

    \centering

    \begin{minipage}[t]{0.24\textwidth}
        \centering
        \(\scalebox{0.7}{\begin{ytableau}
*(blue) \scalebox{1.5}{\color{green}$c$} & *(green) \scalebox{1.5}{\color{blue}$b$} & *(green) \scalebox{1.5}{\color{blue}$a$} \\
*(green) \scalebox{1.5}{\color{blue}$a$} & *(blue) \scalebox{1.5}{\color{green}$b$} \\
*(blue) \scalebox{1.5}{\color{green}$c$} & *(green) \scalebox{1.5}{\color{blue}$b$} \\
*(green) \scalebox{1.5}{\color{blue}$a$} \\
*(blue) \scalebox{1.5}{\color{green}$c$} \\
*(green) \scalebox{1.5}{\color{blue}$a$}
\end{ytableau}} \hspace{2em} \scalebox{0.7}{\begin{ytableau}
*(blue) \scalebox{1.5}{\color{green}$c$} & *(green) \scalebox{1.5}{\color{blue}$b$} & *(green) \scalebox{1.5}{\color{blue}$a$} \\
*(green) \scalebox{1.5}{\color{blue}$a$} \\
*(blue) \scalebox{1.5}{\color{green}$c$} \\
*(green) \scalebox{1.5}{\color{blue}$a$}
\end{ytableau}}\)

        \medskip
        {\small (1a)}
    \end{minipage}
    \begin{minipage}[t]{0.24\textwidth}
        \centering
\(\scalebox{0.7}{\begin{ytableau}
*(blue) \scalebox{1.5}{\color{green}$c$} & *(green) \scalebox{1.5}{\color{blue}$b$} & *(green) \scalebox{1.5}{\color{blue}$a$} \\
*(green) \scalebox{1.5}{\color{blue}$a$} & *(blue) \scalebox{1.5}{\color{green}$b$} & *(blue) \scalebox{1.5}{\color{green}$c$} \\
*(blue) \scalebox{1.5}{\color{green}$c$} & *(green) \scalebox{1.5}{\color{blue}$b$} \\
*(green) \scalebox{1.5}{\color{blue}$a$}
\end{ytableau}} \hspace{2em} \scalebox{0.7}{\begin{ytableau}
*(green) \scalebox{1.5}{\color{blue}$a$}  \\
*(blue) \scalebox{1.5}{\color{green}$c$} \\
*(green) \scalebox{1.5}{\color{blue}$a$}
\end{ytableau}}\)

        \medskip
        {\small (1b)}
    \end{minipage}
    \begin{minipage}[t]{0.24\textwidth}
        \centering
\(\scalebox{0.7}{\begin{ytableau}
*(blue) \scalebox{1.5}{\color{green}$c$} & *(green) \scalebox{1.5}{\color{blue}$b$} & *(green) \scalebox{1.5}{\color{blue}$a$} \\
*(green) \scalebox{1.5}{\color{blue}$a$} & *(blue) \scalebox{1.5}{\color{green}$b$} \\
*(blue) \scalebox{1.5}{\color{green}$c$} \\
*(green) \scalebox{1.5}{\color{blue}$a$}
\end{ytableau}} \hspace{2em} \scalebox{0.7}{\begin{ytableau}
*(green) \scalebox{1.5}{\color{blue}$a$}  \\
*(blue) \scalebox{1.5}{\color{green}$c$} \\
*(green) \scalebox{1.5}{\color{blue}$a$}
\end{ytableau}}\)

        \medskip
        {\small (1c)}
    \end{minipage}
    \begin{minipage}[t]{0.24\textwidth}
        \centering
\(\scalebox{0.7}{\begin{ytableau}
*(green) \scalebox{1.5}{\color{blue}$a$} & *(blue) \scalebox{1.5}{\color{green}$b$} & *(green) \scalebox{1.5}{\color{blue}$a$} \\
*(blue) \scalebox{1.5}{\color{green}$c$} & *(green) \scalebox{1.5}{\color{blue}$b$} & *(blue) \scalebox{1.5}{\color{green}$c$} \\
*(green) \scalebox{1.5}{\color{blue}$a$} & *(blue) \scalebox{1.5}{\color{green}$b$} & *(green) \scalebox{1.5}{\color{blue}$a$} \\
*(blue) \scalebox{1.5}{\color{green}$c$} & *(green) \scalebox{1.5}{\color{blue}$b$} & *(blue) \scalebox{1.5}{\color{green}$c$} \\
*(green) \scalebox{1.5}{\color{blue}$a$} & *(blue) \scalebox{1.5}{\color{green}$b$} \\
*(blue) \scalebox{1.5}{\color{green}$c$}
\end{ytableau}} \hspace{2em} \scalebox{0.7}{\begin{ytableau}
*(blue) \scalebox{1.5}{\color{green}$c$} & *(green) \scalebox{1.5}{\color{blue}$b$} & *(green) \scalebox{1.5}{\color{blue}$a$} \\
*(green) \scalebox{1.5}{\color{blue}$a$} & *(blue) \scalebox{1.5}{\color{green}$b$} & *(blue) \scalebox{1.5}{\color{green}$c$} \\
*(blue) \scalebox{1.5}{\color{green}$c$} & *(green) \scalebox{1.5}{\color{blue}$b$} & *(green) \scalebox{1.5}{\color{blue}$a$} \\
*(green) \scalebox{1.5}{\color{blue}$a$} & *(blue) \scalebox{1.5}{\color{green}$b$} & *(blue) \scalebox{1.5}{\color{green}$c$} \\
*(blue) \scalebox{1.5}{\color{green}$c$}
\end{ytableau}}\)

        \medskip
        {\small (1d)}
    \end{minipage}

    \caption{Inductive steps for cases (1a), (1b), (1c), and (1d).}
    \label{fig:sec4_1abcd}
\end{figure}

\item Assume the $(i+1)$th string of \( b \)-parts is preceded by an $a$-part $(k-1)_a$ and followed by a $c$-part $(k+1)_c$. This case is illustrated in Figure~\ref{fig:sec4_1abcd}~(1b).

In $\pi$, the column of $b$s has a green cell at the bottom and the next column has a $c$-cell at the bottom, which contributes the factor $01$ to $\type(\pi)$. In $\tilde{\pi}$, these two columns disappear and all other column endings are unchanged. Therefore, \(\type(\pi)\sim\type(\tilde{\pi})\).

\item Assume the $(i+1)$th string of \( b \)-parts is preceded by a \( c \)-part $(k-1)_c$ and followed by an \( a \)-part $(k+1)_a$. This case is illustrated in Figure~\ref{fig:sec4_1abcd}~(1c). In $\pi$, we have two consecutive columns ending in $b$ and $a$, where the $b$-cell is blue. This contributes the factor $10$ to $\type(\pi)$. In $\tilde{\pi}$, these two columns disappear. Since the factor $10$ is preceded by an odd number $k-1$ of letters, \(\type(\pi)\sim\type(\tilde{\pi})\) whenever $\pi$ or $\tilde{\pi}$ is a starting point.

\item Assume the $(i+1)$th string of \( b \)-parts lies between two \( c \)-parts $(k-1)_c$ and $(k+1)_c$. Since the difference between two consecutive \(c\)-parts is \(2\), the number of columns in \(\pi\) and \(\tilde{\pi}\) is the same. Because \(a\) and \(c\) alternate in the odd columns, the bottom \(b\)-cell is blue in both cases. This is illustrated in Figure~\ref{fig:sec4_1abcd}~(1d). Therefore, \(\type(\pi) = \type(\tilde{\pi})\).

\end{enumerate}
Finally, assume that $\pi$ ends with the $(i+1)$th string of $b$s. We distinguish two cases, according to the label of the preceding odd part.
\begin{enumerate}[label=(2\alph*)]
\item Assume the last string of $b$s follows an $a$-part $(k-1)_a$. Then the last column of $b$s in $\pi$ has a green cell at the bottom, which contributes a letter $0$ to $\type(\pi)$. Since the largest part of $\pi$ is $k_b$, the first row of $\pi$ ends in a $b$-cell, which contributes another letter $0$ at the end of $\type(\pi)$. In $\tilde{\pi}$, the largest part is $(k-1)_a$, so the first row of $\tilde{\pi}$ ends in an $a$-cell, which contributes $00$ at the end of $\type(\tilde{\pi})$. Therefore, \(\type(\pi) = \type(\tilde{\pi})\). This case is illustrated in Figure~\ref{fig:sec4_2ab}~(2a).

\begin{figure}[h]
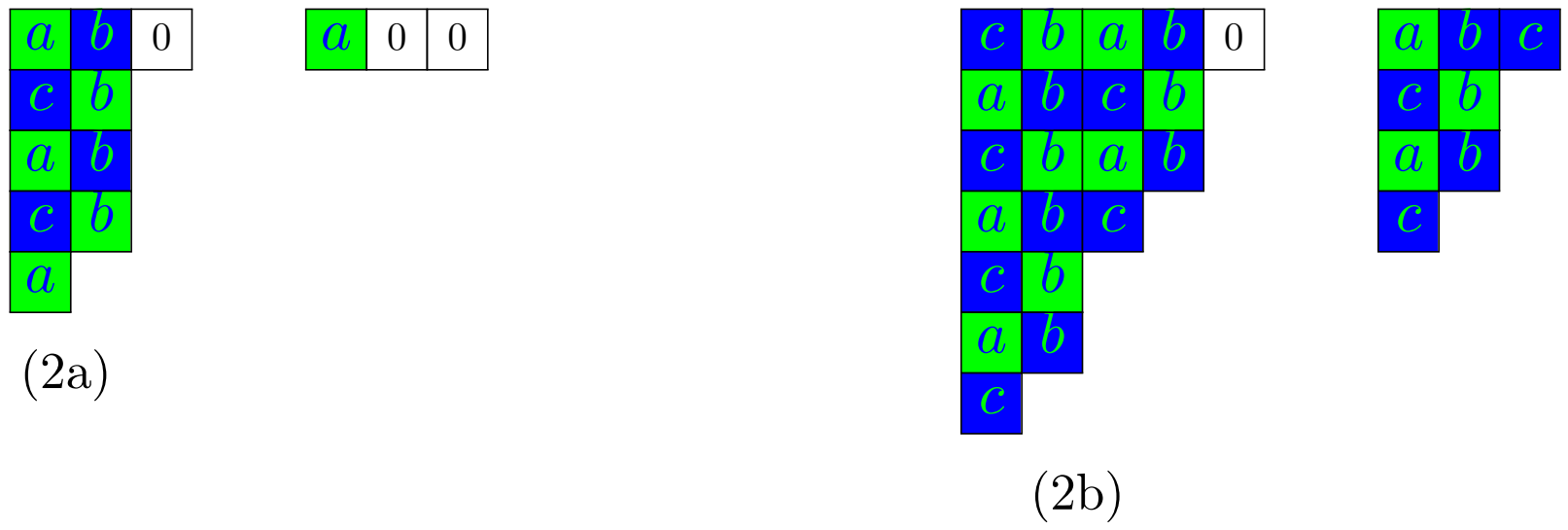

    \centering

    \begin{minipage}[t]{0.38\textwidth}
        \centering
        \hspace{0.35\textwidth}
        \(\scalebox{0.7}{\begin{ytableau}
*(green) \scalebox{1.5}{\color{blue}$a$} & *(blue) \scalebox{1.5}{\color{green}$b$} & 0 \\
*(blue) \scalebox{1.5}{\color{green}$c$} & *(green) \scalebox{1.5}{\color{blue}$b$}  \\
*(green) \scalebox{1.5}{\color{blue}$a$} & *(blue) \scalebox{1.5}{\color{green}$b$} \\
*(blue) \scalebox{1.5}{\color{green}$c$} & *(green) \scalebox{1.5}{\color{blue}$b$} \\
*(green) \scalebox{1.5}{\color{blue}$a$}
\end{ytableau}} \hspace{2em} \scalebox{0.7}{\begin{ytableau}
*(green) \scalebox{1.5}{\color{blue}$a$} & 0 & 0
\end{ytableau}}\)

        \medskip
        {\small (2a)}
    \end{minipage}
    \begin{minipage}[t]{0.38\textwidth}
        \centering
        \hspace{0.35\textwidth}
        \(\scalebox{0.7}{\begin{ytableau}
*(blue) \scalebox{1.5}{\color{green}$c$} & *(green) \scalebox{1.5}{\color{blue}$b$} & *(green) \scalebox{1.5}{\color{blue}$a$} & *(blue) \scalebox{1.5}{\color{green}$b$} & 0 \\
*(green) \scalebox{1.5}{\color{blue}$a$} & *(blue) \scalebox{1.5}{\color{green}$b$} & *(blue) \scalebox{1.5}{\color{green}$c$} & *(green) \scalebox{1.5}{\color{blue}$b$}  \\
*(blue) \scalebox{1.5}{\color{green}$c$} & *(green) \scalebox{1.5}{\color{blue}$b$} & *(green) \scalebox{1.5}{\color{blue}$a$} & *(blue) \scalebox{1.5}{\color{green}$b$} \\
*(green) \scalebox{1.5}{\color{blue}$a$} & *(blue) \scalebox{1.5}{\color{green}$b$} & *(blue) \scalebox{1.5}{\color{green}$c$} \\
*(blue) \scalebox{1.5}{\color{green}$c$} & *(green) \scalebox{1.5}{\color{blue}$b$} \\
*(green) \scalebox{1.5}{\color{blue}$a$} & *(blue) \scalebox{1.5}{\color{green}$b$} \\
*(blue) \scalebox{1.5}{\color{green}$c$}
\end{ytableau}} \hspace{2em} \scalebox{0.7}{\begin{ytableau}
*(green) \scalebox{1.5}{\color{blue}$a$} & *(blue) \scalebox{1.5}{\color{green}$b$} & *(blue) \scalebox{1.5}{\color{green}$c$} \\
*(blue) \scalebox{1.5}{\color{green}$c$} & *(green) \scalebox{1.5}{\color{blue}$b$} \\
*(green) \scalebox{1.5}{\color{blue}$a$} & *(blue) \scalebox{1.5}{\color{green}$b$} \\
*(blue) \scalebox{1.5}{\color{green}$c$}
\end{ytableau}}\)

        \medskip
        {\small (2b)}
    \end{minipage}

    \caption{Inductive steps for cases (2a) and (2b).}
    \label{fig:sec4_2ab}
\end{figure}

\item Assume the last string of \(b\)’s follows a \(c\)-part \((k-1)_c\), and that \(k_b\) is the largest part in \(\pi\).  
This case is illustrated in Figure~\ref{fig:sec4_2ab}.  
Then the last column of \(b\)’s has a blue cell at the bottom, and the first row of \(\pi\) contributes an extra letter \(0\) to \(\type(\pi)\).  
Therefore, \(\type(\tilde{\pi})\) is obtained from \(\type(\pi)\) by removing the right factor \(10\).

First, assume that \(\tilde{\pi}\) is a starting point. By the induction hypothesis and Lemma~\ref{lemma:level2decomp_step1}~(ii), the blue string starting from \(\tilde{\pi}\) has odd length, so \(\type(\tilde{\pi})\) has at least one unmatched letter \(0\). Therefore, the letter \(1\) in the right factor \(10\) of \(\type(\pi)\) gets matched. The letter \(0\) in the right factor \(10\) is unmatched. Hence, \(\type(\pi)\sim\type(\tilde{\pi})\); in particular, \(\pi\) is also a starting point.

Next, assume that \(\tilde{\pi}\) is not a starting point. Then \(\type(\tilde{\pi})\) has an unmatched letter \(1\), which remains unmatched in \(\type(\pi)\). Hence \(\pi\) is not a starting point either.\qedhere
\end{enumerate}
\end{proof}

By Lemmas~\ref{lemma:level2decomp_step1} and~\ref{lemma:level2decomp_step2}, a grounded partition is a starting point if and only if its subsequence of labels restricted to the letters \(a\) and \(c\) is the sequence of labels of a pivot. We say that the starting point \emph{corresponds} to this pivot.

We need to compute the size generating function of pivots in $\mathrm{PIV}(n,k)$.

\begin{example}
Consider the pivots with four $a$-parts and two $c$-parts:
\[
\begin{aligned}
&\quad \textcolor{black}{1_a}\,\textcolor{black}{3_a}\,\textcolor{black}{5_a}\,\textcolor{black}{7_a}\,\textcolor{black}{7_c}\,\textcolor{black}{9_c}, \ 
\textcolor{black}{1_a}\,\textcolor{black}{3_a}\,\textcolor{black}{5_a}\,\textcolor{black}{5_c}\,\textcolor{black}{5_a}\,\textcolor{black}{5_c}, \ 
\textcolor{black}{1_a}\,\textcolor{black}{3_a}\,\textcolor{black}{3_c}\,\textcolor{black}{3_a}\,\textcolor{black}{5_a}\,\textcolor{black}{5_c}, \ 
\textcolor{black}{1_a}\,\textcolor{black}{1_c}\,\textcolor{black}{1_a}\,\textcolor{black}{3_a}\,\textcolor{black}{5_a}\,\textcolor{black}{5_c}, \ 
\textcolor{black}{1_a}\,\textcolor{black}{3_a}\,\textcolor{black}{5_a}\,\textcolor{black}{5_c}\,\textcolor{black}{7_c}\,\textcolor{black}{7_a}, \\[1em]
&\quad \textcolor{black}{1_a}\,\textcolor{black}{3_a}\,\textcolor{black}{3_c}\,\textcolor{black}{3_a}\,\textcolor{black}{3_c}\,\textcolor{black}{3_a}, \ 
\textcolor{black}{1_a}\,\textcolor{black}{1_c}\,\textcolor{black}{1_a}\,\textcolor{black}{3_a}\,\textcolor{black}{3_c}\,\textcolor{black}{3_a}, \ 
\textcolor{black}{1_a}\,\textcolor{black}{3_a}\,\textcolor{black}{3_c}\,\textcolor{black}{5_c}\,\textcolor{black}{5_a}\,\textcolor{black}{7_a}, \ 
\textcolor{black}{1_a}\,\textcolor{black}{1_c}\,\textcolor{black}{1_a}\,\textcolor{black}{1_c}\,\textcolor{black}{1_a}\,\textcolor{black}{3_a}.
\end{aligned}
\]
Their size generating function is given by
\[q^8+q^{12}+2q^{16}+q^{20}+2q^{24}+q^{28}+q^{32} = q^8\,Y_{q^4}(4,2).\]
\end{example}

\begin{lemma}\label{lem:wt_pivot_2}
Let \( \pi \in \mathrm{PIV}(n,k)\) be a pivot whose parts are labelled only $a$ or $c$. Let $w\in\mathcal{Y}(n,k)$ denote the word obtained by replacing $a$ by $0$ and $c$ by $1$ in the sequence of labels of $\pi$. Then we have
  \[
    |\pi| = (n+k)^2-4\,\widetilde{\mathrm{maj}}(\overline{w})+2k,
  \]
where $\overline{w}$ denotes the reverse word of $w$.
\end{lemma}
\begin{proof}
We begin with the odd staircase $(1,3,\ldots,2(n+k)-1)$. The size of this staircase is $(n+k)^2$. Each descent \( i \) in $\overline{w}$ decreases the size of the odd staircase by \( 2i \), since, instead of jumping by two, we have two repeated odd parts. Similarly, each ascent \( i \) decreases the size by \( 2i \). Since $w$ is a Yamanouchi word, $\overline{w}$ does not end with $1$, so $\mathrm{maj}(\overline{w})=\widetilde{\mathrm{maj}}(\overline{w})$. The proof is completed by Proposition~\ref{prop:rev_maj}~(i).
\end{proof}

\begin{example}
Figure~\ref{fig:P2b_pivot_wt} illustrates the example $\pi = 1_a 3_a 3_c 3_a 5_a 7_a 7_c$ inside the odd staircase with $7$ parts. In this example, the sequence of labels is $aacaaac$, so the associated binary word is $w = 0010001$. The reverse word is $\overline{w} = 1000100$, which has descents at positions $1$ and $5$ and an ascent at position $4$. In Figure~\ref{fig:P2b_pivot_wt}, we first remove $5$ pairs of cells from the top $5$ rows of the odd staircase because of the descent at position $5$, then remove $4$ pairs of cells from the top $4$ rows because of the ascent at position $4$, and finally remove one pair of cells from the top row because of the descent at position $1$.
\begin{figure}
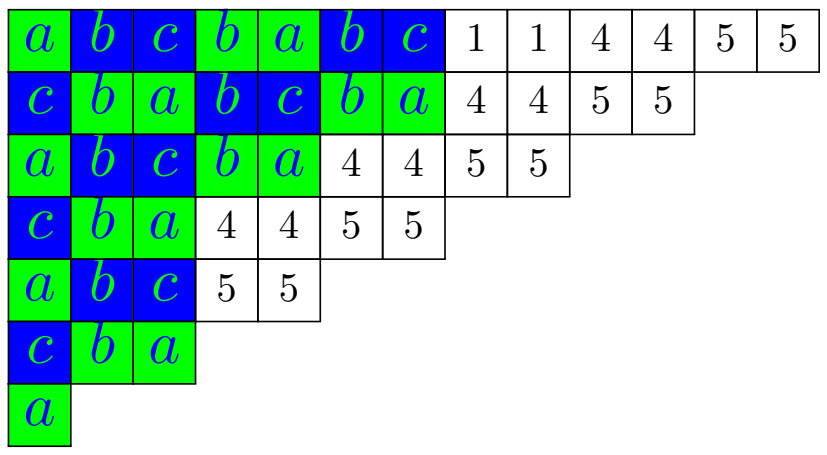

\[
\scalebox{0.7}{
\begin{ytableau}
*(green) \scalebox{1.5}{\color{blue}$a$} & *(blue) \scalebox{1.5}{\color{green}$b$} & *(blue) \scalebox{1.5}{\color{green}$c$} & *(green) \scalebox{1.5}{\color{blue}$b$} & *(green) \scalebox{1.5}{\color{blue}$a$} & *(blue) \scalebox{1.5}{\color{green}$b$} & *(blue) \scalebox{1.5}{\color{green}$c$} & 1 & 1 & 4 & 4 & 5 & 5 \\
*(blue) \scalebox{1.5}{\color{green}$c$} & *(green) \scalebox{1.5}{\color{blue}$b$} & *(green) \scalebox{1.5}{\color{blue}$a$} & *(blue) \scalebox{1.5}{\color{green}$b$} & *(blue) \scalebox{1.5}{\color{green}$c$} & *(green) \scalebox{1.5}{\color{blue}$b$} & *(green) \scalebox{1.5}{\color{blue}$a$} & 4 & 4 & 5 & 5\\
*(green) \scalebox{1.5}{\color{blue}$a$} & *(blue) \scalebox{1.5}{\color{green}$b$} & *(blue) \scalebox{1.5}{\color{green}$c$} & *(green) \scalebox{1.5}{\color{blue}$b$} & *(green) \scalebox{1.5}{\color{blue}$a$} & 4 & 4 & 5 & 5\\
*(blue) \scalebox{1.5}{\color{green}$c$} & *(green) \scalebox{1.5}{\color{blue}$b$} & *(green) \scalebox{1.5}{\color{blue}$a$} & 4 & 4 & 5 & 5\\
*(green) \scalebox{1.5}{\color{blue}$a$} & *(blue) \scalebox{1.5}{\color{green}$b$} & *(blue) \scalebox{1.5}{\color{green}$c$} & 5 & 5\\
*(blue) \scalebox{1.5}{\color{green}$c$} & *(green) \scalebox{1.5}{\color{blue}$b$} & *(green) \scalebox{1.5}{\color{blue}$a$}\\
*(green) \scalebox{1.5}{\color{blue}$a$}
\end{ytableau}}
\]
\caption{Obtaining $\pi=1_a3_a3_c3_a5_a7_a7_c$ from the odd staircase}
\label{fig:P2b_pivot_wt}
\end{figure}
\end{example}

\begin{lemma}\label{lemma:level2decomp_step4}
  The size generating function of starting points corresponding to a
  pivot \(\pi\in\mathrm{PIV}(n,k)\) with \(n\) letters \(a\) and \(k\)
  letters \(c\) is given by
  \[
    q^{|\pi|}\,(q^2;q^2)_{n+k}^{-1}.
  \]
\end{lemma}
\begin{proof}
  In the proof of~\cite[Theorem~1.11]{DJ25}, we defined a bijection between the set of grounded partitions \(\lambda\in\mathcal{P}_{2,b}\)
  and the set of pairs \((s,\lambda_e)\) where \(s\) is a sequence of
  labels \(a\) and \(c\) of length \(\ell\) and \(\lambda_e\) is a
  partition into even parts, where the largest part is at most
  \(2\ell\). Moreover, the subsequence of labels restricted to the
  letters \(a\) and \(c\) in \(\lambda\) is \(s\). The size of \(\lambda\) is \(|\pi_s|+|\lambda_e|\), where \(\pi_s\) is the grounded partition with sequence of labels \(s\).

  Fix \(s\) to be the sequence of labels of the pivot \(\pi\). This
  yields a bijection \(\mu\mapsto\lambda_e\) between grounded
  partitions \(\mu\) whose subsequence of labels restricted to the
  letters \(a\) and \(c\) is \(s\) and partitions into even
  parts, where the largest part is at most \(2(n+k)\), with \(|\mu|=|\pi|+|\lambda_e|\). By Lemma~\ref{lemma:level2decomp_step2}, these grounded partitions \(\mu\) are exactly the starting points corresponding to \(\pi\). The size generating function for
  partitions into even parts, where the largest part is at most \(2(n+k)\), is
  \((q^2;q^2)_{n+k}^{-1}\).
\end{proof}

Now we are ready to prove the identity stated in~Theorem~\ref{thm:decomp2}.

\begin{theorem} We have
\begin{align*}
\frac{(-q;q^2)_\infty}{(q;q^2)_\infty} = \sum_{\substack{n\geq k\geq 0}}q^{2k+(n-k)^2}\,Y_{q^4}(n,k)\,(q^2;q^2)_{n+k}^{-1}\,[2(n-k)+2]_q.
\end{align*}
\end{theorem}
\begin{proof}
  Using Lemma~\ref{lem:wt_pivot_2} and
  Proposition~\ref{prop:rev_maj}~(ii) and~(iii) with \(q\) replaced by
  \(q^{-4}\), the size generating function of pivots in
  \(\mathrm{PIV}(n,k)\) can be written as
  \[
    \sum_{\pi \in \mathrm{PIV}(n,k)} q^{|\pi|}
    = q^{(n+k)^2}\cdot q^{-4k}\,Y_{q^{-4}}(n,k)\cdot q^{2k} 
    = q^{(n+k)^2-2k}\,Y_{q^{-4}}(n,k).
  \]
  The minimal possible size in \(\mathrm{PIV}(n,k)\) is
  \(2k + (n-k)^2\), realised by the grounded partition
  \[
    1_a1_c1_a1_c\cdots1_a3_a5_a\cdots(2(n-k)-1)_a
  \]
  when \(n>k\), and by \(1_a1_c\cdots1_a1_c\) when \(n=k\). Moreover,
  the \( q \)-Yamanouchi number \(Y_{q^4}(n,k)\) is a palindromic
  polynomial of degree \(4k(n-1)\) with constant term \(1\).
  Therefore, we obtain
  \begin{align}\label{eq:P2_pivot_genfn}
    \sum_{\pi \in \mathrm{PIV}(n,k)} q^{|\pi|}
    = q^{2k+(n-k)^2}\,Y_{q^4}(n,k).
  \end{align}
  By Lemma~\ref{lemma:level2decomp_step4}, we multiply
  (\ref{eq:P2_pivot_genfn}) by the factor \((q^2;q^2)_{n+k}^{-1}\), to
  obtain the generating function of the starting points corresponding
  to pivots in \(\mathrm{PIV}(n,k)\). Finally, by
  Lemmas~\ref{lemma:level2decomp_step1}
  and~\ref{lemma:level2decomp_step2}, if the pivot of a starting point
  lies in \(\mathrm{PIV}(n,k)\), then the length of the corresponding
  blue string is \(2(n-k)+1\). Each arrow increases the size by one.
  Hence we multiply the generating function by \([2(n-k)+2]_q\) and
  sum over \(n\geq k\geq 0\), to obtain the desired identity.
\end{proof}

\begin{remark}
In the affine crystal graph $\mathcal{B}_{2,b}^{\mathcal{P}}$, erasing the blue arrows instead of the green arrows yields the same restriction formula, unlike in the case of $\mathcal{B}_{1,b}^{\mathcal{P}}$, discussed in Subsection~\ref{subsection:decomp1}. This is clear from a representation-theoretic point of view: the affine crystal graph $\mathcal{B}_{2,b}^{\mathcal{P}}\cong\mathcal{B}_{2,b}^{\mathcal{A}}$ has highest weight $\Lambda_0+\Lambda_1$, which is preserved by the involution induced by swapping $\Lambda_0$ and $\Lambda_1$.
\end{remark}

% Restore original \clearpage
\let\clearpage\oldclearpage

% Bibliography
\printbibliography

@article{AK94,
  author  = {Ariki, S. and Koike, K.},
  title   = {A Hecke algebra of $(\mathbb{Z}/r\mathbb{Z})\wr S_n$ and construction of its irreducible representations},
  journal = {Advances in Mathematics},
  volume  = {106},
  number  = {2},
  year    = {1994},
  pages   = {216--243}
}

@article{KMOTU00,
  author        = {Kuniba, Atsuo and Misra, Kailash C. and Okado, Masato and Takagi, Taichiro and Uchiyama, Jun},
  title         = {Paths, {D}emazure crystals, and symmetric functions},
  journal       = {Journal of Mathematical Physics},
  volume        = {41},
  number        = {9},
  pages         = {6477--6486},
  year          = {2000},
  doi           = {10.1063/1.1286284},
  eprint        = {q-alg/9612018},
  archivePrefix = {arXiv},
  note          = {Theorem numbers refer to the arXiv version}
}

@incollection{Kas95,
  author    = {Kashiwara, M.},
  title     = {On crystal bases},
  booktitle = {Representations of Groups ({B}anff, {AB}, 1994)},
  series    = {CMS Conference Proceedings},
  volume    = {16},
  publisher = {American Mathematical Society},
  address   = {Providence, RI},
  year      = {1995},
  pages     = {155--197}
}

@article{Kas90,
  author  = {Kashiwara, M.},
  title   = {Crystalizing the $q$-analogue of universal enveloping algebras},
  journal = {Communications in Mathematical Physics},
  volume  = {133},
  number  = {2},
  pages   = {249--260},
  year    = {1990}
}

@article{Kas91,
  author  = {Kashiwara, M.},
  title   = {On crystal bases of the $q$-analogue of universal enveloping algebras},
  journal = {Duke Mathematical Journal},
  volume  = {63},
  number  = {2},
  pages   = {465--516},
  year    = {1991}
}

@article{Kl95,
  author  = {Kleshchev, A. S.},
  title   = {Branching rules for modular representations of symmetric groups, {II}},
  journal = {J. reine angew. Math.},
  volume  = {459},
  pages   = {163--212},
  year    = {1995}
}

@article{GK76,
  author  = {Greene, C. and Kleitman, D. J.},
  title   = {Strong versions of {S}perner's theorem},
  journal = {J. Combin. Theory Ser. A},
  volume  = {20},
  pages   = {80--88},
  year    = {1976}
}

@book{Mac95,
  author    = {Macdonald, I. G.},
  title     = {Symmetric Functions and {H}all Polynomials},
  edition   = {Second},
  series    = {Oxford Mathematical Monographs},
  publisher = {Oxford University Press},
  address   = {Oxford},
  year      = {1995}
}

@article{You28,
  author  = {Young, A.},
  title   = {On quantitative substitutional analysis (third paper)},
  journal = {Proceedings of the London Mathematical Society (2)},
  volume  = {28},
  pages   = {255--292},
  year    = {1928},
  doi     = {10.1112/plms/s2-28.1.255}
}

@article{Mac13,
  author  = {MacMahon, P. A.},
  title   = {The indices of permutations and the derivation therefrom of functions of a single variable associated with the permutations of any assemblage of objects},
  journal = {American Journal of Mathematics},
  volume  = {35},
  number  = {3},
  year    = {1913},
  pages   = {281--322},
  doi     = {10.2307/2370312}
}

@incollection{Kra15,
  author    = {Krattenthaler, C.},
  title     = {Lattice Path Enumeration},
  booktitle = {Handbook of Enumerative Combinatorics},
  editor    = {B{\'o}na, M.},
  series    = {Discrete Mathematics and Its Applications},
  publisher = {Chapman \& Hall / CRC Press},
  year      = {2015},
  pages     = {589--678}
}

@article{Whi79,
  author  = {Whitworth, W. A.},
  title   = {Arrangements of $m$ things of one sort and $n$ things of another sort, under certain conditions of priority},
  journal = {Messenger of Mathematics},
  volume  = {8},
  pages   = {105--114},
  year    = {1879}
}

@article{AKT08,
  author       = {Ariki, S. and Kreiman, V. and Tsuchioka, S.},
  title        = {On the tensor product of two basic representations of $U_v(\widehat{\mathfrak{sl}}_e)$},
  journal      = {Advances in Mathematics},
  volume       = {218},
  number       = {1},
  pages        = {28--86},
  year         = {2008},
  doi          = {10.1016/j.aim.2007.11.018}
}

@article{AM00,
  author  = {Ariki, S. and Mathas, A.},
  title   = {The number of simple modules of the Hecke algebras of type $G(r,1,n)$},
  journal = {Mathematische Zeitschrift},
  volume  = {233},
  number  = {3},
  year    = {2000},
  pages   = {601--623},
  doi     = {10.1007/s002090050489}
}

@article{Ari96,
  author  = {Ariki, S.},
  title   = {On the decomposition numbers of the Hecke algebra of $G(m,1,n)$},
  journal = {Journal of Mathematics of Kyoto University},
  volume  = {36},
  number  = {4},
  year    = {1996},
  pages   = {789--808},
  doi     = {10.1215/kjm/1250518452}
}

@article{Br98,
  author  = {Brundan, J.},
  title   = {Modular branching rules and the Mullineux map for Hecke algebras of type A},
  journal = {Proceedings of the London Mathematical Society},
  series  = {3},
  volume  = {77},
  number  = {3},
  year    = {1998},
  pages   = {551--581},
  doi     = {10.1112/S0024611598000562}
}

@book{BS17,
  author    = {Bump, D. and Schilling, A.},
  title     = {Crystal Bases: Representations and Combinatorics},
  publisher = {World Scientific Publishing Company},
  year      = {2017},
  isbn      = {9814733466},
}

@article{CLSXY24,
  author       = {Chern, S. and Li, Z. and Stanton, D. and Xue, T. and Yee, A. J.},
  title        = {The Ariki–Koike algebras and Rogers–Ramanujan type partitions},
  journal      = {Journal of Algebraic Combinatorics},
  volume       = {60},
  pages        = {491--540},
  year         = {2024},
  doi          = {10.1007/s10801-024-01340-z}
}

@article{DHK25,
  title       = {{Partition identities from higher level crystals of $A_1^{(1)}$}},
  author      = {Dousse, J. and Hardiman, L. and Konan, I.},
  url         = {https://hal.science/hal-03448464},
  journal     = {{Proceedings of the American Mathematical Society}},
  publisher   = {{American Mathematical Society}},
  volume      = {153},
  number      = {4},
  pages       = {1363-1382},
  year        = {2025},
  month       = Feb,
  doi         = {10.1090/proc/16417},
  hal_id      = {hal-03448464},
  hal_version = {v1}
}

@misc{DJ25,
  author        = {Dombos, B. and Jang, J.},
  title         = {{$A_1^{(1)}$-Grounded partitions at levels $1$ and $2$, Part I: bijections}},
  year          = {2025},
  eprint        = {2508.02664},
  archivePrefix = {arXiv},
  primaryClass  = {math.CO},
  doi           = {10.48550/arXiv.2508.02664}
}

@misc{DK19b,
      title={Generalisations of Capparelli's and Primc's identities, II: perfect $A_{n-1}^{(1)}$ crystals and explicit character formulas}, 
      author={Dousse, J. and Konan, I.},
      year={2020},
      eprint={1911.13189},
      archivePrefix={arXiv},
      primaryClass={math.QA},
      url={https://arxiv.org/abs/1911.13189}, 
}

@article{FLOTW99,
  author  = {Foda, O. and Leclerc, B. and Okado, M. and Thibon, J.-Y. and Welsh, T. A.},
  title   = {Branching functions of $A_{n-1}^{(1)}$ and {J}antzen--{S}eitz problem for {A}riki--{K}oike algebras},
  journal = {Advances in Mathematics},
  volume  = {141},
  number  = {2},
  year    = {1999},
  pages   = {322--365},
  doi     = {10.1006/aima.1998.1783}
}

@book{HK02,
  author    = {Hong, J. and Kang, S.-J.},
  title     = {Introduction to Quantum Groups and Crystal Bases},
  publisher = {Amer. Math. Soc.},
  series    = {Graduate Studies in Mathematics},
  volume    = {42},
  year      = {2002}
}

@article{JMMO91,
  author    = {Jimbo, M. and Misra, K. and Miwa, T. and Okado, M.},
  title     = {Combinatorics of representations of $U_q(\widehat{\mathfrak{sl}}_n)$ at $q=0$},
  journal   = {Communications in Mathematical Physics},
  volume    = {136},
  number    = {3},
  pages     = {543--566},
  year      = {1991},
  doi       = {10.1007/BF02100899},
  publisher = {Springer}
}

@article{KMN2a,
author = {Kang, S. and Kashiwara, M. and Misra, K. and Miwa, T. and Nakashima, T. and Nakayashiki, A.},
year = {1992},
month = {04},
pages = {449--484},
title = {Affine Crystals and Vertex Models},
volume = {7},
journal = {Int. Journ. Mod. Phys. A},
doi = {10.1142/S0217751X92003896}
}

@article{KMN92,
  author  = {Kang, S. and Kashiwara, M. and Misra, K. and Miwa, T. and Nakashima, T. and Nakayashiki, A.},
  title   = {Perfect crystals of quantum affine Lie algebras},
  journal = {Duke Mathematical Journal},
  volume  = {68},
  number  = {3},
  pages   = {499--607},
  year    = {1992},
  doi     = {10.1215/S0012-7094-92-06821-9}
}

@article{KN94,
  author  = {Kashiwara, M. and Nakashima, T.},
  title   = {Crystal graphs for representations of the $q$-analogue of classical Lie algebras},
  journal = {Journal of Algebra},
  volume  = {165},
  number  = {2},
  pages   = {295--345},
  year    = {1994},
  doi     = {10.1006/jabr.1994.1114}
}

@article{Kra89,
  author  = {Krattenthaler, C.},
  title   = {Counting lattice paths with a linear boundary I},
  journal = {Sitz.ber. d. ÖAW Math.-naturwiss. Klasse},
  volume  = {198},
  year    = {1989},
  pages   = {87-107}
}

@article{Kra89b,
  author  = {Krattenthaler, C.},
  title   = {Counting lattice paths with a linear boundary {II}},
  journal = {Sitz.ber. d. ÖAW Math.-naturwiss. Klasse},
  volume  = {198},
  year    = {1989},
  pages   = {171--199}
}

@book{Lep78,
  author    = {Lepowsky, J.},
  title     = {Lectures on Kac--Moody Lie Algebras},
  year      = {1978},
  publisher = {Université Paris VI},
  address   = {Paris},
  note      = {Spring 1978}
}

@incollection{Lepowskyproduct,
  author    = {Lepowsky, J.},
  title     = {Affine {L}ie algebras and combinatorial identities},
  booktitle = {Lie Algebras and Related Topics ({N}ew {B}runswick, {NJ}, 1981)},
  editor    = {Winter, D.},
  series    = {Lecture Notes in Mathematics},
  volume    = {933},
  publisher = {Springer},
  address   = {Berlin},
  year      = {1982},
  pages     = {130--156}
}

@article{LLT96,
  author  = {Lascoux, A. and Leclerc, B. and Thibon, J-Y.},
  title   = {Hecke algebras at roots of unity and crystal bases of quantum affine algebras},
  journal = {Communications in Mathematical Physics},
  volume  = {181},
  number  = {1},
  year    = {1996},
  pages   = {205--263},
  doi     = {10.1007/BF02101678}
}

@article{LS78,
  author  = {Lascoux, A. and Sch{\"u}tzenberger, M.-P.},
  title   = {Sur une conjecture de {H}. {O}. {F}oulkes},
  journal = {C. R. Acad. Sci. Paris S{\'e}r. A-B},
  volume  = {286},
  pages   = {323--324},
  year    = {1978}
}

@incollection{LS81,
  author    = {Lascoux, A. and Sch\"utzenberger, M.-P.},
  title     = {Le mono\"ide plaxique},
  booktitle = {Noncommutative structures in algebra and geometric combinatorics (Naples, 1978)},
  series    = {Quad. ``Ricerca Sci.''},
  volume    = {109},
  pages     = {129--156},
  publisher = {CNR},
  address   = {Rome},
  year      = {1981}
}

@article{MM90,
  author    = {Misra, K. and Miwa, T.},
  title     = {Crystal base for the basic representation of $U_q(\widehat{\mathfrak{sl}}_n)$},
  journal   = {Communications in Mathematical Physics},
  volume    = {134},
  number    = {1},
  pages     = {79--88},
  year      = {1990},
  doi       = {10.1007/BF02102090},
  publisher = {Springer}
}

@article{NY97,
  author  = {Nakayashiki, A. and Yamada, Y.},
  title   = {Kostka polynomials and energy functions in solvable lattice models},
  journal = {Selecta Mathematica. New Series},
  year    = {1997},
  volume  = {3},
  pages   = {547--599},
  month   = dec,
  doi     = {10.1007/s000290050020}
}

@book{Sta1,
  author    = {Stanley, R. P.},
  title     = {Enumerative Combinatorics, Volume 2},
  series    = {Cambridge Studies in Advanced Mathematics},
  volume    = {62},
  publisher = {Cambridge University Press},
  address   = {Cambridge, UK},
  year      = {1999},
  isbn      = {0-521-56069-1},
  doi       = {10.1017/CBO9780511609589}
}

@book{Wak01,
  author    = {M. Wakimoto},
  title     = {Lectures on Infinite-dimensional Lie Algebra},
  publisher = {World Scientific Publishing Co. Pte. Ltd.},
  year      = {2001},
  address   = {Singapore},
  isbn      = {9814494003; 9789814494007}
}

\end{document}